%% file: dde-highorder.tex
\documentclass[english]{article}
\usepackage[T1]{fontenc}		
\usepackage[utf8]{inputenc}	
\usepackage{hyperref} 	
\usepackage{textcomp}		
\usepackage{amsmath}		
\usepackage{graphicx}   	
\usepackage{amssymb}		
\usepackage[skip=0pt,font=small,labelfont=bf,margin=24pt]{caption}	
\usepackage{url}			
\usepackage{geometry} 	
\usepackage{enumerate}   
\usepackage{fancyvrb}    
\usepackage{babel} 		
\usepackage[title]{appendix}

\graphicspath{{./figures/pdf/}} 

\newtheorem{theorem}{Theorem}
\newtheorem{definition}{Definition}
\newtheorem{lemma}[theorem]{Lemma}
\newtheorem{rem}[theorem]{Remark}
\newtheorem{pro}[theorem]{Proposition}

\newcommand{\Z}{\ensuremath{\mathbb{Z}}}
\newcommand{\N}{\ensuremath{\mathbb{N}}}
\newcommand{\R}{\ensuremath{\mathbb{R}}}
\newcommand{\I}{\ensuremath{\mathbb{I}}}
\newcommand{\X}{\ensuremath{\mathcal{X}}}
\newcommand{\Y}{\ensuremath{\mathcal{Y}}}
\newcommand{\epsi}{\ensuremath{\varepsilon}}
\newcommand{\Rem}{\ensuremath{\mathtt{Rem}}}
\newcommand{\roughEncl}{\ensuremath{\mathtt{roughEncl}}}
\newcommand{\hull}{\ensuremath{\mathtt{hull}}}
\newcommand{\midpoint}{\ensuremath{\mathtt{m}}}

\DeclareMathOperator{\bd}{\ensuremath{\partial}}
\DeclareMathOperator{\closure}{cl}
\DeclareMathOperator{\diam}{\mathbf{diam}}
\DeclareMathOperator{\domain}{\ensuremath{\mathtt{Dom}}}
\DeclareMathOperator{\interior}{int}
\newcommand{\cover}[1]{\stackrel{#1}{\Longrightarrow}}
\DeclareMathOperator{\Ball}{\ensuremath{\mathbf{B}}}

\newcommand{\eval}{\ensuremath{\mathtt{E}}}
\newcommand{\taylor}{\ensuremath{\mathtt{T}}}
\newcommand{\summa}{\ensuremath{\mathtt{S}}}
\newcommand{\ixi}{\ensuremath{[\xi]}}
\newcommand{\matrices}{\ensuremath{\mathcal{M}}}
\newcommand{\dimension}{\ensuremath{\mathtt{dim}}}

\def\qed{{\hfill{\vrule height5pt width5pt depth0pt}\medskip}}

\begin{document}
\begin{center}
{\bf \LARGE
High-order Lohner-type algorithm for 
rigorous computation of Poincar\'e maps in 
systems of Delay Differential 
Equations with several delays}

\vskip 0.5cm

{\large Robert Szczelina$^{1, 2}$,  Piotr Zgliczy\'nski$^{1, 3}$}

\vskip 0.25cm

\today

\vskip 0.25cm
\end{center}
{\small
\begin{tabular}{ll}
1. & Jagiellonian University, Faculty of Mathematics and Computer Science,  \\
   & Łojasiewicza 6, 30-348 Kraków, Poland \\
2. & corresponding author: robert.szczelina@uj.edu.pl \\
3. & umzglicz@cyf-kr.edu.pl
\end{tabular}
}

\begin{abstract}
We present a Lohner-type algorithm for rigorous
integration of systems of Delay Differential Equations (DDEs) with multiple delays,
and its application in computation of Poincar\'e maps, to study the dynamics
of some bounded, eternal solutions. The algorithm is based on a piecewise 
Taylor representation of the solutions in the phase-space and it exploits the 
smoothing of solutions occurring in DDEs to produces enclosures of solutions of a high order.
We apply the topological techniques to prove various 
kinds of dynamical behaviour, for example, existence of (apparently) unstable periodic 
orbits in Mackey-Glass Equation (in the regime of parameters where chaos is 
numerically observed) and persistence of symbolic dynamics 
in a delay-perturbed chaotic ODE (the R\"ossler system). 
\end{abstract}

\vskip 0.25cm

\input sections/1_introduction.tex

\input sections/2_phasespace.tex

\input sections/3_integrator.tex
\input sections/4_covrel.tex
\input sections/5_applications.tex

\section{Acknowledgements}
Research has been supported by Polish
National Science Centre grant no. 2016/22/A/ST1/00077.

\bibliographystyle{plain} 	
\bibliography{dde}

\begin{appendices}
\input sections/6_appendix_A_lohner.tex

\input sections/7_appendix_B_data.tex
\end{appendices}

\end{document}

%% file: sections/1_introduction.tex
\section{\label{sec:intro}Introduction}

We consider a system of Delay Differential Equations (DDEs) with constant delays
and the initial condition of the following form:
\begin{equation}
  \begin{cases}
	x'(t) = f\left(x(t), x(t - \tau_1), x(t - \tau_2), \ldots, x(t - \tau_m)\right), & t \ge 0\\
	x(t)  = \psi(t), & t \in [-\tau, 0],
  \end{cases} \label{eq:dde-many}
\end{equation}
where $m \in \N$ and $\tau = \tau_1 > \tau_2 > \ldots > \tau_m \ge 0$ are
the delays, $x'$ is understood as a right derivative
and $\psi : [-\tau, 0] \to \R^d$ is of class $C^0$ on $[-\tau, 0]$,
$x(t) \in \R^d$, $f : \R^{(m+1)d} \to \R^d$.

In \cite{nasza-praca-focm} we have presented a method of producing rigorous
estimates on the function $x(t)$ for $t \ge 0$ for the simplest scalar ($d = 1$) 
DDE with a single delay:
\begin{equation}
x'(t) = f(x(t), x(t-\tau)) \\
\label{eq:dde}
\end{equation}
The algorithm presented in \cite{nasza-praca-focm} is an explicit Taylor
method with piecewise Taylor representation of the solution over a fixed step size grid
and with Lohner-type control of the wrapping effect encountered in
interval arithmetics \cite{capd-lohner}. The method consists
of two algorithms: one computing Taylor coefficients of the solutions 
at equally spaced grid points (together with the rigorous estimate on the error size), 
and the second one to compute enclosures of the solution
segment after an arbitrary step of size $\epsi$, smaller than the grid step size $h$. 
The method is suited to construct Poincar\'e maps in the phase space of the
DDEs and it was successfully
applied to prove several (apparently stable) periodic solutions to scalar DDEs 
\cite{nasza-praca-focm,szczelina-ejde}
(among them to Mackey-Glass equation). However, the second method - $\epsi$ step part - 
is not optimal in the sense of the local error order. 
Essentially, the local error of some of the 
coefficients in the Taylor representation of solution is $O(h)$.
The reason is that some of the coefficients are computed using 
just an explicit Euler method with very rough estimates on the derivative. 
With this apparent loss of
accuracy, the images of Poincar\'e maps 
in \cite{nasza-praca-focm} 
are computed with less than optimal quality, and are not well
suited to handle more diverse spectrum of dynamical results.
 
In this work, we provide effective way to decrease the local order
of the full-step algorithm after each full delay
of the integration procedure to significantly reduce 
the error size later on, when applying the 
second $\epsi$ step procedure. 
Under some additional but reasonable assumptions about the integration time
being long enough (see Definition~\ref{def:long-enough-time} and Section~\ref{sec:epsi-steps}), 
the modification allows to decrease the local error size of the $\epsi$ step method for 
all coefficients to $O(h^{n+1})$ where $n$ is the order of representation 
of the initial function and $h$ the step size of interpolation grid 
(compared to $O(h)$ for the previous version from \cite{nasza-praca-focm}).
 
All those enhancements are done without a significant increase in computational
complexity of the most time consuming part of the algorithm: 
a Lohner-type method of controlling the wrapping effect.
What is more, we present an elegant and more general Lohner-type
method of wrapping effect control to handle both systems of equations 
and many delays such as Eq.~\eqref{eq:dde-many}. 
We also employ more elaborate Lohner sets 
to further reduce undesirable effects of interval arithmetic.
With all those improvements, the method produces estimates
on solutions of several orders of magnitude better than the
previous one. 
 
As a presentation of effectiveness of the new method, we give proofs of the
existence of periodic solutions to Mackey-Glass equation for a wider spectrum of
parameters than in \cite{nasza-praca-focm}. The proofs are done 
for parameters in the chaotic regime and the orbits are apparently unstable. 
To this end, we need to expand on the theory, so we extend the concept of covering 
relations \cite{covrel-GiZ-1} to infinite Banach spaces, and we use Fixed Point Index in
Absolute Neighbourhood Retracts (ANRs) \cite{granas-book} to prove a 
Theorem~\ref{thm:symbolic-dynamics}
about the existence of orbits for compact mappings in infinite dimensional spaces 
following chains of covering
relations. We use this technique to show existence of symbolic dynamics in a perturbed 
model $x'(t) = f(x(t)) + \epsi \cdot g(x(t-\tau))$, where $f$ is a chaotic ODE in
three dimensions (R\"ossler system) and for a couple of $g$'s which are some 
explicitly bounded functions. We hope similar techniques will allow to
prove chaos in Mackey-Glass equation \cite{mackey-glass}.

The paper is organized as follows: in Section~\ref{sec:phasespace} we present
some basic theory for DDEs with constant delays, and we recall shortly the basic
structure of (p,n)-functions sets to represent objects in the phase-space of
\eqref{eq:dde-many}. We also generalize this structure and we discuss its properties.
In Section~\ref{sec:integrator} we recall algorithm from \cite{nasza-praca-focm}
within a new, more general notation and we introduce several modifications that
will be crucial for the complexity and accuracy of the algorithm. This new
algorithm will form a base to some improvements in the construction of
Poincar\'e maps in the phase space, especially to enhance the quality of the
estimates. We present some benchmarks to show how the new estimates are in
comparison to the old algorithm. In Section~\ref{sec:covrel}, we present
topological tools to prove existence of a special kind of solutions to
DDEs~\eqref{eq:dde-many}. We go beyond the Schauder Fixed Point Theorem used in
\cite{nasza-praca-focm}: we use Fixed Point Index on ANRs \cite{granas-book}
and we adapt the notion of covering relations \cite{covrel-GiZ-1} to the setting of
(p,n)-functions sets describing the infinite dimensional phase space of the
DDEs. The compactness of the solution operator in the phase space for times
bigger than delay allows to apply the Schauder Fixed Point index in our case. We
establish theorems to prove existence of symbolic dynamics conjugated to the
sequences of covering relations on (p,n)-functions sets. In
Section~\ref{sec:applications} we apply presented methods to prove existence of
(apparently unstable) periodic solutions to the Mackey-Glass equation, for
the value original value of parameters for which Mackey and Glass observed numerically
chaotic attractor \cite{mackey-glass}. We also prove existence of symbolic 
dynamics in a delay-perturbed chaotic ODE (R\"ossler system).


\subsection{Our results in the perspective of current research in the field}

There are many important works that establish the existence and the shape of the
(global) attractor under some assumptions on \eqref{eq:dde}, for example if it
is of the form $x'=-\mu x(t) + f (x(t - 1))$ and under the assumption that $f$
is strictly monotonic, either positive or negative, or if $f$ has a simple
explicit formula, usually piece-wise constant or affine. We would like here to
point out some results, but the list is for sure not exhaustive (we refer to the
mentioned works and references therein). Mallet-Paret and Sell used discrete
Lyapunov functionals to prove a Poincar\'e-Bendixson type of theorem for special
kind of monotone systems \cite{mallet-paret-sell-bendixon}. Krisztin, Walther
and Wu have conducted an intensive study on systems having a monotone positive
feedback, including studies on the conditions needed to obtain the shape of a
global attractor, see \cite{krisztin-book} and references therein. Krisztin and
Vas proved that in the case of a monotonic positive feedback $f$, under some
assumptions on the stationary solutions, there exists large amplitude slowly
oscillatory periodic solutions (LSOPs) which revolve around more than one
stationary solution \cite{dde-lsop-krisztin}. Vas continued this work and showed
a method to construct $f$ such that the structure of the global attractor may be
arbitrarily complicated (containing an arbitrary number of unstable LSOPs)
\cite{vas-big-attractor}. On the other hand, Lani-Wayda and Walther were able to
construct systems of the form $x'=f(x(t - 1))$ for which they proved the
existence of a dynamic which is conjugate to a symbol shift (Smale's horseshoe)
\cite{lani-wayda-walther-chaos}. Srzednicki and Lani-Wayda proved the existence
of multiple periodic orbits and the existence of chaos for some periodic,
tooth-shaped (piecewise linear) $f$ by the use of the generalized Lefshetz fixed
point theorem \cite{srzednicki-wayda}. A nice review of works that deal with the
question of existence of chaos in Mackey-Glass and similar systems are compiled
in Walther review \cite{walther-chaos-review}. Recently, a new approach have
been used to prove the existence of some periodic orbits to the Mackey-Glass
equation in a limiting case when $n \to \infty$ \cite{krisztin-mackey-glass}.

While impressive, all mentioned analytic/theoretic results are usually hard to
apply in the context of general functions $f$, so we might look for other means
of obtaining rigorous results in such cases, for example, by employing computers
for this task.  
In recent years, there were many computer assisted proofs of
various dynamical properties for maps, ODEs and (dissipative) Partial
Differential Equations ((d)PDEs) by application of the theory of dynamical
systems with estimates obtained from rigorous numerical methods and interval
arithmetic, see for example \cite{capd-article} and references therein. A big
achievement of the rigorous computations are proofs of the existence of chaos
and strange attractors, for example the paper by Tucker
\cite{tucker-lorenz}, and recently to prove chaos in Kuramoto-Shivasinski PDE
\cite{chaos-kuramoto}. The application of rigorous numerical methods to DDEs
started to appear a few years ago and are steadily getting more attention.
Probably the first method used to prove existence of periodic orbits by the
expansion in Fourier modes was given in \cite{zalewski}, and then in a more
general framework and by a different theoretical approach in
\cite{lessard2010,lessard1}. Other methods, strongly using the form of r.h.s.
$f$ in \eqref{eq:dde}, were used in \cite{dde-lsop-krisztin} to prove the
structure of the global attractor; then in \cite{proof-wright-conjecture} to
close a gap in the proof of the Wright conjecture; and finally recently in
\cite{feri-krisztin-mackey-glass} to show the existence of many stable periodic
orbits for a DDE equation that is the limiting case of Mackey-Glass equation
when $n \to \infty$. To the author's knowledge, the results from our work
\cite{nasza-praca-focm} are the first application of rigorous integration
(forward in time) of DDEs in the full phase-space for a general class of
problems to prove the existence of some dynamics, namely the existence of
apparently stable periodic orbits in Mackey-Glass equation. A different approach
to one presented in our work \cite{nasza-praca-focm} was recently published
which uses Chebyshev polynomials to describe solutions in the phase space and a
rigorous fixed point finding argument to produce estimates on the solutions to
DDEs forward in time, together with estimates on the Frech\'et derivative of the
time-shift operator $\varphi(\tau, \cdot)$ \cite{cheb-fun-dde}, however the
presented approach has one disadvantage: it can find solutions only on full
delay intervals, therefore cannot be used directly to construct Poincar\'e maps.
Recently, the extension of those methods was used to prove persistence of periodic solutions
under small perturbations of ODEs \cite{lessard-ode-perturb-dde}, and
a similar approach was used 
in a rigorous method of numerically solving initial value
problems to State-Dependent DDEs \cite{church-sdde-rigorous}. This last work
uses similar technique as our work to subdivide the basic interval into smaller
pieces and piecewise polynomial interpolation of the functions in the phasespace, 
but instead of Taylor it uses Chebyshev polynomials 
and a fixed-point finding argument to prove
existence of a true solution nearby.  
On the other hand, the parametrization method was used to prove the persistence of periodic
orbits in delay-perturbed differential equations, including the state-dependent
delays \cite{delallave-sdde-ode}, however it assumes that $\tau$ is relatively
small. Our method has an advantage over those methods, as it allows for 
a larger amplitude of the perturbation and to prove theorems beyond the 
existence of periodic orbits, as we are showing persistence of symbolic dynamics 
in a perturbed ODE. Finally, there are also some methods to obtain rigorous bounds
on the solutions, e.g. \cite{simple-exp-delay}, however, as authors say, 
they do not produce estimates of quality good enough to prove theorems. 

\subsection{Notation}
For reader's convenience we include here all the basic notions used in this paper.
We will also remind them the first time they are used in the text, if necessary.

We will denote by $C^k([-\tau, 0], \R^d)$
the set of functions which are $C^{k}$ on $(-\tau, 0)$ and right and left
derivatives up to $k$ exist at $t = -\tau$ and $t = 0$, respectively.
For short, we will usually write $C^k$ to denote $C^k([-\tau, 0], \R^d)$
when $d$ and $\tau$ is known from the context.

We use standard convention in DDEs to denote
\emph{the segment} of $x : (a-\tau, b) \to \R^d$ at $t \in (a, b)$
by $x_t$, where $x_t(s) = x(t+s)$ for all $s \in [-\tau, 0]$. Then, we will denote
by $\varphi$ the semiflow generated by DDE~\eqref{eq:dde-many} on the
space $C^0$, $\varphi(t, x_0) := x_t$ for a solution $x$ of Eq.~\eqref{eq:dde-many}
with initial data $x_0$.

The algorithms presented in this paper produce estimates on various quantities,
especially, we often work with sets of values that
only encloses some quantity. Therefore, for convenience, by $\I$ we will denote
the set of all closed intervals $[a, b] : a \le b, a, b \in \R$ and we will
denote sets by capital letters like $X, Y, Z...$ etc., and values by lower case
letters $x, y, z$, etc. Usually, the value $x \in X$ for easier reading, but it
will be always stated explicitly in the text for clarity.

Sometimes, instead of using subscripts $x_i$, we will write projections to 
coordinates as $\pi_i x$ or $\pi_{\mathcal{X}} x$ (projection on some  
subspace $\mathcal{X}$ of some bigger space. This will be applied to
increase readability of formulas. 

Let $Z \subset \R^M$. By $\hull(Z)$ we denote the interval hull of $Z$, that is,
the smallest set $[Z] \in \I^M$ such that $Z \subset [Z]$. By $\overline{Z}$ we
denote the closure of set $Z$, by $\interior Z$ we denote the interior of $Z$
and by $\bd Z$ we denote boundary of $Z$. If $Y$ is some normed vector space and
$Z \subset Y$, then we will write $\bd_Y Z, \interior_Y Z, \closure_Y {Z}$ to
denote boundary, interior and closure of $Z$ in space $Y$.
By $\domain f$ we denote the domain of $f$. 

For multi-index vectors $\eta, \zeta \in \N^p$ we will write $\eta
\ge \zeta$ iff $\eta_i \ge \zeta_i$ for all $i \in \{1,\dots,p\}$.

By $\matrices(k, l)$ we denote the set of matrices of dimensions
$k \times l$ (rows $\times$ columns), while by $Id_{d \times d}$
the identity matrix and by $0_{d \times d}$
the zero matrix in $\matrices(d, d)$. When $d$ is known
from the context we will drop the subscript in $Id$.

By $\Ball^{\|\cdot\|}_{D}(p, r)$ we denote the (open) ball
in $\R^D$ in the given norm $\|\cdot\|$ at a point $p\in R^D$ with
radius $r$. In the case when the norm is known from the context, 
we simply use $\Ball_{D}(p, r)$, and eventually $\Ball_{D}(r)$
for $0$-centered balls.

%% file: sections/2_phasespace.tex
\section{\label{sec:phasespace}Finite dimensional description of the phase space}

In the beginning we will work with Eq.~\ref{eq:dde} (single delay) for 
simplicity of presentation, but all the facts can be applied to a 
more general Eq.~\ref{eq:dde-many}.

As we are interested in \emph{computer assisted proofs} of dynamical phenomena
for \eqref{eq:dde}, we assume that $f$ is a simple/elementary function, so that
it and its derivatives can be given/obtained automatically as \emph{computer
programs (subroutines)}. Many equations encountered in theory and applications
are of this form, two well-known examples that fit into this category are Wright
and Mackey-Glass equations. We will also assume that $f$ is sufficiently smooth,
usually $C^\infty$ in both variables. Under this assumptions, the solution
$x(t)$ of \eqref{eq:dde} with $x_0 = \psi \in C^0$ exists forward in time (for
some maximal time $T_{max}(\psi) \in [0, +\infty]$) and is unique, see e.g.
\cite{dde-ode-book-driver}.

The crucial property of DDEs with $f$ smooth (for simplicity we assume $f \in
C^\infty$) is the \emph{smoothing of solutions} \cite{dde-ode-book-driver}. If
the solution exists for a long enough time, then it is of class at least $C^{k}$
on the interval $(-\tau + \tau \cdot k, \tau \cdot k)$ and it is of class
at least $C^{k}$ at $t = \tau \cdot k$. If $\psi$ is of class $C^{m}$ then $x$ is of
class $C^{m+k}$ on any interval $(-\tau + \tau \cdot k, \tau \cdot k)$.
Moreover, the solutions on the global attractor of \eqref{eq:dde} must be of
class $C^\infty$ (for $f \in C^\infty$). From the topological methods point of
view, the smoothing of solutions implies the semiflow
$\varphi(t, \cdot) : C^0 \to C^0$ is a compact operator for $t \ge \tau$,
essentially by the Arzela-Ascoli Theorem, see e.g. \cite{nasza-praca-focm} (in
general, $\varphi(t, \cdot) : C^k \to C^k$ is well defined and compact in $C^k$
if $t \ge (k+1) \cdot \tau$).

On the other hand, the solution can still be of a lower class, in some cases -
even only of class $C^0$ (at $t = 0$). It happens due to the very nature of the
DDE~\eqref{eq:dde}, as the right derivative at $t = 0$ is given by
\eqref{eq:dde} whereas the left derivative of the initial data $\psi$ at $0$ can
be arbitrary. This discontinuity propagates in time so the solution $x$, in
general, is only of class $C^k$ at $t = k \cdot \tau$. In other words, a
solution to DDE with an initial segment of higher regularity can sometimes
,,visit'' the lower regularity subset of the phase-space. This behaviour
introduces some difficulties in the treatment of the solutions of DDEs and the
phase-space, especially when one is interested in finding $\varphi(t, x)$ for
$t \ne m \cdot \tau$, $m \in \N$.

In the rest of this section we will recall the notion of (p,n)-functions sets
from \cite{nasza-praca-focm} used in our method to represent functions in the
phase space of DDE~\eqref{eq:dde}. However, we use a slightly different notation
and we introduce some generalizations that will be suitable for the new
integration algorithm in Section~\ref{sec:integrator}.

\subsection{Basic definitions}

The algorithm we are going to discuss in Section~\ref{sec:integrator} is a
modified version of the \emph{(explicit) Taylor rigorous method} for ODEs, that
is, we will be able to produce the Taylor coefficients of the solution at given
times using only the well known recurrent relation resulting from the
successively differentiating formula \eqref{eq:dde} w.r.t. $t$. For this
recurrent formula (presented later in the text, in
Eq.~\eqref{eq:jet-reccurence}) it is convenient to use the language of jets.

Let $m \in \N$ and let $g : \R^m \to \R$ be of class $C^{n}$ and
$z \in \R^d$. We denote by $\alpha$ the $m$-dimensional multi
index $\alpha = \left( \alpha_1, \dots, \alpha_m \right) \in \N^m$ and
we denote $z^{\alpha} = \Pi_{i=1}^{m} z_i^{\alpha_i}$,
$|\alpha| = \sum_{i=1}^{m} \alpha_i$,
$\alpha! = \Pi_{i=1}^{n} \alpha_{i}!$, and
\begin{equation*}
g^{(\alpha)} = \frac{\partial^{|\alpha|}g}{{{\partial z_1^{\alpha_1}}\dots{\partial z_m^{\alpha_m}}}}.
\end{equation*}
By $J^{[n]}_{z}{g}$ we denote the $d$-dimensional
jet of order $n$ of $g$ at $z$, i.e.:
\begin{equation}
\left(J^{[n]}_z{g}\right)(y) = \sum_{|\alpha| \le n} \frac{g^{(\alpha)}(z)}{\alpha!} \cdot (y - z)^\alpha.
\label{eq:jet-general}
\end{equation}
We will identify $J^{[n]}_{z}{g}$ with the collection of the
\emph{Taylor coefficients} $J^{[n]}_{z}{g} = \left(g^{[\alpha]}(z)\right)_{|\alpha| \le n}$,
where
\begin{equation*}
g^{[\alpha]}(z) := \frac{g^{(\alpha)}(z)}{\alpha!}.
\end{equation*}
We will use $J^{[n]}_{z}{g}$ either as a function defined by
\eqref{eq:jet-general} or a collection of numbers depending on the context.
For a function $g : \R^m \to \R^d$ the jet
$J^{[n]}_z(g) = \left(J^{[n]}_z{g_1}, \ldots, J^{[n]}_z{g_d}\right)$
is a collection of jets of components of $g$.

In the sequel we will use extensively the following properties of jets:
\begin{pro}
\label{pro:jet-algebra}
The following are true:
\begin{enumerate}
\item
	if $g : \R \to \R$ then
	$J^{[k]}_z\left({J^{[n]}_z{g}}\right) = J^{[k]}_z{g}$ for $k \le n$;
\item
	if $f = g \circ h : \R \to \R$ for $g : \R^d \to \R$
	and $h : \R \to \R^d$, then
	\begin{equation}
	\label{eq:jet-composition}
		J^{[n]}_{t_0} f =
			J^{[n]}_{t_0} \left(
				\left(J^{[n]}_{h(t_0)} g\right)
				\circ
				\left( J^{[n]}_{t_0}h_1, \dots, J^{[n]}_{t_0} h_d\right)
			\right).
	\end{equation}
\end{enumerate}
\end{pro}
In other words, Equation \eqref{eq:jet-composition} tells us that, in order to
compute $n$-th order jet of the composition, we only need to compose jets
(polynomials) of two functions and ignore terms of order higher than $n$.
For a shorter formulas, we will denote by $\circ_J$ the composition of jets
in \eqref{eq:jet-composition}, i.e. if $a = J^{[n]}_{h(t_0)} g$ and
$b_i = J^{[n]}_{t_0} h_i$, for $i \in \{1,\dots,d\}$ then:
\begin{eqnarray}
\label{eq:jet-composition-short}
a \circ_J b & := & J^{[n]}_{t_0} \left(
		a \circ b
	\right) \\
\notag
	& = &
	J^{[n]}_{t_0} \left(
		\left(J^{[n]}_{b_{[0]}} g\right)
		\circ
		\left( J^{[n]}_{t_0}h_1, \dots, J^{[n]}_{t_0} h_d\right)
	\right).
\end{eqnarray}

\begin{rem}
Operation from Eq.~\eqref{eq:jet-composition} can be effectively implemented
in an algorithmic and effective way by means of Automatic Differentiation \cite{AD1, AD2}.
\end{rem}

From the Taylor's Theorem with integral form of the remainder it follows:
\begin{equation}
\label{eq:taylor-integral-jet}
x(t) = \left(J^{[n]}_{a}x\right)(t) + (n+1) \cdot \int_{a}^{t} x^{[n+1]}(s) \cdot (t-s)^n ds.
\end{equation}
Eq.~\eqref{eq:taylor-integral-jet} motivates the following:
\begin{definition}
\label{def:forward-taylor-representation}
We say that a function $x : \R \to \R$ has
\emph{a forward Taylor representation of order $n$}
on interval $I = [a, a+\delta)$, $\delta > 0$ iff
formula \eqref{eq:taylor-integral-jet} is valid for $x|_I$.

We say that $x : \R \to \R^d$ has a forward Taylor representation on $I$, iff
each component $x_j : \R \to \R$ has the representation on $I$.
\end{definition}
Mostly, we will be using jets to describe (parts of) 
functions $g : I \to \R^d$ with forward Taylor representations, 
therefore, in such cases we understand that in
\begin{equation*}
\left(J^{[n]}_z g\right)(y) = \sum_{k=0}^{n} \frac{g^{(k)}(z)}{k!} \cdot (y - z)^k
\end{equation*}
the $g^{(k)}$ is computed as a right-side derivative. 

It is easy to see and it will be often used in the algorithms: 
\begin{pro}
\label{pro:forward-taylor-derivative}
Assume $x : \R \to \R$ has a forward Taylor representation over
$[t, t+\delta)$ of order $n$. Then for $k \in \{0, \dots, n\}$
the function $x^{[k]} = \frac{x^{(k)}}{k!}$ has a forward
Taylor representation over $[t, t+\delta)$ of order $m = n-k$ and
\begin{align*}
J^{[m]}_{t}(x^{[k]}) & = \left(c^0,\dots,c^{n-k}\right) & \\
(x^{[k]})^{[m+1]}(s) & = \binom{n+1}{k} \cdot  x^{[n+1]}(s)   & s \in [t,t +\delta)
\end{align*}
where
\begin{align*}
c^l & = \binom{l+k}{k} \cdot x^{[l+k]}(t), & l \in 0,\dots,n-k.
\end{align*}
\end{pro}

\begin{pro}
\label{pro:forward-taylor-cmp}
Assume $x : \R \to \R$ has a forward Taylor representation over
$I = [t, t+\delta)$ of order $n$. Then for $k=0,\dots,n$
\begin{align}
\label{eq:kth-tay-form}
x^{[k]}(t + \epsi) & = \sum_{l=0}^{n-k} \binom{l+k}{k} \cdot \left(J^{[n]}_{t} x\right)_{[l+k]} \cdot \epsi^l \ + \\
\nonumber
& + (n+1-k) \cdot \int_{0}^{\epsi} \binom{n+1}{k} \cdot  x^{[n+1]}(t+s) \cdot (\epsi-s)^{n-k} ds
\end{align}
for $\epsi \in [0, \delta)$.
\end{pro}

\begin{rem}[On treating jets as vectors and vice-versa]
As mentioned earlier, for $g : \R \to \R$ the Taylor series $J^{[n]}_{t_0} g$
(which is formally also a function $\R \to \R$) can be uniquely identified
with the collection of the Taylor coefficients 
$\left(g^{[k]}(t_0)\right)_{0 \le k \le n}$, and this collection might 
be identified with a vector in $\R^{n+1}$. One have a freedom how to
organize the sequence into the vector (up to a permutation of coefficients),
but in computer programs we will use the standard ordering
from $k = 0$ at the first coordinate of the vector and 
$k = n$ at the last coordinate. Conversely, for any vector $j \in \R^n$:
\begin{equation}
\label{eq:default-jet-organization}
j = \left(j_{[0]}, j_{[1]}, \ldots, j_{[n]}\right),
\end{equation}
we can build a jet (at some point $t_0$) given by 
\begin{equation}
\label{eq:default-jet-interpretation}
\left(J_{t_0}^{[n]}g\right)(t) = \sum_{k=0}^{n} j_{[k]} (t - t_0)^k.
\end{equation} 
This notion will be convenient when we would have some estimates on the jet,
in particular, we can write that a jet $J^{[n]}_{t_0}{g} \in X \subset \R^{n+1}$,
meaning, that there exists vector $j \in X$ such that \eqref{eq:default-jet-interpretation}
is true for $j$ interpreted as a jet at a given $t_0$.
Also, we can use the convention to do algebraic operations
on jets, such as vector-matrix multiplication to describe jets
in suitable coordinates, etc.  

We will use convention with square brackets $j_{[k]}$ to denote the
relevant coefficient from the sequence $j = J^{[n]}_{z}{g}$, and to underline the
fact that we are using the vector $j$ as its jet interpretation. 

For $g : \R \to \R^d$ the jet $J^{[n]}_{t_0} g$ can be represented as
 a vector in a high dimensional space $\R^M$, where $M = d \cdot (n+1)$.
We organize such jets into vectors in the same manner as in 
Eq.~\eqref{eq:default-jet-organization},
but each $j_{[k]}$ represents $d$ consecutive values.
\end{rem}

\subsection{Outline of the method and the motivation for phase space description}
\label{subsec:integrator-overview}

In a numerical Taylor method for ODEs one produces the jet of solution at the
current time $t_0$ by differentiating the equation $x'(t) = f(x(t))$ w.r.t. $t$ on
both sides at $t_0$, as long as the differentiation makes sense.
For $f \in C^\infty$ we can get any order of the jet at $t_0$ and the situation
is similar in the case of DDE~\eqref{eq:dde}. If $f$ has a jet at
$z = \left(x(t_0), x(t_0 - \tau)\right)$ and $x$ has a jet at $(t_0 - \tau)$,
both of order $n$, then we can proceed as in the case of ODEs to obtain
jet at $t_0$. In the following Lemma we underline the fact that this jet
can be computed from $x(t_0)$ and $J^{[n]}_{t_0 - \tau} x$:
\begin{lemma}
\label{lem:taylor-reccurent-formula}
Let $t_0$ be fixed and $z$ be a solution to \eqref{eq:dde} with $f$
of class at least $C^{n}$. Assume $z$ exists on $[t_0 - \tau, t_0 + \delta]$, and  
$z$ is of class $C^{n}$ on some \emph{past interval} 
$I = [t_0 - \tau, t_0 - \tau + \delta)$ 
for some $\delta > 0$.
Then $z$ is of class $C^{n+1}$ on $I = [t_0, t_0 + \delta)$,  
$J^{[n+1]}_{t_0} z$ exists and it is given explicitly in terms of 
$z(t_0)$, $J^{[n]}_{t_0 - \tau} z$ and r.h.s. $f$ of Eq.~\eqref{eq:dde}.
\end{lemma}
\textbf{Proof: }
The continuity $C^{n+1}$ on $[t_0, t_0 + \delta)$ follows
directly from \eqref{eq:dde}, since $x'$ is of class $C^n$ on $[t_0, t_0 + \delta)$.
Let $F(t) := f(z(t-\tau), z(t))$ and denote the coefficients
of jets $J^{[n]}_{t_0}F$, $J^{[n+1]}_{t_0}z$ and 
$J^{[n]}_{t_0 - \tau}z$ by $F_{[0]},\ldots,F_{[n]}$,
$x_{[0]},\ldots,x_{[n]},x_{[n+1]}$ and $y_{[0]},\ldots,y_{[n]}$ respectively, that is
\begin{eqnarray*}
\left(J^{[n]}_{t_0} F\right) (t) &=& F_{[0]} + F_{[1]} \cdot (t-t_0) + \dots + F_{[n]} \cdot (t-t_0)^n, \\
\left(J^{[n+1]}_{t_0} z\right) (t) &=& x_{[0]} + x_{[1]} \cdot (t-t_0) + \dots + x_{[n]} \cdot (t-t_0)^n + x_{[n+1]} \cdot (t-t_0)^{n+1} \\
\left(J^{[n]}_{t_0 - \tau} z\right) (t) &=& y_{[0]} + y_{[1]} \cdot (t-t_0) + \dots + y_{[n]} \cdot (t-t_0)^n
\end{eqnarray*}
Now Eq.~\eqref{eq:dde} implies that
\begin{equation*}
\left(J^{[n+1]}_{t_0} z\right)' = J^{[n]}_{t_0}{F},
\end{equation*}
or more explicitly:
\begin{equation*}
\left(x_{[0]}+ x_{[1]}(t-t_0)+ \dots + x_{[n+1]}(t-t_0)^{n+1}\right)' = F_{[0]} + F_{[1]}(t-t_0) + \dots + F_{[n]}(t-t_0)^n.
\end{equation*}
Using the obvious fact that $(z^{(k)})' = z^{(k+1)}$,
we have $(z^{[k]})' = (k+1) z^{[k+1]}$ and matching
coefficients of the same powers we end up with:
\begin{equation}
x_{[k]}=\frac{1}{k} F_{[k-1]}.
\label{eq:taylor-jet-x}
\end{equation}
Finally, using Proposition~\ref{pro:jet-algebra} on $J^{[n]}_{t_0}{F}$ we get:
\begin{eqnarray*}
J^{[n]}_{t_0}{F}	& = & \left(J^{[n]}_{(z(t_0), z(t_0-\tau))}{f}\right) \circ_J \left( J^{[n]}_{t_0}{z}, J^{[n]}_{t_0-\tau}{z} \right) \\
                  	& = & \left(J^{[n]}_{(x_{[0]}, y_{[0]})}{f}\right) \circ_J \left( x, y \right).
\end{eqnarray*}

Now, we get the following recurrent formula:
\begin{eqnarray}
F^{[0]}(x_{[0]}, y) & := & f(x_{[0]}, y_{[0]}), \notag \\
F^{[k]}(x_{[0]}, y) & := & \left( \left(J^{[k]}_{(x_{[0]}, y_{[0]})}{f}\right) \circ_J \left( \left(x_{[0]}, w_k * F^{[k-1]}(x_{[0]}, y) \right), \left(y_{[0]},\ldots,y_{[k]}\right) \right) \right),
\label{eq:jet-reccurence}
\end{eqnarray}

for $1 \le k \le n$ with operation $w_n * j$ defined for
a jet $j$ as:
\begin{equation*}
w_n * j := \left(\frac{1}{1} j_{[0]}, \frac{1}{2} j_{[1]}, \ldots, \frac{1}{n} j_{[n-1]}\right).
\end{equation*}
Obviously $F^{[k]}(x_{[0]}, y) = (F_{[0]},\ldots, F_{[k]}) = J^{[k]}_{t_0} F$, and together
with \eqref{eq:taylor-jet-x} we get:
\begin{equation}
\label{eq:jet-reccurence-2}
\left(x_{[0]}, \ldots, x_{[n]}, x_{[n+1]}\right) = \left(x_{[0]}, w_{n+1} \cdot F^{[n]}(x_{[0]}, y)\right),
\end{equation}
that depends only on the formula for $f$, $x_{[0]} = z(t_0)$ and the jet $y = J^{[n]}_{t_0 - \tau} z$.
\qed

We note two important facts. Firstly, the \emph{a priori} existence of the solution $z$
over $[t_0, t_0 + \delta)$
is assumed in Lemma~\ref{lem:taylor-reccurent-formula} and, when doing the integration step, 
it needs to be achieved by some other means - we will later show one way to do that. Secondly,  
Eq.~\eqref{eq:jet-reccurence-2} gives recipe to produce
$J^{[n+1]}_{t_0}$ - a jet of order one higher than the order of the input jet
$y = J^{[n]}_{t_0 - \tau} x$. This simple observation will lead to a significant
improvement to the rigorous integration algorithm in comparison to the first
version presented in \cite{nasza-praca-focm}. To have a complete rigorous method
we will need also formulas to estimate Taylor remainder in
\eqref{eq:taylor-integral-jet} - we will do this later in
Section~\ref{sec:integrator}.

As the jet at $t_0 - \tau$ and the value at $t_0$ allows to compute the jet of
the solution $x$ at $t_0$, the reasonable choice for the description of
functions in the phase-space is to use piecewise Taylor representation of the
solutions at grid points that match the step size of the method. \emph{Uniform step
size over the integration time} will assure that the required jets of the solution
in the formula \eqref{eq:jet-reccurence-2} are always present in the description
of the solution. This approach have been proposed in \cite{nasza-praca-focm} 
with the \emph{uniform order}
of the jets at each grid point. Now, we are going to elaborate how to implement
and use the extra derivative we get in Eq.~\eqref{eq:jet-reccurence-2} to improve
the method. For this, we will need a representation of solutions with \emph{non-uniform
order of jets}.

\subsection{Representation of the phase-space}
\label{subsec:rep-phase-space}

Previously, in \cite{nasza-praca-focm}, we have proposed to describe sets in the
phase space by piecewise Taylor forward representation of a fixed order $n$ on a
uniform grid of points over basic interval $[-\tau, 0]$. Our definition was
stated for $d = 1$ (scalar equations), but the notion can be extended to any
number of dimensions - just by assuming each of the Taylor coefficients in
equations are in fact $d$-dimensional vectors. No formula will be different in
that case. In the rest of the paper we will assume that $d$ is known from the
general context, so we will omit it from the definitions.

We start with a key definition from \cite{nasza-praca-focm} and then
we will propose some generalization that will be relevant to many
important improvements proposed later in this paper.

\begin{definition}
\label{def:Cnp}
Let $p \ge 1$, $n \ge 0$ be given natural numbers. Let $h = \frac{\tau}{p}$
be a grid step, $t_i = -i \cdot h$ be grid points for $i \in \{0, \dots, p\}$
and let intervals $I_i = [t_i, t_{i-1})$ for $i \in \{1, \dots, p\}$.

We define $C^n_p([-\tau, 0], \R^d)$ to be a set of functions
$x : [-\tau, 0] \to \R^d$ such that $x$ has a forward Taylor
representation of order $n$
on all $I_i$ and such that $x^{(n+1)}$ (understood as a right derivative)
is bounded over whole $[-\tau, 0]$.
\end{definition}
From now on we will assume that $\tau$ is fixed and we will write $C^n_p$ and
$C^{k}$ to denote $C^n_p([-\tau, 0], \R^d)$ and $C^k([-\tau, 0], \R^d)$,
respectively. Moreover, whenever we use $p$ and $h$ without additional
assumption, we assume that $h$ is given by $h = h(p, \tau) = \frac{\tau}{p}$
as in Def.~\ref{def:Cnp}.

Note that $x \in C^{n}_{p}$ might be discontinuous at $t = t_i$,
$i \in \{-p, \dots, 0\}$. However, $C^n_p \cap C^k$ is a
linear subspace of $C^k$ for any $k \in \N$
and if $k > n$ then obviously $C^k_p \subset C^n_p$ (see \cite{nasza-praca-focm}).
Therefore $\mathcal{X} = C^n_p \cap C^0$ can be used as a suitable
subspace of the phase space $C^0$ for solutions of Eq.~\eqref{eq:dde}. In fact,
following two lemmas, proved in \cite{nasza-praca-focm},
state that $\varphi(h, \cdot)$ and $\varphi(t, \cdot)$ for $t$
large enough are well defined maps $\mathcal{X} \to \mathcal{X}$:
\begin{lemma}
\label{lem:Cnp-smoothing}
Assume $f$ in \eqref{eq:dde} is $C^\infty$ (or smooth enough).
Let $\psi \in C^{n}_p$ be an initial function to \eqref{eq:dde}.
If $\varphi(h, \psi)$ exists
then $\varphi(h, \psi) \in C^{n}_{p}$.
Moreover, if $\psi \in C^n_p \cap C^0$ and
$i = k \cdot p$ for some $k \in \N$ then
$\varphi(i \cdot h, \psi) \in C^{n+k}_{p} \cap C^k$.
\end{lemma}

\begin{lemma}
\label{lem:Cnp-long-enough-time}
Assume $f$ in \eqref{eq:dde} is $C^\infty$ (or smooth enough).
Let $\psi \in C^{n}_p \cap C^0$ be initial function so that the
solution to \eqref{eq:dde} exists up to some $t \ge T$, where
$T = T(n, \tau) = (n+1) \cdot \tau$.
Then $\varphi(t, \psi) \in C^{n}_{p} \cap C^0$.
\end{lemma}
Time $T(n, \tau)$ will be important when constructing Poincar\'e maps later
in the paper, so to underline its importance, we state the following:
\begin{definition}
\label{def:long-enough-time}
We call $T(n)$ in Lemma~\ref{lem:Cnp-long-enough-time}
a \emph{long enough integration time}.
\end{definition}

In the current work we generalize the notion of the space $C^n_p$ to allow
different order of the jets at different points of the grid. This will be
beneficial to the final estimates later, as the representation of functions
will take advantage of the smoothing of solutions:
\begin{definition}
\label{def:Cetap}
Let $p$ be fixed, $\eta = (n_1, \dots, n_p) \in \N^p$ and let $t_i, I_i, h$ be
as in Definition~\ref{def:Cnp}. We define space of functions $C^\eta_p$ so that
$x \in C^\eta_p$ iff $x$ has a forward Taylor representation of order $n_i$ on
$I_i$ and $x^{(n_i + 1)}(I_i)$ is bounded for $i \in \{1,\dots,p\}$.
\end{definition}
The discussion from Section~\ref{subsec:integrator-overview} about the smoothing
of solutions of DDEs shows that if we have $n$-th order Taylor representation at
$t = -\tau$ then we can obtain $(n+1)$-th order representation of $x$ at
$t = 0$. Therefore, the order of the representation of solution will not decrease
during the integration, and it can increase, in general, only by one at a time (after
integration for a full delay). Therefore we introduce the following special
class of $C^\eta_p$ spaces. Let $q \in \{0, \dots, p\}$ by $C^n_{p, q}$ we will
denote the space $C^{\eta}_{p}$ with
\begin{equation*}
\eta_i = \begin{cases}
	n+1 	& i \le q \\
	n 	& i > q
\end{cases},
\end{equation*}
that is, the Taylor representation would be of order $n$ on grid points
$-\tau = t_p, t_{p-1}, \dots, t_{q-1}$ and of order $n+1$ on
$t_{q}, t_{q+1}, \dots, t_1 = h$. Among all $C^\eta_p$ spaces, spaces $C^n_{p,q}$ 
will be used most extensively in the context of rigorous integration of DDEs, 
but we keep the general notation of Definition~\ref{def:Cetap}
for simplicity of formulas later.

Now, it is easy to see that $C^{n}_{p, p} = C^{n+1}_{p, 0}$ and so that
$C^n_p = C^n_{p,0}$. Analogously we can write for $q > p$ that
$C^{n}_{p, q} = C^{n+\overline{q}}_{p, \underline{q}}$ with
$\overline{q} = \left\lfloor \frac{q}{p} \right\rfloor$ and
$\underline{q} = q \bmod p$. With that in mind the analogue of
Lemma~\ref{lem:Cnp-smoothing} can be stated as:
\begin{lemma}
\label{lem:Cetap-smoothing}
Let $\psi \in C^{n}_{p, q}$ be an initial function to \eqref{eq:dde} and let
$h$ be as in Def.~\ref{def:Cnp}.
If $\varphi(h, \psi)$ exists
then $\varphi(h, \psi) \in C^{n}_{p,q+1}$.
Moreover, if $\psi \in C^n_{p,q} \cap C^0$ and
$m = k \cdot p$ for some $k \in \N$ then
$\varphi(m \cdot h, \psi) \in C^{n+k}_{p,q} \cap C^k$.
\end{lemma}
\textbf{Proof:}
It follows from the smoothing of solutions,
the definition of $C^{n}_{p,q}$, equality of
spaces $C^{n}_{p,p} = C^{n+1}_{p,0}$
and by applying method of steps
(see e.g. \cite{dde-ode-book-driver}) to solve \eqref{eq:dde}.
\qed

In the rigorous method we will use Lemma~\ref{lem:Cetap-smoothing} as follows:
we will start with some set $X_0 \subset C^n_p = C^n_{p,0}$ defined with a
finite number of constraints. Then we will in sequence produce representations
of sets $X_i = \varphi(h, X_{i-1}) \in C^{n}_{p,i} = C^{n + \overline{i}}_{p,\underline{i}}$.
Finally, to compare sets defined in different $C^\eta_p$ spaces we would need the
following simple fact:
\begin{pro}
\label{pro:Cetap-subset}
$C^{\eta}_{p} \subset C^{\zeta}_{p}$ iff $\eta_i \ge \zeta_i$ for all $i \in \{1,\dots,p\}$.
\end{pro}
Now we  show how to describe sets in $C^{\eta}_p$.
Obviously, by the Taylor's theorem, we have that $x \in C^\eta_{p}$
is uniquely described by a tuple $\bar{x} = \left(z(x), j(x), \xi(x)\right)$,
where
\begin{itemize}
\item $z(x) := x(0) \in \R^d$,
\item $j(x) := (j_1(x),\dots,j_p(x))$ with $j_i(x) := J^{[n_i]}_{t_i}(x) \in \R^{d \cdot (n_i+1)}$,
\item $\xi(x) := (\xi_1(x),\dots,\xi_p(x))$ and $\xi_i(x) := x^{[n_i+1]}|_{I_i} \in C^0(I_i, \R^d)$ are bounded.
\end{itemize}
Please note, that the subscript $i$ denotes the grid point here, not the
component of the $x$ in $\R^d$. We will usually use subscript $j$ for this
purpose and we will write $z(x)_j$, $j_i(x)_j$, etc., but for now, 
all formulas can be interpreted simply for $d = 1$, generalization
to many dimensions being straightforward. We will use notation of $z(x)$,
$j(x)$, $\xi(x)$ etc. for a shorthand notation in formulas, sometimes dropping
the argument $x$ if it is known from the context. For example, we will say that
we have a solution described by a tuple $(z, j, \xi) \in \R^M \times (C^{0})^{p\cdot d}$,
then we will know how to interpret them to get the function $x$. Here
$M = M(p, \eta, d) = d \cdot (1 + \sum_{i=1}^{p} (\eta_i + 1))$.
A direct consequence is that:
\begin{pro}
The space $C^\eta_p$ is a Banach space isomorphic  
to $\R^M \times (C^{0})^{p \cdot d}$ by $x \mapsto (z(x), j(x), \xi(x))$, and
with a natural norm on $x$ given by
\begin{equation*}
\|x\|_{C^\eta_p} := \left\|\left(z(x), j(x)\right)\right\| + \sum_{i = 1}^{p}\sum_{j = 1}^{d} \sup_{t \in I_i} |\xi(x)_j(t)|,
\end{equation*}
where $\| \cdot \|$ denotes any norm in $\R^M$ (all equivalent).
We will use $\max$ norm in $\R^M$.
\end{pro}

Let now $\I$ be a set of all closed intervals over $\R$. We define:
\begin{eqnarray}
\label{eq:ixi_i(x)_j} \ixi_i(x)_j & := &
	\left[ \min_{\epsi \in [0, h]} \xi_i(x)_j(\epsi), \max_{\epsi \in [0, h]} \xi_i(x)_j(\epsi) \right] \in \I, \\
\notag \ixi_i(x) & := & \ixi_i(x)_1 \times \cdots \times \ixi_i(x)_d \in \I^{d}
\end{eqnarray}
and $\ixi(x) = \left(\ixi_1(x),\dots,\ixi_p(x)\right) \in \I^{d \cdot p}$.
That is a very complicated way to say $\ixi(x)$ is the collection of bounds
on the remainder terms in the Taylor representation of $x$.
The interval $\ixi_i(x)_j$ is well defined, since we assumed each $x^{(n_i+1)}$
bounded in Definition~\ref{def:Cetap}. Now, we can describe $x \in C^{\eta}_{p}$
by the following \emph{finite set} of numbers:
\begin{definition}
\label{def:peta-representation}
Let $M = M(p, \eta, d) = d \cdot \left(1 + \sum_{i=1}^{p} (\eta_i+1) \right)$.

We say that $\bar{x} = (z(x), j(x), \ixi(x)) \in \R^M \times \I^{d \cdot p}$
is a \emph{(p,$\eta$)-representation} of $x \in C^\eta_p$.

Given $\bar{x} \in \R^M \times \I^{d \cdot p}$ by
$X(\bar{x}) \subset C^\eta_p$ we denote the set of all functions whose
$\bar{x}$ is their \emph{(p,$\eta$)-representation}.
\end{definition}
The number $M$ is called the size of the representation and we will omit
parameters if they are known from the context. We will use shorthand notation of
$\R^{n}_{p}$, $\R^{\eta}_{p}$ or $\R^{n}_{p,q}$ to denote appropriate $\R^{M}$
in context of spaces $C^{n}_{p}$, $C^{\eta}_{p}$ and $C^{n}_{p,q}$,
respectively. We will write $\I_p$ to denote $\I^{p \cdot d}$. Note, that we are
dropping $d$ because it is always well known from the context.

Observe that, in general, $X(\bar{x})$ contains infinitely many functions. We
will identify $\bar{x}$ and $X(\bar{x})$, so that we could use notion of
$z(\bar{x})$, $j(\bar{x})$, etc. Moreover, we will further generalize the notion
of $X(\bar{x})$:
\begin{definition}
\label{def:pn-fset}
Let $A \subset \R^\eta_p$,
$R \in \I_p$ be a product of closed intervals.
We define set $X(A, R)$ as
\begin{equation*}
X(A, R) = \left\{ x \in C^{\eta}_{p} : (z(x), j(x)) \in A, \ixi(x) \subset R \right\}
\end{equation*}
We call $X(A, R)$ a \emph{(p,$\eta$)-functions set} (or (p,$\eta$)-fset
for short) and $(A, R)$ its \emph{(p,$\eta$)-representation}.
\end{definition}
If $A$ is convex then $X(A, R)$ is also a convex subset of
$C^{\eta}_p$, so $X(A, R) \cap C^k$ is also convex for any $k \in \N$,
see~\cite{nasza-praca-focm}. For a space $C^{n}_{p,q}$ we will use the
term (p,q,n)-representation and (p,q,n)-fsets when needed, but usually
we will use just names like ,,fset'' and ,,representation''.

Finally, we introduce the following shorthand symbols used for evaluation of terms:
\begin{eqnarray}
\label{eq:taylor-part} \taylor^n(j; \epsi) 		& := & \sum_{k=0}^{n} j_{[k]} \cdot \epsi^k\ , \\
\label{eq:summa-part}  \summa^n(\xi; \epsi) 	& := & (n+1) \cdot \int_0^{\epsi} \xi_i(s) \cdot (\epsi-s)^n ds, \\
\label{eq:eval}        \eval^n(j, \xi; \epsi)	& := & \taylor^n(j; \epsi) + \summa^n(\xi; \epsi),
\end{eqnarray}
for any function $\xi \in C^{0}([0, h), \R^d)$ and any jet
$j \in \R^{N \cdot d}$ of order $N \ge n$. The letters should be coined to the
terms $\taylor$ - (T)aylor sum, $\summa$ - (S)umma, formal name for the integral
symbol, $\eval$ - (E)valuation of the function. We use superscript $n$ to
underline order to which the operation applies, but in general, it can be simply
inferred from the arguments (for example - maximal order of the jet $j$ in
$\taylor$). Also, the superscript argument might be used to truncate computation
for higher order jets, e.g. let $j = J^{[2n]}_t x$ and consider applying
$\taylor^{n}(j)$ to Taylor-sum only part of the jet. This will be used in
algorithms later. If we omit the parameter $n$ then it is assumed that we use
the biggest possible $n$ (for that argument, inferred from the representation itself).

Then we will write formally for any $x \in C^\eta_p$:
\begin{eqnarray*}
\taylor^n(x; t) 	& := & \taylor^n(j_i(x); \epsi) , \\
\summa^n(x; t) 	& := & \summa^n(\xi_i(x); \epsi), \\
\eval^n(x; t)		& := & \taylor^n(j_i(x); \epsi) + \summa^n(\xi_i(x)(\cdot - t_i); \epsi),
\end{eqnarray*}
where $t = t_i + \epsi$, $\epsi \in [0, h)$.
For $X = X(A,R)$ we will write $a(X) = A$ and $\ixi(X) = R$ and for
$x \in X(A, R)$ we will write $a(x) = (z(x), j(x)) \in A$. We will
also extend the notion of operators $\taylor$, $\summa$ and $\eval$
to (p,$\eta$)-fsets:
\begin{eqnarray*}
\taylor(X; t) 	& := & \taylor^{\eta_i}(j_i(X); \epsi), \\
\summa(X; t) 	& := & \ixi_i(X) \cdot \epsi^{\eta_i+1}, \\
\eval(X; t) 	& := & \taylor^{\eta_i}(j_i(X); \epsi) + \summa^{\eta_i}(\ixi_i(X); \epsi).
\end{eqnarray*}
where $t = t_i + \epsi$, $\epsi \in [0, h)$.
Note, that $\taylor(x; t) = \taylor(\bar{x}; t)$,
$\summa(x; t) \in \summa(\bar{x}; t)$ and of course
$\eval(x; t) \in \eval(\bar{x}; t)$. In the rigorous
computation we as well might use intervals or whole
sets in the computation (e.g. $t = [t] = t_i + [0, \epsi]$) -
in such circumstances we will get sets representing
all possible results and in that way an estimate for
the true value. From now on, we will also drop
bar in $\bar{x}$ wherever we treat $x$ as an element
of $C^{\eta}_p$ with a known bounds in form
of some $X(A, R)$.

Finally, we make an observation that for $x$ - a solution to DDE~\eqref{eq:dde}
such that $x_{t_0} \in C^\eta_p$ - the $k$-th derivative $x_{t_0}^{[k]}$ must also
by representable by piecewise Taylor representation. In fact, since we know
$x(0)$ and all jets of the representation of $x_{t_0}$ we can obtain
$x_{t_0}^{[k]}(0)$ by applying Lemma~\ref{lem:taylor-reccurent-formula}, namely
Eq.~\eqref{eq:jet-reccurence-2}. Then, the value of all other jets and
remainders follows from Proposition~\ref{pro:forward-taylor-derivative}:
\begin{pro}
\label{pro:derivative-p-n-representation}
Let $x \in C^{\eta}_{p}$ be a segment of a solution to DDE~\eqref{eq:dde}
and for $k \in \N$ define $\eta - k := (\eta_1-k, \dots, \eta_p -k)$.
Then for $1 \le k \le \min_i{\eta_i}$ the derivative $x^{[k]}(t)$
(interpreted as a right derivative)
exists for $t \in [-\tau, 0]$ and
$x^{[k]} \in C^{\eta-k}_{p}$,
with a (p,$\eta-k$)-representation given in terms of the (p,$\eta$)-representation of $x$:
\begin{align}
j_i(x^{[k]}) & = \left(c^0_i,\dots,c^{\eta_i-k}_i\right) \notag \\
\xi_i(x^{[k]}) & = \binom{\eta_i+1}{k} \cdot \summa^{\eta_i-k}(\xi_i(x); \ \cdot \ ) \notag \\
\ixi_i(x^{[k]}) & \subset \binom{\eta_i+1}{k} \cdot \ixi_i(x), \notag \\
z(x^{[k]}) & = \frac{1}{k} \cdot \left(F_{k-1} \left( j_p(x), z(x) \right)\right)_{[k-1]}, \label{eq:z(x^[k])}
\end{align}
for $i \in \{1,\dots,p\}$, where
\begin{align*}
c^l_i & = \binom{l+k}{k} \cdot j_i(x)_{[l+k]}, & l \in 0,\dots,\eta_i-k.
\end{align*}
\end{pro}

%% file: sections/3_integrator.tex
\section{\label{sec:integrator}Rigorous integrator: basic algorithms and some improvements}

Now we are ready to show how to obtain estimates on the representation $Y$ of
$\varphi(h, X)$ for a given set of initial functions $X \in C^{\eta}_{p}$. Due
to the finite nature of the description of the set $Y$ we will have only the
relation $\varphi(h, X) \subset Y$, in general.

First, we want to recall in short the details of the integrator from
\cite{nasza-praca-focm} as those are crucial in the improvements presented
later. Then, we will show how to incorporate new elements: the extension of
the representation from \eqref{eq:jet-reccurence-2} and the spaces $C^\eta_p$,
the generalization to systems of equations (i.e. $d > 1$), and to multiple
delays (under the assumption that they match the grid points). Then, we will
discuss the Lohner-type method for the generalized algorithm.

\subsection{ODE tools}

We start with describing some ODE tools to be used in rigorously solving
\eqref{eq:non-autonomous-ode} using the computer. For this we will need a method
to find rigorous enclosures of the solution $x$ (and its derivatives w.r.t. $t$)
over compact intervals $[0, h]$. A straightforward method here is to consider
Eq.~\eqref{eq:dde} on $[t_0, t_0 + h]$, $h \le \tau$ as a non-autonomous ODE,
just as in the case of method of steps~\cite{dde-ode-book-driver}. If we plug-in
a \emph{known} initial function $x_{t_0}$ into \eqref{eq:dde} and we denote
$\hat{f}(z, t) := f(z, x_{t_0}(t-\tau))$ for $t\in [0, h]$ we end up with
non-autonomous ODE:
\begin{equation}
  \begin{cases}
    	z'(t) = \hat{f}(z, t), & t \in [0, h], \\
	   	z(0) = x_{t_0}(0).
  \end{cases} \label{eq:non-autonomous-ode}
\end{equation}
Please note that $t - \tau \in \domain\left(x_{t_0}\right) = [-\tau, 0]$ so $\hat{f}$ is
well defined, and $\hat{f}$ is of class $C^{k}$ as long as the solution segment
$x_{t_0}$ is of class $C^{k}$ (for $f$ sufficiently smooth). Therefore, in view
of \eqref{eq:taylor-jet-x} and \eqref{eq:jet-reccurence}, to find estimates on
the Taylor coefficients of $x$ over $I_{t_0} = [t_0, t_0 + h)$ it suffices only
to ascertain the existence of $z$ over $I_{t_0}$ and to have some finite
\emph{a priori} bounds $Z$ on it, as the estimates on the higher order
coefficients will follow from recurrent formulas \eqref{eq:taylor-jet-x} and
\eqref{eq:jet-reccurence}. Luckily, the existence of the solution to
Eq.~\eqref{eq:non-autonomous-ode} and a good a priori bounds over $I_{t_0}$ can
be obtained using existing tools for ODEs \cite{capd-zgliczynski,capd-lohner} as
was shown~in~\cite{nasza-praca-focm} and efficient implementations are already
available \cite{capd-article,www-capd}. We have the following:
\begin{lemma}[see Theorem~1 in \cite{capd-lohner}]
\label{lem:rough-encl}
We consider $\hat{f}$ as in
non-autonomous ODE~\eqref{eq:non-autonomous-ode}.

Let $B \subset \R^d$ be a compact set. 
If a set $W \subset \R^d$ is such that
\begin{equation*}
B + [0,\epsi] \cdot \hat{f}(W, [t_0,t_0+\epsi]) =: Z \subset W,
\end{equation*}
then, any solution $z$ of \eqref{eq:non-autonomous-ode} such that $z(t_0) \in B$
has $z(t_0 + \delta) \in Z$ for all $\delta \in [0,\epsi]$.
\end{lemma}
By $\roughEncl$ we denote a procedure (heuristic) to find the set $Z$:
\begin{equation*}
\roughEncl(f, B, t_0, \epsi) := Z, \textrm{ as in Lemma~\ref{lem:rough-encl}}.
\end{equation*}
We do not go into the details of this algorithm nor
the proof of Lemma~\ref{lem:rough-encl}, but we
refer to \cite{capd-lohner,capd-zgliczynski,www-capd} and
references therein.

\begin{rem}
Please note that finding a rough enclosure is a heuristic procedure and
therefore it is the point where the algorithm can fail (in fact the only one).
If that happens, we must abort computations or apply some strategy to overcome
the problem. In the ODE context it is possible to shorten the step or to
subdivide the set of initial conditions. Those strategies can be difficult to
adopt in the DDE context: we cannot shorten step because of the definition of
$C^\eta_p$ spaces and the loss of continuity problems discussed earlier; and we
could not afford extensive subdivision as we work with very high-dimensional
representations (projections) of functions. This makes obtaining the
higher order methods even more useful.
\end{rem}

Consider now $x_{t_0} \in C^n_p$, so that $\hat{f} \in C^{n+1}$. Applying
Eqs.~\eqref{eq:taylor-jet-x} and \eqref{eq:jet-reccurence} allows to obtain
$J^{[n+1]}_{t_0} x$, where rough enclosure procedure gives $Z$ such that
$x(I_{t_0}) \subset Z$. In what follows we will sum up all the formulas needed
to obtain (guaranteed enclosures on) the forward Taylor representation of $x$
on the interval $I_{t_0}$ of order $n+1$.

\subsection{\label{sec:full-step-integrator}The rigorous integrator in $C^{n}_{p,q}$}

Assume now that we are given some $x_0 \in C^n_{p,q}$. 
We will show how to compute rigorous estimates on a set
$X(A_h, R_h) \subset C^n_{p,q+1}$, with an explicitly given 
$A_h \subset \R^n_{p,q+1}$ and $R_h \in \I^{p\cdot d}$,
representing $\varphi(h, x_0)$, i.e. $\varphi(h, x_0) \in X(A_h, R_h)$.
The sets $A_h$ and $R_h$ will be computed using only data available
in $(z(x), j(x), \xi(x))$.
The subscript $h$ in $A_h$, $R_h$ is used to underline that we are making a full step
$h = \frac{1}{p}$. In what follows we will use the convention that $X_h = X(A_h, R_h)$.

This is an analogue to the algorithm described in Section~2.2 in
\cite{nasza-praca-focm}, but we account for the effect of smoothing of the
solutions in DDEs (Lemma~\ref{lem:Cetap-smoothing}), so that
$\varphi(h, x_0) \in C^{n}_{p,q+1}$ (and we remind that
$C^n_{p,p} = C^{n+1}_{p,0} = C^{n+1}_p$):
\begin{theorem}
\label{thm:integrator-h}
Let $x \in C^n_{p,q}$, with $0 \le q < p$ and the representation
$(z(x), j(x), \ixi(x)) \in \R^n_{p,q} \times \I^{d \cdot p}$.

We define the following quantities:
\begin{align}
\hat{f}		& := \textrm{ as in Eq.~\eqref{eq:non-autonomous-ode} } 	\notag \\
\label{eq:enclosure-delayed}
[c]_{[k]}		& := \eval \left(x^{[k]}; [-\tau, -\tau+h]\right) = \eval \left(x^{[k]}, [t_p, t_p + h]\right) \in \I^d, 	& 0 \le k \le n 	\\
\notag 						
[c]_{[n+1]} 	& := \ixi_p(x) \in \I^d 	\\
\notag 						
[Z] 			& := \roughEncl(\hat{f}, z(x), t_0, h) \in \I^d \\
\label{eq:rough-encl-recurrence} 	
[F] 			& := F^{[n+1]}\left([Z], [c]\right)
\end{align}
Then, we have for $y = \varphi(x, h)$  the following:
\begin{align}
\label{eq:shift-taylor} 		j_i(y) 		& = 		j_{i-1}(x)  	\ =: \ j_i(X_h) 									& i \in \{2,\dots,p\} \\
\label{eq:shift-remainder} 	\ixi_i(y) 	& = 		\ixi_{i-1}(x) \ =: \ \ixi_i(X_h) 								& i \in \{2,\dots,p\} \\
\label{eq:forward-jet} 		j_1(y) 		& = 		\left(z(x), w_{n+1} * F^{[n]}\left(z(x), j_p(x)\right) \right) \ =: \ j_1(X_h) 	&  	\\
\label{eq:forward-remainder} \ixi_1(y) & \subset \frac{1}{n+2} \cdot [F]_{[n+1]} \ =: \ \ixi_1(X_h),			&	\\
\label{eq:forward-taylor} 	z(y) 		& \in 		\taylor(j_1(y); h) + \left([F]_{[n+1]} \cdot [0, h]\right) \cdot h^{n+1} \ =: z(X_h) \ &
\end{align}
or, in other words, $y \in X_h \subset C^{n}_{p,q+1}$.
\end{theorem}
\textbf{Proof:}
Eq.~\eqref{eq:shift-taylor},\eqref{eq:shift-remainder} are representing the
shift in time by $h$ (one full grid point): from segment $x_0$ to segment
$x_h$ (of the solution $x$), therefore we simply reassign appropriate jets $j_i$
and remainders $\ixi_i$, as the appropriate grid points in both representations
overlap. The rest of formulas are an easy consequence of
Lemmas~\ref{lem:Cetap-smoothing}~and~\ref{lem:rough-encl}, the recurrence
relation \eqref{eq:taylor-jet-x} for $F^{[n]}$ and
Proposition~\ref{pro:derivative-p-n-representation} to obtain estimates on
$x^{[k]}$ over intervals $[-\tau, -\tau + h)$ in \eqref{eq:enclosure-delayed}.
Note, that the second term in \eqref{eq:forward-taylor} is formally given by the
integral remainder in Taylor formula~\eqref{eq:taylor-integral-jet}, namely for
$s \in [0, h)$ we have $\xi_1(y)(s) \in [\xi]_1(y) = \frac{1}{n+2} \cdot [F]_{[n+2]}$
(by the recurrence formula \ref{eq:taylor-jet-x}) and
\begin{equation*}
\summa \left(\xi_1(y), s\right) \in \summa \left(\ixi_1(y), [0, h]\right)
= (n+2) \cdot \left(\ixi_1(y) \cdot [0, h]\right) \cdot h^{n+1} = \left([F]_{[n+2]} \cdot [0, h]\right) \cdot h^{n+1}.
\end{equation*}
\qed

We denote the procedure of computing $X(A_h, R_h)$ for a given initial data $x
\in C^n_{p,q}$ by $\mathcal{I}$, i.e. $\mathcal{I}(x) = X(A_h, R_h)$. Clearly, it
is a multivalued function $\mathcal{I} : C^n_{p,q} \rightrightarrows C^n_{p,q+1}$.
We are abusing the notation here, as $\mathcal{I}$ is a family of maps (one for each
domain space $C^n_{p,q}$), but it is always known from the context 
(inferred from the input parameters).

We would like to stress again, that the increase of the order of representation
at $t = h$ in the solution $x$ will be very important for obtaining better
estimates later. It happens in Eq.~\eqref{eq:forward-jet}, as the resulting jet
is of order $n+1$ instead of order $n$ as it was in \cite{nasza-praca-focm}.
Please remember that $F^{[n]}$ is a recurrent formula for computing whole jet of
order $n$ of function $F = f \circ (x(\cdot), x(\cdot - \tau))$ at the current
time $t$, so it produces a sequence of coefficients, when evaluating
\eqref{eq:rough-encl-recurrence} and \eqref{eq:forward-jet}. Obviously, each of
those coefficients belongs to $\R^d$.

The nice property of the method is that the Taylor coefficients at $t = 0$, i.e.
$j_1(y)$ are computed \emph{exactly}, just like in the corresponding Taylor
method for ODEs (or in other words, if $x$ is a true solution to \eqref{eq:dde}
and $X_h = \mathcal{I}(x)$ then $J^{(n+1)}_h(x) = j_1(X_h)$).
It is easy to see, as formulas \eqref{eq:shift-taylor} and
\eqref{eq:forward-jet} does not involve a priori any interval sets (bracketed
notation, e.g. $\ixi$,$[Z]$, etc.). Therefore, to assess local error made by
the method we need only to investigate Eq.~\eqref{eq:forward-taylor}, which is
essentially the same as in the Taylor method for ODEs. As the interval bounds
are only involved in the remainder part $\left([F]_{[n+1]} \cdot [0, h]\right) \cdot h^{n+1}$
therefore, the local error of the method is $O(h^{n+2})$. Since $J^{(n+1)}_h(x)
= j_1(X_h)$ for a true solution $x$, this error estimation also applies to \emph{all} the
coefficients in the $j_1(y)$ computed in Eq.~\eqref{eq:forward-jet} in \emph{the next 
integration step}, when computing $X_{2h} = \mathcal{I}(X_h)$, as they 
depend on $z(X_h)$ that already contains the error. It will be also 
easily shown in numerical experiments (benchmarks) presented at the end of
this section.

\subsection{\label{sec:integrator-many-delays}Extension to many delays}

Now, we are in position to show how our algorithm can be generalized to include
the dependence on any number of delays $\tau_i$ as in Eq.~\eqref{eq:dde-many}, 
as long as they match with the grid points: $\tau_i = i \cdot h$. 
Therefore, we consider the following:
\begin{equation}
\label{eq:dde-many-delays}
x'(t) = f\left(x(t), x(t - p_1 h), x(t - p_2 h), \ldots, x(t - p_m h)\right),
\end{equation}
where $1 \le m \le p$ and $p = p_1 > p_2 > \ldots > p_m \ge 1$. We will denote
by $u(x_t) = u_f(x_t) := (x(t), x(t - p_1 h), x(t - p_2 h), \ldots, x(t - p_mh))$
the set of variables that are actually used in the evaluation of the r.h.s. $f$ in
Eq.~\eqref{eq:dde-many-delays} (as opposed to ,,unused'' variables, those at grid 
points not corresponding to any delays $\tau_i = p_i \cdot h$ 
in \eqref{eq:dde-many-delays}). 
This distinction will be important to obtain good computational complexity later on.
In case of Eq.~\eqref{eq:dde}, we have $u(x_t) = (x(t), x(t-\tau))$. Please
note, that since $u(x_t)$ contains variables at grid points, it is
easy to obtain $J_{t}^{(n)} u$ of appropriate order $n$. If $u = u(x_t)$,
we will use subscripts $u_{0}$, $u_{p_1}$, etc. to denote respective projections
onto given delayed arguments, and we use $u_{p_i,[k]}$ to denote their appropriate 
coefficients of the jet $J^{(n)}_{-\tau_i} x_t$. 

In order to present the method for many delays we need to redefine
$F(t) = (f \circ u \circ x) (t)$ and investigate
Eqs.~\eqref{eq:enclosure-delayed}-\eqref{eq:forward-taylor}. It is easy to see,
that the only thing which is different is $F$ and computation of its jets. Thus,
we rewrite the algorithm $F^{[n]}$ from Eq.~\eqref{eq:jet-reccurence} in terms
of $u = u(x)$:
\begin{eqnarray}
F^{[0]}(u) & := & f(u), \notag \\
F^{[k]}(u) & := & \left(J^{[k]}_{(u)}{f}\right) \circ_J \left( \left(u_0, w_k * F^{[k-1]}(u) \right), \left(u_{{p_1},[l]}\right)_{0 \le l \le k}, \ldots, \left(u_{{p_m},[l]}\right)_{0 \le l \le k} \right).
\label{eq:jet-reccurence-general}
\end{eqnarray}
Now, the algorithm from \eqref{eq:enclosure-delayed}-\eqref{eq:forward-taylor}
for an $x \in C^{\eta}_p$ consists of two parts.
First, the enclosure of the solution and all used variables over the basic
interval $[0, h]$:
\begin{align}
\notag 		\hat{f}	(t, z)		& := 		f(z, x(t-p_1 h), x(t-p_2 h),\ldots, x(t-p_m h))					& \\
\label{eq:neta}n				& := 		\min_{1 \le i \le m } \eta_{p_i}		=:   n(\eta,f)				& \\
\notag 		[U]_{{p_i},[k]}	& := 		\eval \left(x^{[k]}, [t_{p_i}, t_{p_i} + h]\right) \in \I^d, 	& 1 \le i \le m, \quad 0 \le k \le n 	\\
\notag 		[U]_{{p_i},[n+1]} & := 		\ixi_{p_i}(x) \in \I^d 												& 1 \le i \le m \\	
\label{eq:U0} 	[U]_0			& := 		\roughEncl(\hat{f}, z(x), t_0, h) \in \I^d 						& \\
\label{eq:FU} 	[F]				& := 		F^{[n+1]}\left([U]\right),											& 
\end{align}
then, building the representation after the step $h$:
\begin{align}
\notag 		j_i(y) 				& = 		j_{i-1}(x)  \ =: \ j_i(X_h) 											& i \in \{2,\dots,p\} \\
\notag 		\ixi_i(y) 			& = 		\ixi_{i-1}(x) \ =: \ \ixi_i(X_h) 										& i \in \{2,\dots,p\} \\
\notag 		j_1(y) 				& = 		\left(z(x), w_{n+1} * F^{[n]}\left(u(x)\right) \right) \ =: \ j_1(X_h) 	&  	\\
\notag 		\ixi_1(y) 			& \subset 	\frac{1}{n+2} \cdot [F]_{[n+1]} \ =: \ \ixi_1(X_h),					&	\\
\notag 		z(y) 				& \in 		\taylor(j_1(y); h) + \left([F]_{[n+1]} \cdot [0, h]\right) \cdot h^{n+1} \ =: z(X_h) \	&
\end{align}
Please note that we used in \eqref{eq:U0} symbol $[U]_0$ to denote enclosure of
$x$ over $[0, h]$ (computed by the $\roughEncl$ procedure). All other components
of $[U]$ are computed estimates on jets $j_{p_i}(x)$ over the same interval $[0, h]$
using Proposition~\ref{pro:forward-taylor-derivative}.
That way, we can think of $[U]$ as the enclosure of $u$ over interval $[0, h]$.
We have also generalized the algorithm to be valid for any $C^\eta_p$ by
introducing the notion of $n(\eta,f)$ in Eq.~\eqref{eq:neta}. The $n(\eta,f)$
depends on $f$ in the sense, the minimum is computed only for $n_i$ that are
actually used in computations.

\subsection{\label{sec:epsi-steps}Steps smaller than $h$}

In this section we consider computation of the (p,n)-representations of
$\varphi(t, x_0)$ where $t$ is not necessary the multiple of the basic step size
$h = \frac{\tau}{p}$, and for the initial $x_0 \in C^{\eta}_{p}$,
where the apparent connection between $\eta$, $n$ and $t$ will be discussed soon.
This problem arises naturally in the construction of Poincar\'e maps. Roughly speaking, the
Poincar\'e map $P$ for a (semi)flow in the phase space $\mathcal{X}$ is defined
as $P(x) = \varphi(t_P(x), x)$, where $x \in S \subset \mathcal{X}$ and
$t_P : S \to (0, \infty)$ - the return time to the section $S$ - is a
\emph{continuous} function such that $\varphi(t_P(x), x) \in S$ (we skip the
detailed definition and refer to \cite{nasza-praca-focm}). We see that the
algorithm presented so far is insufficient for this task, as it can produce
estimates only for \emph{discrete times} $t = i \cdot h$, $i \in \N$, not for a
possible continuum of values of $t_P(S)$. It is obvious that we can express
$t = m \cdot h + \epsi$ with $m \in \N$ and $0 < \epsi < h$ and the computation
of $\varphi(t, x_0)$ can be realized as a composition $\varphi(\epsi, \varphi(m
\cdot h, x_0))$. Therefore, we assume that the initial function is given as
$x_m = \varphi(m \cdot h, x)$ and we focus on the algorithm to compute
(estimates on) $x_\epsi = \varphi(\epsi, x_m)$.

First, we observe that, for \emph{a general} $x_m$ in some $(p, \eta)$-fset, we
\emph{cannot expect} that $x_\epsi \in C^{\zeta}_p$ for any $\zeta$.
The reason is that the solution $x$ of DDE~\eqref{eq:dde} with initial 
data in $C^\eta_p$ can be of class
as low as $C^0$ at $t=0$, even when the r.h.s. and the initial data is smooth
(as we have discussed in the beginning of Section~\ref{sec:phasespace}). The
discontinuity appears at $t = 0$ due to the very nature of Eq.~\eqref{eq:dde}.
This discontinuity is located at $s = -\epsi$ in the segment $x_\epsi$ of the
solution and, of course, we have $-\epsi \in [-h, 0]$. Therefore, the function
$x_\epsi$ does not have any Taylor representation (in the sense of
Def.~\ref{def:forward-taylor-representation}) on the interval $I_1 = [-h, 0]$,
as the first derivative of $x$ is discontinuous there.

On the other hand, we are not working with a general initial function,
but with $x_m = \varphi(m \cdot h, x_0)$, with $x_0 \in C^\eta_p$.
From Lemma~\ref{lem:Cetap-smoothing}
we get that $x_m \in C^{\eta + n + 1}_{p} \cap \mathcal{C}^{n + 1}$,
where $n \in \N$ be the largest value
such that $m \ge (n + 1) \cdot p$. Moreover, the same is true for
$x_{m+1} = \varphi(h, x_m)$. Therefore
$x_\epsi = \varphi(\epsi, x_m) \in \mathcal{C}^{n+1}$,
so that it has a $C^{n}_{p}$ representation.

Now, the question is:
can we estimate this (p,n)-representation in terms of
the coefficients of representations of $x_m$ (and maybe $x_{m+1}$)?
The answer is positive, and we have:
\begin{lemma}
\label{lem:epsi-step-formulas}
Assume $x$ is a solution to \eqref{eq:dde} with
a segment $x_0 \in C^{\eta}_{p} \cap \mathcal{C}^{0}$. Let $t \in \R$
be given with $t = m \cdot h + \epsi$, $m \in \N$, $0 < \epsi < h$.
Let $n = \lfloor \frac{m}{p}\rfloor - 1$ and assume $n \ge 0$, i.e.
$m \ge p$ and $t \ge \tau$.

Let denote $x_m = \varphi(m \cdot h, x_0)$ and
$x_{m+1} = \varphi(m \cdot h + h, x_0)$
and for $i \in \{1,\ldots, p\}$ let
\begin{align}
\label{eq:[L]} [L]_{i} &= \eval \left(j_i(x_m^{[n+1]}), \ixi_i(x_m^{[n+1]}), [0, h]\right), \\
\label{eq:[R]} [R]_{i} &= \eval \left(j_i(x_{m+1}^{[n+1]}), \ixi_i(x_{m+1}^{[n+1]}), [0, \epsi]\right).
\end{align}
Then we have $x_t \in X_\epsi \subset C^{n}_{p} \cap C^{n+1}$ for $X_\epsi$ given by:
\begin{align}
\label{eq:epsi-step-z} z\left(X_\epsi\right) & := \taylor\left(j_1(x_{m+1}); \epsi\right) + \summa\left(\ixi_1(x_{m+1}); \epsi\right), & \\
\label{eq:epsi-step-j} j_{i,[k]}\left(X_\epsi\right) & := \taylor\left(j_i(x_m^{[k]}); \epsi\right) + \summa\left(\ixi_{i}(x_m^{[k]}); \epsi\right), & i \in \{1,\ldots, p\}, k \in \{0, \ldots, n\}, \\
\label{eq:epsi-step-R-0} \ixi_i\left(X_\epsi\right) & := \hull\left([L]_{i}, [R]_{i}\right), & i \in \{1,\ldots, p\}.
\end{align}
\end{lemma}
Before the proof, we would like to make a small comment. The 
representation of $x_{m+1}$ is used for optimization and simplification 
purposes, as usually we have it computed nevertheless (when finding the
crossing time of the Poincar\'e map). It contains the representation
of $x$ over $[mh, mh + h)$ in $j_1$. Otherwise we would need to expand 
the jet of solution $x$ at $t = 0$ to compute $[R]_1$ and $z$ in \eqref{eq:epsi-step-z}. 
Also, the formula \eqref{eq:[R]} would be less compact.

\textbf{Proof of Lemma~\ref{lem:epsi-step-formulas}:} It is a matter of simple calculation.
To focus the attention on the $\epsi$ step,
let us abuse notation and denote $x_t = x_\epsi = \varphi(\epsi, x_m)$.
We have
\begin{equation*}
x^{[k]}_\epsi(-i \cdot h) = x^{[k]}_m(-i \cdot h + \epsi), \quad i \in \{1, \ldots, p \}
\end{equation*}
so we get a straightforward formula:
\begin{align}
\label{eq:jet-epsi-step-eval-form}
j_i(x_\epsi)_{[k]} & = \eval\left(x_m^{[k]}; -i \cdot h + \epsi\right) \ = \
\taylor^{\eta_i + n + 1 - k}\left(j_i(x^{[k]}_m); \epsi\right) + \summa^{\eta_i + n + 1 - k}\left(\xi_i(x_m^{[k]}); \epsi\right)
\end{align}
where representations of $x^{[k]}_{m}$ are obtained
by applying Proposition~\ref{pro:derivative-p-n-representation}.
Similarly, one can find that
\begin{equation}
\label{eq:z-epsi-step-eval-form} z(x_\epsi) = x_\epsi(0) = x_{m+1}(-h + \epsi) = \ \taylor(j_1(x_{m+1}); \epsi) + \summa(\ixi_1(x_{m+1}); \epsi),
\end{equation}
and for $s \in [0, h)$, $i \in \{1, \ldots, p \}$:
\begin{equation}
\label{eq:xi-epsi-step-eval-form} \xi_i(x_\epsi)(s) = x^{[n+1]}_\epsi(-i \cdot h + s) =
\begin{cases}
x^{[n+1]}_m(-i \cdot h + \epsi + s) = \eval\left(x^{[n+1]}_m; \epsi + s\right)							&	\epsi + s < h \\
x^{[n+1]}_{m+1}(-i \cdot h + (\epsi + s - h)) = \eval\left(x^{[n+1]}_{m+1}; (\epsi + s - h)\right) 	&	\epsi + s \ge h
\end{cases}.
\end{equation}
Note, in the second case of Eq.~\eqref{eq:xi-epsi-step-eval-form}
we have $0 \le (\epsi + s - h) < h$. Now, we exchange each $\xi$ with $\ixi$
in Eqs.~\eqref{eq:jet-epsi-step-eval-form}-\eqref{eq:xi-epsi-step-eval-form}
to get the corresponding estimates in Eqs.~\eqref{eq:epsi-step-z}-\eqref{eq:epsi-step-R-0}.
\qed

This algorithm is valid for any number of dimensions and for any
number of delays (i.e. for any definition of used variables $u(n,f)$) -
in fact, there is no explicit dependence on the r.h.s of \eqref{eq:dde} in the formulas
- the dynamics is ,,hidden'' implicitly in the already computed jets $j_{i}(x_m)$ and $j_1(x_{m+1})$.
This form of the algorithm will allow in the future to make general improvements
to the method, without depending on the actual formula for the projection of used variables
$u(n,f)$ in the r.h.s. of DDE~\eqref{eq:dde-many}, or even when 
constructing methods for other forms of Functional Differential Equations.
We will denote the $\epsi$ step algorithm given by 
\eqref{eq:epsi-step-z}-\eqref{eq:epsi-step-R-0} by $\mathcal{I}_\epsi$.

As a last remark, similarly to the discussion in the last 
paragraph of Section~\ref{sec:full-step-integrator}, let us consider the 
order of the local error in the method $\mathcal{I}_\epsi$.
This local error will have a 
tremendous impact on the computation of Poincar\'e maps, and thus
on the quality of estimates in computer assisted proofs. 
To see why, set the order $n$ and let us consider two maps: $T = \varphi(m h, \cdot)$ 
and $T_\epsi = \varphi(m h + \epsi, \cdot)$, where, without loss of generality,
we choose $m = p \cdot (n+1)$ (in applications, return time in Poincar\'e maps will 
be required to be greater than this) and we fix some $0 < \epsi < h$. 
It is of course sufficient to use full step method $\mathcal{I}$ to rigorously compute
map $T$, while $T_\epsi$ is a good model
of computing estimates on a real Poincar\'e Map and will require usage of
$\mathcal{I}_\epsi$ in the last step. Let us denote $x_m = T(x_0)$,
$x_{m+1} = \varphi(h, T(x_0))$ and $x_\epsi = T_{\epsi}(x_0)$. 
Obviously we have $x_{m+1} = \varphi(h, x_m) \in \mathcal{I}(x_m)$
and $x_\epsi = \varphi(\epsi, x_m) \in \mathcal{I}_\epsi(x_m)$
Assume $x_0 \in C^\eta_p$ with uniform order on all grid points,
$\eta = n$. From Lemma~\ref{lem:Cetap-smoothing} for both maps we end up with
$x_m \in C^{2n + 1}_{p} \cap \mathcal{C}^{n+1}$, 
$x_{m+1} \in C^{2n + 2}_{p,1} \cap \mathcal{C}^{n+1}$
and $x_\epsi \in C^{n+1}_{p}  \cap \mathcal{C}^{n+1}$. From discussion
in the last 
paragraph of Section~\ref{sec:full-step-integrator}, we can infer
that the local error introduced in $\mathcal{I}(x_m)$
is of order $O(h^{2n + 2})$, as the only term with non-zero
Taylor remainder is $z(\mathcal{I}_\epsi(x_m))$. Therefore,
we can expect that 
the accumulated error of estimating map $T$ with $\mathcal{I}^m$ ($m$ steps
of the full step integrator  $\mathcal{I}$) is
of order $O(h^n)$ \cite{book-adaptative-steps}, 
as this is the accumulated error of covering the first 
delay interval $[0, \tau)$ in the beginning of the integration process. 
Later, thanks to smoothing of solutions
and expanded space, the subsequent errors would be of higher order. 
This in general should apply \emph{even if we do not 
expand the representation}, as in such case the local error
in each step (even after $[0, \tau]$) is still 
$\mathcal{I}$ is still just $O(h^{n+1})$.

I comparison, algorithm $I_\epsi$ evaluates Taylor
expansion  with non-zero remainder not only at $z(\cdot)$ in \eqref{eq:epsi-step-z},
but at \emph{every} grid point and \emph{every} coefficient order 
of the representation in \eqref{eq:epsi-step-j}. What is more, the
impact of the remainder term $[\xi]$ is of different order
at different Taylor coefficients. 
Here we use Proposition~\ref{pro:derivative-p-n-representation}
to get that $k$-th Taylor coefficient $x_m^{[k]}$ has a (p,l)-representation with
$l = 2n + 1 - k$, so the local error of $j_{i,{[k]}}(X_\epsi)$ is of order $O(h^{2n + 1 - k})$.
Since $k \in \{0, \ldots, n\}$, then in the worst case of $k = n$, 
the local error size is $O(h^{n+1})$. 
This is of course worse than $O(h^{2n + 2})$
of the full step method, but it is a significant improvement over the first
version of the algorithm presented in \cite{nasza-praca-focm}, where
the local error of the last $\epsi$ step was $O(h)$ (basically, because
$x^{[n]}_\epsi$ was computed by explicit Euler method in the non-expanded representation
of $x_m \in C^n_p$). Current error is of the order comparable to the accumulated error
over the course of a long time integration $\mathcal{I}^m$, therefore has
a lot less impact on the resulting estimates. 

Exemplary computations, supporting
the above discussion, are presented in 
Section~\ref{sec:benchmark}.

\subsection{Computation of Poincare maps}

In this section we would like to discuss shortly some minor changes
to the algorithm of computing image of Poincar\'e map
using algorithms $\mathcal{I}$ (full step $h$) and $\mathcal{I}_\epsi$
($\epsi < h$), particularly, we discuss the case when the estimate on $t_P(S)$
has diameter bigger than $h$ - this will be important in one of
the application discussed in this paper.

In the context of using rigorously computed images of Poincar\'e maps in
computer assisted proofs in DDEs, we will usually do the following
(for details, see \cite{nasza-praca-focm}):
\begin{enumerate}
\item \label{item:cap-choose-space}
	We choose subspace of the phase-space of the semiflow $\varphi$
	as $C^n_p \cap \mathcal{C}^0$ with $p$, $n$ fixed.
\item \label{item:cap-choose-section}
	We choose sections $S_1, S_2 \subset C^{n}_{p}$, usually as some hyperplanes
	$S_i = \{ x \in C^{n}_{p} : S_i(x) := (s_i\ .\ a(x)) - c_i = 0 \}$, with $s_i \in \R^{M(d, p, n)}$, 
	$c \in \R$ and $(\ .\ )$ denoting the standard scalar product in $\R^{M(d, p, n)}$
	(we remind $a(x) = (z(x), j(x)) \in \R^{M(d,p,n)}$, $M(d,p,n) = d \cdot (1 + (n+1) \cdot p)$).
	Of course, in the simplest case, we can work only with a single section, $S_1 = S_2$.
\item  \label{item:cap-choose-set}
	We choose some initial, closed and convex set $X_0 \subset S_1 \subset C^n_p$
	on the section $S_1$.
\item \label{item:cap-tp}
	We construct $[t] \in \I$  such that $t_P(X_0) \subset [t]$, where
	$t_P : X_0 \to \R_+$ is the return time function from $X_0$ to $S_2$,
	so that $\varphi(t_P(x_0), x_0) \in S_2 \subset C^n_p$ for
	all $x_0 \in X_0$. This is done usually alongside the computation of the image $P(X_0)$,
	by successive iterating $X_{j+1} = \mathcal{I}(X_{j})$ until $X_m$ is
	\emph{before} and $X_{m+1}$ is \emph{after} the section $S_2$ (i.e. $S_2(X_m) < 0$
	and $S_2(X_{m+1}) > 0$ or $S_2(X_m) > 0$ and $S_2(X_{m+1}) < 0$).
	In such a case $[t] = m \cdot h + [\epsi]$, where $[\epsi] \subset [0, h)$.
	
	In view of
	Lemma~\ref{lem:epsi-step-formulas} we require
	$t_P(X_0) \ge (n+1) \cdot \tau$ - the return time to the section is
	\emph{long enough}. Moreover, $X_m$
	and $X_{m+1}$ are already computed to be used in the formulas
	\eqref{eq:epsi-step-z}-\eqref{eq:epsi-step-R-0}.
	The tight estimates on $[t]$ can be obtained for example with
	the \emph{binary search} algorithm, in the same manner as it was
	done in \cite{nasza-praca-focm}.
	
	Finally, using formulas from Lemma~\ref{lem:epsi-step-formulas}
	we get $X_\epsi \subset C^n_p$ such that $\varphi([\epsi], X_m) \subset X_\epsi$. 
\item \label{item:cap-theorem}
	We use sets $X_0$ and $X_\epsi$ together with the estimates on
	$P(X_0) \subset \varphi([t], X_0)$
	to draw conclusion on existence of some interesting dynamics.
	For example, if $S_1 = S_2$ and $P(X_0) \subset X_0$
	we can use Schauder Fixed Point Theorem to show existence of a periodic
	point of $P$ (the compactness of the operator $P$ plays here a crucial role).
\end{enumerate}
Now, we have already mentioned that the computation of the Poincar\'e map
$P(x_0) = \varphi(t_P(x_0), x_0)$ can be done by splitting the return time
$t_P(x_0) = m(x_0) \cdot h + \epsi(x_0)$ with $m(x_0) \in \N$ and
$\epsi(x_0) \in (0, h)$. This leads to a rough idea of rigorous algorithm
to compute estimates on $P(x_0)$ in the following form:
\begin{equation}
\label{eq:P-rigorous-single-point}
P(x_0) \in \mathcal{I}_{\epsi(x_0)} \circ \mathcal{I}^{m(x_0)} \left(x_0\right).
\end{equation}
However, in the case of computing (estimates on) $P(X_0)$
for a whole set $X_0 \subset C^n_p$, we can face the following problem:
for $x, y \in X_0$ we can have $m(x) \ne m(y)$, especially, when $X_0$ is large.
In \cite{nasza-praca-focm} we have simply chosen $X_0$ so small, such that
$m(x)$ is constant in $X_0$. Then, we have $[t] = m \cdot h + [\epsi]$, with
$[\epsi] = [\epsi_1, \epsi_2]$, $0 < \epsi_1 \le \epsi_2 < h$.
In such a situation formula
\eqref{eq:P-rigorous-single-point} could be applied with $m(x) = m$ and $\epsi = [\epsi]$.
In the current work we propose to take the advantage of all the data already stored in
the $(p,\eta)$-fsets and to extend the algorithm in Lemma~\ref{lem:epsi-step-formulas}
to produce rigorous estimates on $\varphi([\epsi], x_0)$ for
$[\epsi] = [\epsi_1, \bar{m} + \epsi_2]$, $0 < \epsi_i < h$
$\bar{m} \in \N$. It is not difficult to see that we have the following:
\begin{pro}
\label{pro:epsi-step-formulas-long-epsi}
Let $[t] = m \cdot h + [\epsi_1, \bar{m} \cdot h + \epsi_2]$ with
$0 < \epsi_1, \epsi_2 < h$, $m, \bar{m} \in \N$ with $\bar{m} > 0$.
Let assume $X_j$ are such that $\varphi(j \cdot h, X_0) \subset X_M$ for
$j = m, m+1, \ldots, m+\bar{m}+1$. Finally, let $n$ be as in
Lemma~\ref{lem:epsi-step-formulas}.

We define ($k \in \{ 0, \ldots, n+1 \}$, $j \in \{0, \ldots, \bar{m} + 1\}$, $i \in \{ 1, \ldots, p\}$):
\begin{align}
\notag [L]^{[k]}_{i,j} &= \eval \left(j_i(x_{m+j}^{[k]}), \ixi_i(x_{m+j}^{[k]}), [\epsi_1, h]\right),  \\
\notag [C]^{[k]}_{j,i} &= \eval \left(j_i(x_{m+j}^{[k]}), \ixi_i(x_{m+j}^{[k]}), [0, h]\right),  \\
\notag [R]^{[k]}_{i,j} &= \eval \left(j_i(x_{m+j}^{[k]}), \ixi_i(x_{m+j}^{[k]}), [0, \epsi_2]\right),
\end{align}
and a set $X_\epsi$ given by:
\begin{align}
z(X_\epsi) & := \hull \left([L]^{[0]}_{1,1}, [C]^{[0]}_{1,2}, \ldots, [C]^{[0]}_{1,\bar{m}}, [R]^{[0]}_{1,\bar{m} + 1}\right) , & \\
j_{i,[k]}(X_\epsi) & := \hull \left( [L]^{[k]}_{i,0}, [C]^{[k]}_{i,1}, \ldots, [C]^{[k]}_{i,\bar{m}-1}, [R]^{[k]}_{i,\bar{m}} \right), & i \in \{1,\ldots, p\}, k \in \{0, \ldots, n\}, \\
\ixi_i(X_\epsi) & := \hull\left([L]^{[n+1]}_{i,0}, [C]^{[n+1]}_{i,1}, \ldots, [C]^{[n+1]}_{i,\bar{m}}, [R]^{[n+1]}_{i,\bar{m}+1}\right), & i \in \{1,\ldots, p\}.
\end{align}

Then for all $t \in [t]$ we have  $x_t \in X_\epsi \subset C^{n}_{p} \cap C^{n+1}$.
\end{pro}
Of course, in the case $m(X_0) = const$ we use algorithm from
Lemma~\ref{lem:epsi-step-formulas}.

\subsection{The Lohner-type control of the wrapping effect}
\label{subsec:lohner-control-wrapping}
An important aspect of the rigorous methods using interval arithmetic is an
effective control of the \emph{wrapping effect}. The wrapping effect
occur in interval numerics, when the result of some non-linear operation 
or map needs to be enclosed in an interval box. When this box is
chosen naively, then a huge overestimates may occur, see 
Figure~\ref{fig:wrapping-effect} in Appendix~\ref{app:lohner}. 

To control wrapping effect in our computations we employ the Lohner
algorithm \cite{capd-lohner}, by representing sets in a good local coordinate 
frame: $X = x_0 + C \cdot r + E$, where $x_0$ is a vector in $\R^M$, $C \in \matrices(M, N)$,
$r_0 \in \I^N$ - an interval box centred at $0$, and $E$ some representation of local
error terms. As it was shown in \cite{nasza-praca-focm}, taking $E \in \I^M$ 
(a interval form of the error terms) was enough to prove existence of periodic orbits. 
Moreover, taking into account the form of the algorithm given by
\eqref{eq:shift-taylor}-\eqref{eq:forward-taylor} (especially the shift part
\eqref{eq:shift-taylor}-\eqref{eq:shift-remainder}) to properly
reorganize computations was shown to be crucial 
to obtain an algorithm of optimal computational complexity. 

In this work, we not only adopt this optimized Lohner algorithm to the 
systems of equations and to many delays, but we also propose
another form of the error term $E$ to get better estimates on the solutions
in case of systems of equations, $d > 1$, much in the same way
it is done for systems of ODEs \cite{capd-lohner, capd-article}. 
The proposed algorithm does not sacrifice
the computational complexity to obtain better estimates. We use this
modified algorithm in our proof of the symbolic dynamics in a delay-perturbed
R\"ossler system. 

The details of the algorithm are highly technical, so we decided to put 
them in the Appendix~\ref{app:lohner}, to not overshadow the presentation 
of the theoretical aspects, but on the other hand to be accessible for people 
interested in actual implementation details and/or in re-implementing presented 
methods on their own. 

\subsection{\label{sec:benchmark}Benchmarks}

As the last remark in this section, we present the
numerical experiment showing the effect of using
the new algorithm with expanding representation
in comparison with the old algorithm in \cite{nasza-praca-focm}.
As a test, we use a constant initial function $x_0(t) = 1.1$ for
$t \in [-\tau, 0]$ and the Mackey-Glass equation with parameter values 
$\beta = 2$, $\gamma = 1$, $n = 8$ and $\tau = 2$.
The configuration of (d,p,n)-fset $X_0$ has $n = 4$ (order 4 method), 
$p = 128$, $d = 1$ (scalar equation). The initial diameter
of the set $X_0$ is $0$. The test does integration over the $3n$ full delays 
(so that the final solution is smoothed enough). Then an 
$\epsi$-step is made, with the step $\epsi = \frac{h}{2}$, where
$h = \frac{\tau}{p}$ is the grid size (full step). 
In the Table~\ref{tab:benchmark} we present the maxima
over all diameters of the coefficients of the sets: $X_{3n} = \mathcal{I}^{3n}(X_0)$
that contains the segment $x_{3n}$ of the solution, and 
$X_{3n + \epsi} = \mathcal{I}_{\epsi} \left(\mathcal{I}^{3n}(X_0)\right)$.
We remind that $\mathcal{I}$ denotes the full-step integrator method
that does one step of size $h$, while $\mathcal{I}_{\epsi}$ is the 
$\epsi$-step method. Each maximum diameter is computed over all
Taylor coefficients of a given order $0 \le k \le 4$. We also
show the maximum diameter of the $\Xi$ part (order $k = 5$).

We test several maximal orders of the expanded representations: $2n$, $2n+1$ and $3n$.
The last one is the maximal order obtainable with the $3n$ full-delay integration steps, 
while the first one is the minimal reasonable one - taking into account the 
long enough integration time, see Def.~\ref{def:long-enough-time} and 
Lem.~\ref{lem:Cnp-long-enough-time}. 

\begin{rem}
Using the diameter $0$ of the set $X_0$ in the test will
show how the local errors of the method at each step affect
the final outcome.  
\end{rem}
 
\input tables/table1-content.tex

From Table~\ref{tab:benchmark} we see that the diameters of
the sets integrated with the new algorithm are far superior to
the old one. One can observe in a) that for the fixed number of full steps 
both methods produce results with coefficients of all orders
of a comparable diameter. This indicates that both methods 
are of order $h^4$. However, new algorithm produces
estimates of three orders of magnitude better. This is
because internally, the algorithm becomes of higher order 
after each full delay. After $k$ full delays, the actual 
order of the method is $n + k$. The second big advantage  
is shown in the b) part, where we have diameters of coefficients
after a small $\epsi$ step. This simulates for example computation
of a Poincar\'e map. The old algorithm produces estimates that
depend on the order of coefficient: the coefficient
$0$ has a diameter proportional to $h^n$, however, other
coefficients are computed with worse accuracy. 
The $4$'th order coefficient is computed with the lowest
accuracy of order $h^1$. On the contrary, the new algorithm
still retains the accuracy of the full step size algorithm
and produce far superior estimates (several orders of magnitude better).

The data and programs used in those computations
are described more in detail in Appendix~\ref{app:data}.

%% file: tables/table1-content.tex
\begin{table}[h!]
\begin{center}
\small
\begin{tabular}{|c|c|c|c|c|c|}
    \hline
    \multicolumn{6}{|l|}{\textbf{a) The set $X_{3n}$ after a fixed number of full steps - 12 full delays}} \\
    \hline
    Order $k$ & $h^k$ &  No expand & Expand $n$ & Expand $n + 1$ & Expand $2n$\\ 
    \hline
    $0$ $\star$    & 1 &  $8.0928124e-07$ & $1.3894812e-09$ & $1.3890612e-09$ & $1.3890594e-09$\\ 
    \hline 
    $1$ $\star$    & 0.015625 &  $2.0313339e-06$ & $3.4887294e-09$ & $3.4876694e-09$ & $3.487666e-09$\\ 
    \hline 
    $2$ $\star$    & 0.00024414062 &  $2.2627332e-06$ & $3.9124373e-09$ & $3.9113023e-09$ & $3.9113028e-09$\\ 
    \hline 
    $3$ $\star$    & 3.8146973e-06 & $2.096601e-06$ & $3.6231229e-09$ & $3.6220176e-09$ & $3.6220075e-09$\\ 
    \hline 
    $4$ $\star$    & 5.9604645e-08 & $3.1646014e-06$ & $5.5100828e-09$ & $5.508467e-09$ & $5.5084535e-09$\\ 
    \hline 
    $5$ $\dagger$    & 9.3132257e-10 & $0.14380491$ & $0.044424773$ & $0.044424773$ & $0.044424773$\\ 
    \hline 
    \multicolumn{6}{|l|}{\textbf{b) The final set $X_{3n + \epsi}$ after applying $\epsi$-step to $X_{3n}$}} \\
    \hline
    Order $k$ & $h^k$ & No expand & Expand $n$ & Expand $n + 1$ & Expand $2n$\\ 
    \hline
    $0$ $\star$ & 1 & $8.254823e-07$ & $1.4173127e-09$ & $1.4168844e-09$ & $1.4168826e-09$\\ 
    \hline 
    $1$ $\star$   & 0.015625 & $2.0673499e-06$ & $3.5503207e-09$ & $3.5492428e-09$ & $3.5492394e-09$\\ 
    \hline 
    $2$ $\star$    & 0.00024414062 & $2.4780643e-06$ & $3.9715719e-09$ & $3.9703922e-09$ & $3.970392e-09$\\ 
    \hline 
    $3$ $\star$    & 3.8146973e-06 & $8.9902593e-05$ & $3.9834122e-09$ & $3.7954002e-09$ & $3.7904426e-09$\\ 
    \hline 
    $4$ $\star$   & 5.9604645e-08 & $0.0056199777$ & $4.8690342e-08$ & $7.2956736e-09$ & $5.8822278e-09$\\ 
    \hline 
    $5$ $\dagger$   & 9.3132257e-10 & $0.17276611$ & $0.066240464$ & $0.066240464$ & $0.066240464$\\ 
    \hline 
\end{tabular}
\end{center}
\caption{Effectiveness of the method in computing rigorous enclosures of solutions in Mackey-Glass equation for parameters $n=[8,8]$, $\tau = [2,2]$, $\gamma=[1,1]$, $\beta=[2,2]$. Table shows statistics of coefficients of a given order computed over all grid points of the solution at a given time. Test setup was $\varepsilon = [0.0078125,0.0078125]$ (full step $h = \frac{\tau}{p}$ = [0.015625,0.015625]), $T = 24$ ($1536$ full steps or $12$ full delays)
 Note: superscript $\star$ means that diameter of coefficients at a grid point are presented (i.e. $j$ part of the f-set), where $\dagger$ means enclosures over intervals of length $h$ are presented ($\Xi$ part used). ,,No expand'' column contains data for the old algorithm, without representation expansion. ,,Expand $n$'' contains data for maximal order of the representation $2n$, 
,,Expand $n+1$'' contains data for maximal order of the representation $2n+1$, etc.}
\label{tab:benchmark}
\end{table}

%% file: sections/4_covrel.tex
\section{\label{sec:covrel}Topological tools}

In \cite{nasza-praca-focm} we have proven the existence of periodic orbits
(apparently stable) using the Schauder Fixed Point Theorem. Here we are interested
in a more general way to prove existence of particular solutions to DDEs with the use
of Poincar\'e maps generated with semiflow $\varphi$ of \eqref{eq:dde}. For this we will
recall the  concept of \emph{covering relations} from \cite{covrel-GiZ-1}, but
we will adopt it to the setting of infinite dimensional spaces and compact
mappings, similarly to a recent work \cite{chaos-kuramoto}. The main theoretical
tool to prove the existence of solutions, in particular the fixed points of continuous
and compact maps in $C^n_p$, will be the Leray-Schauder degree, which is an extension of
Fixed Point index (i.e. the local Brouwer degree of $Id - F$) to infinite dimensional Banach
spaces. We only recall the properties of the degree
that are relevant to our applications. For a broader description of the topic
together with the proofs of presented theorems we point out to
\cite{granas-book,degree-book} and references therein.
In particular, in what follows, we will use the notion of Absolute Neighbourhood Retract
(ANR) \cite{granas-book}. We do not introduce the formal definition
but we only note that (1) any Banach space is ANR and (2) any convex,
closed subset of a Banach space (or a finite sum of such) is an ANR
(Corollary 5.4 and Corollary 4.4 in \S 11. of \cite{granas-book}, respectively).

\subsection{Fixed Point Index for Compact Maps in ANRs}

Let $\X$ be a Banach space.
We recall that a continuous function $f : \X \supset V \to \X$ is a
\emph{compact map} iff $\overline{f(V)}$ is compact in $\X$.
With $Fix(f, U) = \{ x \in U : f(x) = x \}$ we denote
the set of fixed points of $f$ in $U$. Let now $X$ be
an ANR \cite{granas-book}, in particular $\X$ can be $X$,
and let $U$ be open subset of $X$, $f : \overline{U} \to X$.
Following \cite{granas-book}, by $\mathcal{K}\left(\overline{U}, X\right)$ we denote
the set of all compact maps $\overline{U} \to X$, and by
$\mathcal{K}_{\bd U}\left(\overline{U}, X\right)$ the set of all
maps $f \in \mathcal{K}(\overline{U}, X)$ that have no fixed points on $\bd U$,
$Fix(f, \bd U) = \emptyset$. We will denote
$Fix(f) = Fix(f, U) = Fix(f, \overline{U})$.
Let $V \subset \X$ be any set in the Banach space $\X$.
We say that a map $f : V$ is \emph{admissible} in $V$
iff $Fix(f, V)$ is a compact set. The following stronger assumption
that implies admissibility is often used in applications:
\begin{lemma}
\label{lem:fix-bd-empty}
Let $\X$ be a Banach space (can be infinite dimensional)
and $U \subset \X$ be an open set.
Assume $f : \overline{U} \to \X$ is a continuous, compact map.
If $f(x) \ne x$ for all $x \in \bd U$ then $f$ is admissible.
\end{lemma}
\textbf{Proof:} Let $F = (Id - f)^{-1}(\{0\})$ be the set of fixed points of $f$.
By assumption $f(x) \ne x$ on $\bd U$, we have $F \cap \bd U = \emptyset$ so $F \cap \bar{U} = F \cap U$.
The set $F$ is closed as a preimage of the closed set $\{0\}$ under continuous function $Id - f$, and so is $F \cap \overline{U}$.
Therefore $F \cap U$ is closed ant thus compact as
a subset of a compact set $\overline{f(\overline{U})}$:
$F \cap U = F \cap \overline{U} = f(F \cap \overline{U}) \subset \overline{f(\overline{U})}$.
\qed

By Lemma~\ref{lem:fix-bd-empty} we see that all
functions $f \in \mathcal{K}_{\bd U}(\overline{U}, X)$
are admissible, so that the Fixed Point Index is well
defined on them \cite{granas-book}:
\begin{theorem}[Theorem 6.2 in \cite{granas-book}]
\label{thm:fixed-point-index-and-properties}

Let $X$ be an ANR.
Then, there exists an integer-valued fixed point index
function $\iota(f, U) \in \Z$ (\emph{Leray-Schauder degree} of $Id - f$)
which is defined for all $U \subset X$ open and all
$f \in \mathcal{K}_{\bd U}\left(\overline{U}, X\right)$
with the following properties:
\begin{enumerate}[(I)]
\item
	\emph{(Normalization)} If $f$ is constant $f(x) = x_0$ then,
	$\iota(f, U) = 1$ iff $x_0 \in U$ and $\iota(f, U) = 0$ iff
	$x_0 \notin U$.
\item
	\emph{(Additivity)} If $Fix(f) \subset U_1 \cup U_2 \subset U$
	with $U_1, U_2$ open and $U_1 \cap U_2 = \emptyset$, then
	$\iota(f,U) = \iota(f,U_1) + \iota(f,U_2)$.
\item
	\emph{(Homotopy)} If $H : [0, 1] \times \overline{U} \to X$ is an
	\emph{admissible compact homotopy}, i.e. $H$ is continuous,
	$H_t = H(t, \cdot)$ is compact and admissible
	for all $t$, then $\iota(H_t) = \iota(H_0)$ for all $t \in [0, 1]$.
\item
	\emph{(Existence)} If $\iota(f, U) \ne 0$ then $Fix(f) \ne \emptyset$.
\item
	\emph{(Excision)} If $V \subset U$ is open, and $f$ has no fixed points in
	$U \setminus V$ then $\iota(f, U) = \iota(f, U \setminus V)$.
\item
	\emph{(Multiplicativity)} Assume $f_i : \X_i \supset X_i \supset \overline{U}_i \to X_i$,
	$i = 1, 2$ are admissible compact maps, and define
	$f(x_1, x_2) = (f_1(x_1), f_1(x_2)) \in X_1 \times X_1$ for
	$(x_1, x_2) \in U := \overline{U}_1 \times \overline{U}_2$.
	Then $f$ is a continuous, compact and admissible map with
	$\iota(f, U) = \iota(f_1, U_1) \cdot \iota(f_2, U_2)$.
\item
	\emph{(Commutativity)} Let $U_i \subset X_i \subset \X_i$, for $i =1, 2$
	be open and assume $f_i : U_1 \to X_2$,
	$g : U_2 \to \X_1$ and at least one of the maps $f, g$ is compact.
	Define $V_1 = U_1 \cap f^{-1}(U_2)$ and $V_2 = U_2 \cap g^{-1}(U_1)$,
	so that we have maps $g \circ f : \overline{V_1} \to X_1$ and
	$f \circ g : \overline{V_2} \to X_2$.
	
	Then $f \circ g$ and $g \circ f$ are compact and if
	$Fix(g \circ f) \subset V_1$ and $Fix(f \circ g) \subset V_2$ then
	\begin{equation*}
	\iota(g \circ f, V_1) = \iota(f \circ g, V_2).
	\end{equation*}
	
\end{enumerate}
\end{theorem}

For us, the key and the mostly used properties are the Existence,
Homotopy and Multiplicativity properties.
First one states that, if the fixed point index is non-zero, then there must
be a solution to the fixed-point problem $f(x) = x$ in the given set.
The Homotopy allows to relate the fixed point index $\iota(f, U)$
to some other, usually easier and better understood map,
for example $\iota(A, U)$, where $A$ is some linear function in
finite dimensional space. Normalization and Multiplicativity are used
to compute the fixed point index in the infinite dimensional ,,tail part''.

The following is a well-known fact:
\begin{lemma}
\label{lem:index-of-A}
Let $A : \R^n \to \R^n$ be a linear map. Then for any $U \subset \R^n$:
\begin{equation}
\iota(A, U) = \mathrm{sgn}\left(\det\right(Id - A\left)\right).
\end{equation}
\end{lemma}
Applying Commutativity property to $f = F \circ h^{-1}$ and $g = h$ gives:
\begin{lemma}
\label{lem:homeomorphism}
Let $F : U \to \X$ be admissible, continuous, compact map and let $h : \X \to \X'$
be a homeomorphism. Then $h \circ F \circ h^{-1} : \X' \supset h(U) = V \to \X'$
is admissible, and
\begin{equation*}
\iota(F, U) = \iota(h \circ F \circ h^{-1}, V).
\end{equation*}
\end{lemma}

\subsection{\label{sec:covrel-finite}Covering relations in $\R^d$}

In our application we will apply the fixed point index to detect periodic orbits of  some Poincar\'e maps
$P : C^{n}_p \supset U \to C^{n}_p$. We will introduce a concept of
\emph{covering relations}. A covering relation is
a way to describe that a given map $f$ \emph{stretches in a proper way} one set
over another. This notion was formalized in \cite{covrel-GiZ-1} for
finite dimensional spaces and recently extended to infinite spaces in
\cite{chaos-kuramoto} in the case of mappings between compact sets. In the
sequel we will modify this slightly for compact mappings between (not compact)
sets in the $C^n_p$ spaces.

To set the context and show possible applications,
we start with the basic definitions from \cite{covrel-GiZ-1}
in finite dimensional space $\R^d$, and then we will move to extend
the theory in case of $C^n_p$ spaces later in this section.

\begin{definition}[Definition~1 in \cite{covrel-GiZ-1}]
\emph{A h-set} N in $\R^{d_N}$ is an object consisting of the following data:
\begin{itemize}
\item $|N|$ - a compact subset of $\R^{d_N}$;
\item $u_N, s_N \in \N$ such that $u_N + s_N = d_N$;
\item a homeomorphism $c_N : \R^{d_N} \to \R^{d_N} = \R^{u_N} \times \R^{s_N}$ such that
\begin{equation*}
c_N(|N|) = \overline{\Ball_{u_N}}(0,1) \times \overline{\Ball_{s_N}}(0,1).
\end{equation*}
\end{itemize}
We set:
\begin{eqnarray*}
N_c & = & \overline{\Ball_{u_N}}(0,1) \times \overline{\Ball_{s_N}}(0,1) \\
N^-_c & = & \bd \overline{\Ball_{u_N}}(0,1) \times \overline{\Ball_{s_N}}(0,1) \\
N^+_c & = & \overline{\Ball_{u_N}}(0,1) \times \bd \overline{\Ball_{s_N}}(0,1) \\
N^- & = & c_N^{-1}(N^-_c) \\
N^+ & = & c_N^{-1}(N^+_c).
\end{eqnarray*}
\end{definition}
In another words, h-set $N$ is a product of two closed balls
in an appropriate coordinate system. The numbers $u_N$ and $s_N$
stands for the dimensions of exit (nominally unstable) and entry (nominally stable)
directions. We will usually drop the bars from the support
$|N|$ of the h-set, and use just $N$ (e.g. we will write $f(N)$ instead of $f(|N|)$.

The h-sets are just a way to organize the structure of a
support into nominally stable and unstable directions and to give
a way to express the exit set $N^-$ and the entry set $N^+$.
There is no dynamics here yet - until we introduce some
maps that stretch the h-sets across each other in a proper way.
\begin{definition}[Definition~2 in \cite{covrel-GiZ-1}]
\label{def:singlevalued-covering-FINITE}
Assume $N$, $M$ are h-sets, such that $u_N = u_M = u$. Let $P : |N| \to \R^{d_M}$ a continuous map.
We say that $N$ $P$-covers $M$, denoted by:
\begin{equation*}
N \cover{P} M
\end{equation*}
iff there exists continuous homotopy $H : [0,1] \times |N| \to R^{d_M}$
satisfying the following conditions:
\begin{itemize}
\item $H(0, \cdot) = P$;
\item $h([0,1], N^-) \cap M = \emptyset$;
\item $h([0,1], N) \cap M^+ = \emptyset$;
\item there exists a linear map $A : \R^{u} \to \R^{u}$ such that
\begin{eqnarray*}
H_c(1, (p, q)) & = & (A p, 0) \\
A(\bd \Ball_u(0,1)) & \subset & \R^u \setminus \overline{\Ball_u}(0,1)
\end{eqnarray*}
where $H_c(t, \cdot) = c_M \circ H(t, \cdot) \circ c_N^{-1}$ is
the homotopy expressed in good coordinates.
\end{itemize}
\end{definition}
A basic theorem  about  covering relations
is as follows:
\begin{theorem}[Simplified version of Theorem~4 in \cite{covrel-GiZ-1}]
\label{thm:periodic-covering-FINITE} 
Let $X_i \subset \R^d$ be h-sets and let
\begin{equation*}
X_1 \cover{P_1} X_2 \cover{P_2} \ldots \cover{P_k} X_{k+1} = X_1
\end{equation*}
be a covering relations chain. Then there exists $x \in X_1$ such that
\begin{eqnarray*}
x & \in & X_1 \\
(P_{r-1} \circ \ldots \circ P_{1}) (x) & \in & X_r \quad \textrm{ for } 2 \le r \le k, \\
(P_{k} \circ \ldots \circ P_{1}) (x) & = & x.
\end{eqnarray*}
\end{theorem}
Before we move on, we would like to point out what results can be
obtained using Theorem~\ref{thm:periodic-covering-FINITE}:
\begin{itemize}
\item
\textbf{Example 1.} Let $X \cover{P} X$, where $X$ is some h-set on a section
$S \subset \R^d$ and $P$ is a Poincare map $S \to S$ induced by the local flow
$\varphi$ of some ODE $x' = f(x)$. Then, there exists a periodic solution $x$
to this ODE, with initial value $x_0 \in X$. The parameter $u_X$ give the number of apparently
unstable directions for $P$ at $x$.

\item
\textbf{Example 2.}
Let $X_1$, $X_2$ be h-sets on a common section
$S \subset \R^d$, $X_1 \cap X_2 = \emptyset$, and assume $X_i \cover{P} X_j$ for
all $i, j \in \{1, 2\}$ where again $P$ is a Poincar\'e map $S \to S$ induced by
the semiflow $\varphi$ of some ODE. Then this ODE is chaotic in
the sense that there exists a countable many periodic solutions of arbitrary
basic period that visits $X_1$ and $X_2$ in any prescribed order. Also, there
exist non-periodic trajectories with the same property,
see for example \cite{covrel-GiZ-1,PZ-rossler-henon-chaos}.
\end{itemize}
In what follows, we will show the same construction can be done
under some additional assumptions in the infinite dimensional spaces.

\subsection{Covering relations in infinite dimensional spaces}

\begin{figure}
\center{
\includegraphics[width=12.5cm]{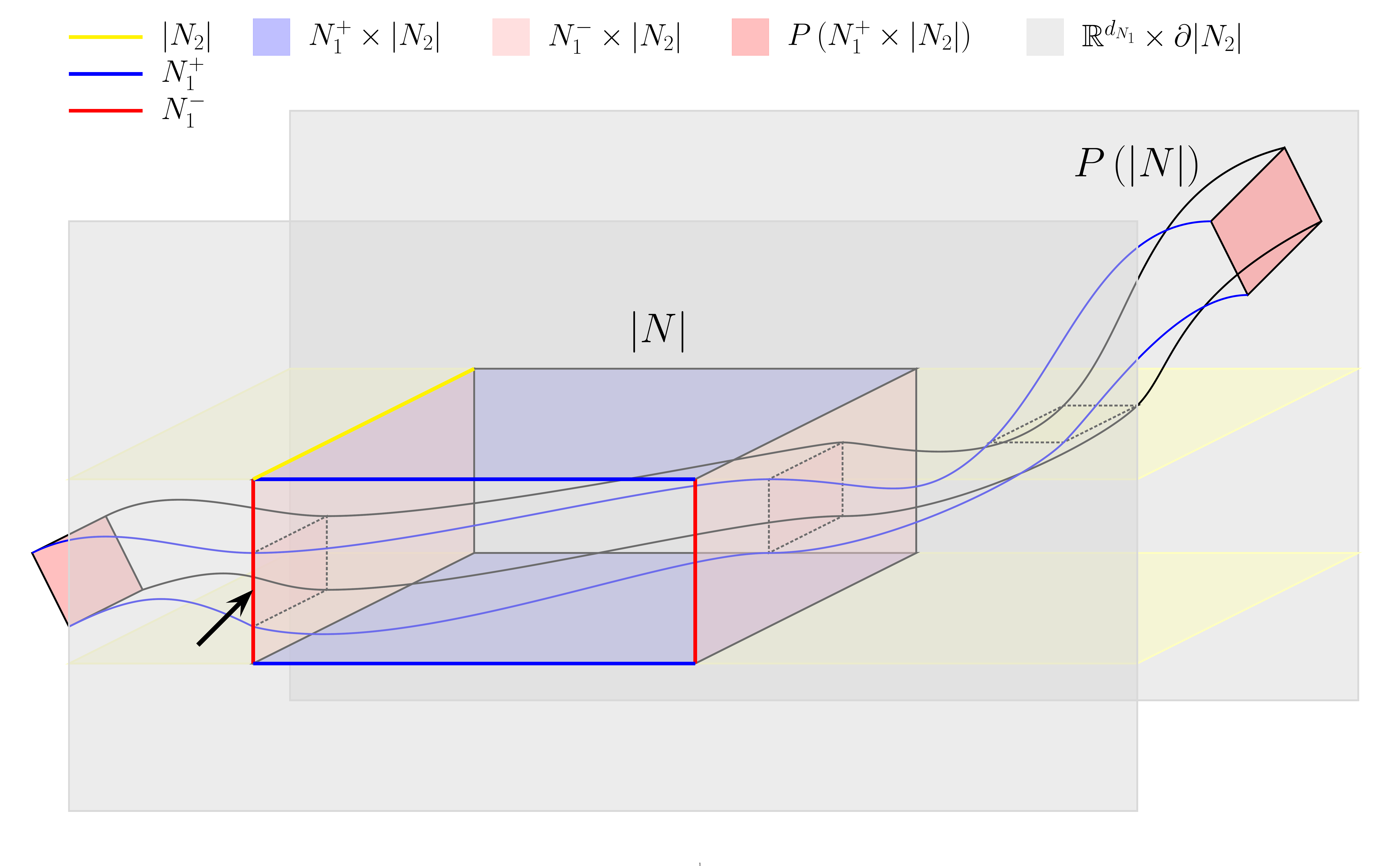}
}
\caption{\label{fig:covrel}An example of a covering relation $N \cover{P} N$
on an h-set with tail $N = (N_1, |N_2|)$, $u_{N_1} = 1, s_{N_1} = 1$. The tail $|N_2|$
is closed and convex in a potentially infinite dimensional space. The legend is as
follows: the set $|N|$ is the parallelepiped in the middle, whereas its image
$P(|N|)$ is stretched across $N$. The finite dimensional part is drawn in (x,y)-plane (width and height of the page), where the tail is drawn in z-coordinate (depth). 
The yellow thick line is one copy of the set
$|N_2|$ (the tail part), blue thick lines mark the set $N_1^+$ (the
,,entrance set'' of the finite dimensional part of $N$), red thick lines mark
the set $N_1^-$ (the ,,exit set'' of the finite dimensional part of $N$),
light-blue and light-red polygons mark the entrance set $N_1^+ \times |N_2|$ and
the exit set $N_1^- \times |N_2|$, respectively. The grey planes denote the
boundary of the strip $\R^{d_{N_1}} \times |N_2|$ - the image of $|N|$ under $P$
is forbidden to extend beyond those planes in z-coordinate due to the condition (C0). The set
$P(|N|)$ does not ,,touch'' the entrance set $N_1^+ \times |N_2|$ - condition
(C3) and the exit set $N_1^+ \times |N_2|$ is mapped outside $|N|$ (red
polytopes on left and right part of the picture) - condition (C2). Please note,
that the image $P(|N|)$ is allowed to touch the boundary $N_1 \times \bd |N_2|$
(place marked with a black arrow) as long as it does not go beyond the grey
planes. It is also allowed to bend in the stable direction of the finite
dimensional part outside the strip bounded by yellow hyperplanes (see right part of the picture). 
It is easy to
see, that the map $P$ can be homotopied, with a straight line homotopy
fulfilling condition (C0), to a map $(x, y, r) \mapsto (2\cdot x, 0, \bar{r})$,
where $\bar{r} \in |N_2|$ (up to the coordinate change $c_N$) - condition (C4).
}
\end{figure}

The crucial tool in proving Theorem~\ref{thm:periodic-covering-FINITE}
is the Fixed Point Index in finite dimensional spaces.
Therefore, similar results are expected to
be valid for maps and sets for which the infinite dimensional
analogue, namely Leray-Schauder degree of $Id - f$, exists. This was used
in \cite{chaos-kuramoto} for maps on compact sets in infinite dimensional
spaces. In this work we do not assume sets are compact, but we use the
assumption that the maps are compact -
the reasoning is almost the same. We will work on spaces
$\X = \X_1 \oplus \X_2$, where
$\X_1$ is finite dimensional (i.e. $\X_1 \equiv \R^M$) and
$\X_2$ will be infinite dimensional (sometimes refereed to as the tail).
In our applications, we will set
$\X = C^n_p = \R^{M(d,p,n)} \times (C^0([0,h], \R^d))^{d \cdot p}$,
with $\X_1 = \R^{M(d,p,n)}$.
We will use the following definitions that are slight modifications
of similar concepts from \cite{chaos-kuramoto}, where the tail
was assumed to be a compact set.
\begin{definition}
\label{def:h-set-with-tail}
Let $\X$ be a real Banach space.

An \emph{h-set with tail} is a pair $N = (N_1, |N_2|)$ where
\begin{itemize}
\item
	$N_1$ is an h-set in $\X_1$,
\item
	$|N_2| \subset \X_2$ is a closed, convex and bounded set.
\end{itemize}
Additionally, we set $u_N = u_{N_1}$, $|N| = |N_1| \times |N_2|$,
$c_N = (c_{N_1}, Id)$ and
\begin{eqnarray*}
N_c & = & c^{-1}_N\left(|N|\right) \ = \ N_{1,c} \times |N_2| \ = \\
    & = & \overline{\Ball_{u_{N_1}}}(0,1) \times \overline{\Ball_{s_{N_1}}}(0,1) \times |N_2|.
\end{eqnarray*}
\end{definition}
The tail in the definition refers to the part $|N_2|$.
We will just say that $N$ is an h-set when context is clear. Please
note that each h-set $N$ in $\R^d$ can be viewed as an h-set with tail,
where the tail is set as the trivial space $\R^{0} = \{ 0 \}$.
\begin{definition}
\label{def:covering-with-tail}
Let $\X$ be as in Def.~\ref{def:h-set-with-tail}. Let $N$, $M$
be h-sets with tails in $\X$ such that $u_N = u_M = u$.
Let $P : N \to \X$ be a continuous and compact mapping in $\X$.

We say that $N$ $P$-covers $M$
(denoted as before in Def.~\ref{def:singlevalued-covering-FINITE}
by $N \cover{P} M$), iff there exists continuous and compact homotopy
$H : [0,1] \times |N| \to \X$ satisfying the conditions:
\begin{itemize}
\item (C0) $H\left(t, |N|\right) \subset \R^{d_{M_1}} \times |M_2|$;
\item (C1) $H\left(0, \cdot\right) = P$;
\item (C2) $H\left(\left[0,1\right], N_1^- \times |N_2|\right) \cap M = \emptyset$;
\item (C3) $H\left(\left[0,1\right], |N|\right) \cap \left(M_1^+ \times |M_2|\right) = \emptyset$;
\item (C4) there exists a linear map $A : \R^{u} \to \R^{u}$ and
a point $\bar{r} \in M_2$ such that for all
$(p, q, r) \in N_c = \overline{\Ball_{u_{N_1}}}(0,1) \times \overline{\Ball_{s_{N_1}}}(0,1) \times |N_2|$
we have:
\begin{eqnarray*}
H_c(1, (p, q, r)) & = & (A p, 0, \bar{r}) \\
A(\bd \Ball_u(0,1)) & \subset & \R^u \setminus \overline{\Ball_u}(0,1)
\end{eqnarray*}
where again $H_c(t, \cdot) = c_M \circ H(t, \cdot) \circ c_N^{-1}$
is the homotopy expressed in good coordinates.
\end{itemize}
\end{definition}
Let us make some remarks on Definition~\ref{def:covering-with-tail}.
In contrary to \cite{chaos-kuramoto}, we do not assume that the
h-sets with tails $N$ and $M$ are compact in $\X$,
but we assume that the map $P$ is compact instead. However, the definition
in \cite{chaos-kuramoto} is a special case of Definition~\ref{def:covering-with-tail},
if we have $u_{N_1} = d_{N_1}$
and $|M_2|$ is a compact set. The additional structure of the finite dimensional
part $N_1$ we assume in Def.~\ref{def:covering-with-tail} allows for a more
general form of covering occurring in the finite dimensional part,
see Figure~\ref{fig:covrel}.

Now we will state theorems, similar to Theorem~\ref{thm:periodic-covering-FINITE},
that joins the sequences of covering relations to the real dynamics happening
in the underlying compact maps. We start with definitions:
\begin{definition}
Let $k > 0$ be fixed integer
and let $B$ be \emph{a transition matrix}: $B \in \matrices(k, k)$
such that $\Ball_{ij} \in \{0, 1\}$. Then define:
\begin{equation*}
\Sigma^+_B = \left\{ s \in \{1,\ldots,k\}^\N : \Ball_{s_i,s_{i+1}} = 1, \quad \forall i \in \N \right\}
\end{equation*}
and \emph{a shift function} $\sigma : \Sigma^+_B \to \Sigma^+_B$ by
\begin{equation*}
\sigma(s)_i = s_{i+1}.
\end{equation*}
The pair $(\Sigma^+_B, \sigma)$ is called \emph{a subshift of finite type
with transition matrix $B$}.
\end{definition}

\begin{definition}
Let $\mathcal{F}$ be a family of compact maps in a real Banach space $\X$.

We say that $\Gamma = ( \mathcal{N}, \mathcal{F}, Cov )$
is \emph{a set of covering relations on $\X$} iff
\begin{itemize}
\item
	$\mathcal{F}$ is a collection of continuous and compact maps on $\X$,
\item
	$\mathcal{N}$ is a collection of h-sets with tails
	$N_i \subset \X$, $i \in \{1,..,k\}$,
\item
	$Cov \subset \mathcal{N} \times \mathcal{F} \times \mathcal{N}$
	is a collection of covering relations, that is
	if $(N_i, P_l, N_j) \in Cov$ then $N_j \cover{P_l} N_j$.
\end{itemize}
A transition matrix $B \in \matrices(k, k)$ associated to $\Gamma$ is defined as:
\begin{equation}
B_{ij} =
\begin{cases}
	1 & \textrm{if there exists covering relation } N_i \cover{P_l} N_j \in Cov \\
	0 & \textrm{otherwise}.
\end{cases}
\end{equation}
\end{definition}

\begin{definition}
A sequence $(x_i)_{i \in \N}$ is called \emph{a full trajectory
with respect to family of maps} $\mathcal{F} = \{ f_i : 1 \le i \le m\}$
if for all $i \in \N$ there is $j(i) \in \{ 1,\ldots,m \}$ such that
$f_{j(i)}(x_i) = x_{i+1}$.
\end{definition}

Now we state two main theorems:
\begin{theorem}
\label{thm:periodic-covering-with-tail}
The claim of Theorem~\ref{thm:periodic-covering-FINITE} is true for
a covering relation chain where sets $X_i$ are h-sets with tail in a
real Banach space $\X$.
\end{theorem}

\begin{theorem}
\label{thm:symbolic-dynamics}
Let $\Gamma = ( \mathcal{N}, \mathcal{F}, Cov )$ be a set of covering
relations and let $B$ be its transition matrix.

Then, for every sequence of symbols $(\alpha_i )_{i \in \N} \in \Sigma^+_B$
there exist $(x_i)_{i \in \N}$ - a full trajectory with respect to
$\mathcal{F}$, such that $x_i \in X_{\alpha_i}$. Moreover, if
$(\alpha_i )_{i \in \N}$ is $T$-periodic, then the corresponding trajectory
may be chosen to be a $T$-periodic sequence too.
\end{theorem}

Before we do the proofs of Theorems~\ref{thm:periodic-covering-with-tail} and
\ref{thm:symbolic-dynamics}, we note that the examples of results that can be
obtained with covering relations on h-sets with tails are the same as
given before in Section~\ref{sec:covrel-finite} in the case of a finite-dimensional space $\R^d$.
In the context of DDEs we will use those theorems for h-sets with tails in the
form of a (p,n)-fset: $N = (N_1, |N_2|) = X(A, R) \subset C^n_p$. The natural decomposition 
is such that $\{ \xi \in \left(C([0, h], \R^d)\right)^p : [\xi] \subset R \} = |N_2|$ (the tail) 
and $N_1 = A \subset \R^{M(d,p,n)}$ (the finite-dimensional
part). In each application presented later in the paper we will
decide on $u_{N_1}$ and on the coordinates $c_{N_1}$ on the finite-dimensional part $A$.

\textbf{Proof of Theorem~\ref{thm:periodic-covering-with-tail}}:
We proceed in a  way, similar to the proof of Theorem~2 in \cite{chaos-kuramoto}.
To focus the attention and get rid of too many subscripts at once, we assume
without loss of generality that $c_{X_i} = Id$ for all $i$ and $X_i = N_i \times R_i$,
where $N_i \in \R^M$ is the finite-dimensional part.

Let now denote $X = X_1 \times \ldots \times X_k$, $N = N_1 \times \ldots \ N_k$
and $R = R_1 \times \ldots \ R_k$. Let also denote by
$\mathcal{Y} = \R^{M \cdot k} \times R$. With a slight abuse of notation we can
write $X \subset \mathcal{Y}$ and that $\mathcal{Y} \subset \X^k$. Since $\X^k$
is a Banach space (with the product maximum norm) so is $\mathcal{Y}$ with
topology inherited from the space $\X^k$. Moreover, we have $X \subset \Y$ with
$\interior_\Y X = \interior N_1 \times R_1 \times \ldots \times \interior N_k \times R_k$.
This will be important for proving that a fixed point problem we
are going to construct is solution-free on the boundary of $X$ in $\mathcal{Y}$.

We construct  zero finding problem:
\begin{equation}
\label{eq:fix-point-finding}
\begin{array}{rcl}
 P_k(x_k) & = & x_1 \\
 P_1(x_1) & = & x_2 \\
& \cdots & \\
P_{k-1}(x_{k-1}) & = & x_{k}, \\
\end{array}
\end{equation}
and we denote the left side of \eqref{eq:fix-point-finding}
by $F(x)$ and we are looking for a solution $x = F(x)$ with
$x = (x_1, x_2, \ldots, x_k) \in X$. With the already mentioned abuse of
notation, we can write $F(a, \xi) = (b, \zeta)$ for $a \in \R^{M \cdot k}$,
$\xi \in R$. In a similar way we construct a homotopy $H$, by pasting together
homotopies from the definition of h-sets with tails $X_i$:
\begin{equation*}
H(t, x) = \left(H_k\left(t, x_k\right), H_1\left(t, x_1\right), \ldots, H_{k-1}\left(t, x_{k-1}\right)\right)
\end{equation*}
It is obvious that $H(0,\cdot) = F$ and we will show that
$H(t, \cdot)$ is fixed point free (admissible) on the boundary
$\bd_{\Y}X$. Indeed, since
$\interior_\Y X = \interior N_1 \times R_1 \times \ldots \times \interior N_k \times R_k$
then for $(b, \zeta) \in \bd_\Y X$ there must be $i \in \{1,\ldots,k\}$ such
that $b_i \in \bd N_i = N_i^+ \cup N_i^-$. If $b_i \in N_i^-$ then (C2) gives
$H_i\left(t, (b_i,\zeta_i)\right) \notin X_{i+1}$ and consequently
$(\Ball_{i+1}, \zeta_{i+1}) \ne H\left(t, (b, \zeta)\right)_{i+1}$
(note, if $i = k$, the we set $i+1 = 1$).
If $b_i \in N_i^+$, then from (C3) it follows that
$H_{i-1}\left(t, (\Ball_{i-1}, \zeta_{i-1})\right) \notin \left(N_i^+ \times |R_i|\right)$
and so $H\left(t, (b, \zeta)\right)_i \ne (b_i, \zeta_i)$
(note, if $i = 1$, the we set $i-1 = k$). Therefore $H$ is admissible,
$H(t, x) \ne x$ for all $x \in \bd_\Y X$. Of course $H$ is also continuous
and compact.

Now, $\Y$ is an ANR (Corollary 4.4 in \S 11. of \cite{granas-book})
so fixed point index $\iota(H(t, \cdot), X)$ is well defined and constant
for all $t \in [0, 1]$. Applying Multiplicativity, Normalization (on the tail part)
and \ref{lem:index-of-A} on $H(1,\cdot)$ we get
$\iota(H(1, \cdot), X) = \Pi \iota\left(A_i, B^{u}(0,1)\right) = \pm 1$
(since $det(Id - A_i) \ne 0$ as $\|A_i\| > 1$ due to (C4)).

Finally, Existence property yields a fixed point $\bar{x}$ to $H(0, x) = F(x) = x$.
\qed

\textbf{Proof of Theorem~\ref{thm:symbolic-dynamics}} is almost the same as
of Theorem~3 in \cite{chaos-kuramoto}, with the exception that the
sets $X_i$ are not compact. This is overcome by considering
the convergence of sequences of points in the images $P_i(X_i)$,
which are pre-compact by the assumption on $P_i$'s.
\qed

We conclude with a lemma that allows to easily check
whether $N \cover{P} M$ in case $u_N = u_M = 1$. We will
check the assumptions of this lemma later in Section~\ref{sec:applications}, 
with the help of a computer.
\begin{lemma}
\label{lem:one-unstable}
For a h-set with tail $N$ let define: 
\begin{itemize}
\item $N_c^{l} = \{-1\} \times \Ball_{s} \times |N|$, $N^l = c_N^{-1}(N_c^{l})$ - the \emph{left edge of $N$}, and 
\item $N_c^{r} = \{1\} \times \Ball_{s} \times |N|$, $N^r = c_N^{-1}(N_c^{r})$ - the \emph{right edge of $N$}.
\end{itemize}

Let $\X$ be a Banach space, $X \subset \X$ be an ANR,
$N = (N_1, |N_2|)$, $M(M_1, |M_2|)$ be h-sets with tails in $X$
with $u_N = u_M = 1$ and $P : |N| \to X$ be a continuous and
compact map such that the following conditions apply
(with $P_c = c_M \circ P \circ c_N^{-1} : N_c \to M_c$):
\begin{enumerate}
\item
	(CC1) $\pi_{\X_2} P\left(|N|\right) \subset |M_2|$;
\item
	Either (CC2A)
	\begin{eqnarray*}
	P_c\left(N_c^l\right) \subset (-\infty, -1) \times \R_{s} \times |M_2|
	& and &
	P_c\left(N_c^r\right) \subset (1, \infty) \times \R_{s} \times |M_2|
	\end{eqnarray*}
	or (CC2B)
	\begin{eqnarray*}
	P_c\left(N_c^l\right) \subset (1, \infty) \times \R_{s} \times |M_2|
	& and &
	P_c\left(N_c^r\right) \subset (-\infty, -1) \times \R_{s} \times |M_2|
	\end{eqnarray*}	
\item
	(CC3)
	$P_c\left(N_c\right) \cap \left(\Ball_s \times |M_2|\right) = \emptyset$
\end{enumerate}
Then $N \cover{P} M$ with the homotopy given as
$H(t, \cdot) = (1-t) \cdot P + t \cdot (A, 0, \bar{r})$,
where $A : \R \to \R$ such that $Ax = 2x$ (CC2A) or $Ax = -2x$ (CC2B)
and $\bar{r}$ is any selected point in $|M_2|$.
\end{lemma}
\textbf{Proof:}
(C0) and (C1) from Definition~\ref{def:covering-with-tail} are
obviously satisfied. We also have (CC2) implies (C2) and (CC3) is the same as (C3).
Therefore, we only need to show (C4), that is, the image of the homotopy
computed on the set $\bd \Ball_u \times \overline{\Ball_s} \times |N_2|$ does
not touch the set $M_c$. This is obvious from the definition of $A$
in both cases (CC2A) and (CC2B).
\qed

Figure~\ref{fig:covrel} presents such a covering in case $u = s = 1$
and $N = M$. The easiest way to assure (CC1) and (CC3) is to assume
$P_c(N_c) \subset \R \times \Ball_s \times |M_2|$ - in fact
we check this in our computer assisted proofs presented in the
next section. 

%% file: sections/5_applications.tex
\section{\label{sec:applications}Applications}

In this section we present applications of the discussed algorithm
to two exemplary problems. First one is a computer assisted proof
of symbolic dynamics in a delay-perturbed R\"ossler system \cite{rossler}. 
The proof is done for two different choices of perturbations. The second application
consists of proofs of (apparently) unstable periodic orbits in the Mackey-Glass
equation for parameter values for which Mackey and Glass observed chaos in
their seminal paper \cite{mackey-glass}. 

Before we state the theorems, we would like to discuss presentation
of floating point numbers in the article. Due to the very nature of
the implementation of real numbers in current computers, numbers like
$0.1$ are \emph{not representable} \cite{ieee-754}, i.e. cannot be stored
in memory exactly.
On the other hand, many numbers representable on the computer could
not be presented in the text of the manuscript in a reasonable way,
unless we adopt not so convenient digital base-2 number representation.
However, the implementation IEEE-754 of the floating point numbers
on computers \cite{ieee-754} guarantees that, for any real number $x$
and its representation $\tilde{x}$ in a computer format, there is always a
number $|\epsi| \le \epsi_{machine}$ such that $\tilde{x} = x(1+\epsi)$.
The number $\epsi_{machine}$ defines the \emph{machine precision},
and, for the \texttt{double} precision \texttt{C++} floating-point numbers
that we use in the applications, it is of the order $10^{-16}$.
Finally, in our computations we use the \emph{interval arithmetic} to produce
rigorous estimates on the results of all basic operations such as
$+$, $-$, $\times$, $\div$, etc. In principle, we operate on intervals
$[a, b]$, where $a$ and $b$ are representable numbers, and the result
of an operation contains all possible results, adjusting
end points so that they are again representable numbers
(for a broader discussion on this topic,
see the work \cite{nasza-praca-focm} and references therein).
For a number $x \in \R$ we will write $[x]$ to denote the
interval containing $x$. If $x \in \Z$ then we have $[z] = [z,z]$,
as integer numbers (of reasonably big value) are representable in floating
point arithmetic.

Taking all that into account we use the following convention:
\begin{itemize}
\item whenever there is an explicit decimal fraction defined
in the text of the manuscript of the form
$d_1 d_2 \cdots d_k . q_1 q_2 \cdots q_m$
then that number appears in the computer implementation as
\begin{equation*}
[d_1 d_2 \cdots d_k q_1 q_2 \cdots q_k] \div [10^m],
\end{equation*}
where $\div$ is computed rigorously with the interval arithmetic.
For example, number $10^{-3} = 0.001$ appears in source codes
as \texttt{Interval(1.) / Interval(1000.)}.

\item whenever we present a \emph{result from the output of the
computer program} $x$ as a decimal number with non-zero fraction part, 
then we have in mind the fact that this represents some other number $y$
- the true value, such that $y = x(1+\epsi)$ with $|\epsi| \le \epsi_{machine}$.
This convention applies also to intervals: if we write interval $[a_1, a_2]$,
then there are some representable computer numbers $b_1$, $b_2$
which are true output of the program, so that
$b_i = a_i (1 +\epsi_i)$.

\item if we write a number in the following manner: $d_1 . d_2 \cdots d_k {}^{u_1 u_2 \cdots u_m}_{l_1 l_2 \cdots l_m}$ with digits $l_i, u_i, d_i \in  \{0,..,9\}$ then it represents
the following interval
\begin{equation*}
\left[ d_1 . d_2 \cdots d_k l_1 l_2 \cdots l_m, d_1 . d_2 \cdots d_k u_1 u_2 \cdots u_m \right].
\end{equation*}
For example $12.3_{456}^{789}$ represents the interval $[12.3456, 12.3789]$
(here we also understand the numbers taking into account the first two conventions).
\end{itemize}

The last comment concerns the choice of various parameters for the proof,
namely, the parameters of the space $C^n_p$ and the initial sets around
the numerically found approximations of the dynamical phenomena under consideration.
The later strongly depends on the investigated phenomena, so we will
discuss general strategy in each of the following sections, whereas
the technical details are presented in Appendices~\ref{app:lohner} and \ref{app:data}.

The choice of parameters $n$ and $p$ corresponds basically to the choice 
of the order of the numerical method and a fixed step size $h = \frac{\tau}{p}$, respectively.

Usually, in computer assisted proofs, we want $n$ to be high, so that the
local errors are very small. In the usual case of ODEs with $f \in C^\infty$
we can use almost any order, and it is easy for example to set $n = 40$. However,
in the context of $C^n_p$ spaces and constructing Poincar\'e maps for DDEs, 
we are constrained with the \emph{long enough time} $T = (n+1) \cdot \tau$ 
(Definition~\ref{def:long-enough-time}) to obtain well defined maps. 
Therefore, the choice of $n$ corresponds usually to the return time
to section $t_P$ for a given Poincar\'e map, satisfying 
$t_P(X_0) > (n+1) \cdot \tau$, for some set of initial data $X_0 \subset C^n_p$.

The choice of the step size $h$ is more involved. It should not be too small, 
to reduce the computational time and cumulative impact of all local errors after 
many iterations, and not so big, as to 
effectively reduce the size of the local error. Also, the dynamics of the system 
(e.g. stiff systems) can impact the size of the step size $h$.  
In the standard ODE setting, there are strategies to set the step size
dynamically, from step to step, e.g. \cite{book-adaptative-steps}, but in the 
setting of our algorithm for DDEs, 
due to the continuity issues described in Section~\ref{sec:integrator}, we must 
stick to the fixed step size $h = \frac{\tau}{p}$. The step size must be also smaller 
than the (apparent) radius of convergence of the forward Taylor representation of the 
solution at each subinterval, but this is 
rarely an issue in comparison to other factors, e.g. the local error estimates.
In our applications we chose $p = 2^{m}$ for a fixed $m \in \N$, so that the 
grid points are \emph{representable floating point numbers} 
(but the implementation can work for any $p$). 

We also need to account for 
the memory and computing power resources. For $d$-dimensional systems \eqref{eq:dde}, and 
with $n$, $p$ fixed, we have that the representation of a Lohner-type set $A = x_0 + C \cdot r_0 + E$ in phase-space of $\varphi$,
where $C \in \matrices{M, M}$, requires 
at least $O(M^2)$, with $M = O(d \cdot n \cdot p)$. Then, doing one step of the full step algorithm
is of $O(d^2 \cdot n^2 \cdot M)$ computational complexity. Due to the long enough time integration, computation
of a single orbit takes usually $O(n \cdot p)$ steps, and we get the computational complexity
of computing image $P(X)$ for a single set $X$ of $O(d \cdot n^2 \cdot d \cdot p \cdot n \cdot M) = O(d \cdot n^2 \cdot M \cdot M) = O(M^2)$ 
(if we assume $n, d << M$). Therefore, we want to keep 
$M^2 = (d \cdot n \cdot p)^2$ of reasonable size, both because of time and
memory constraints. Our choice here is $M \le 10^3$.

\subsection{Symbolic dynamics in a delay-perturbed R\"ossler system}

In the first application, we use R\"ossler ODE of the form \cite{rossler}:
\begin{eqnarray}
x' &=& -(y + z) \notag \\
y' &=& x + ay \label{eq:rossler} \\
z' &=& b + z(x-c).
\end{eqnarray}
In what follows we will
denote r.h.s. of \eqref{eq:rossler} by $f$ and by $v \in \R^3$ we
denote vector $v = (x, y, z)$. By $\pi_x$ we denote projection onto
$x$ coordinate, similarly for $\pi_y, \pi_z$.

\begin{figure}
\center{
\includegraphics[width=8.0cm]{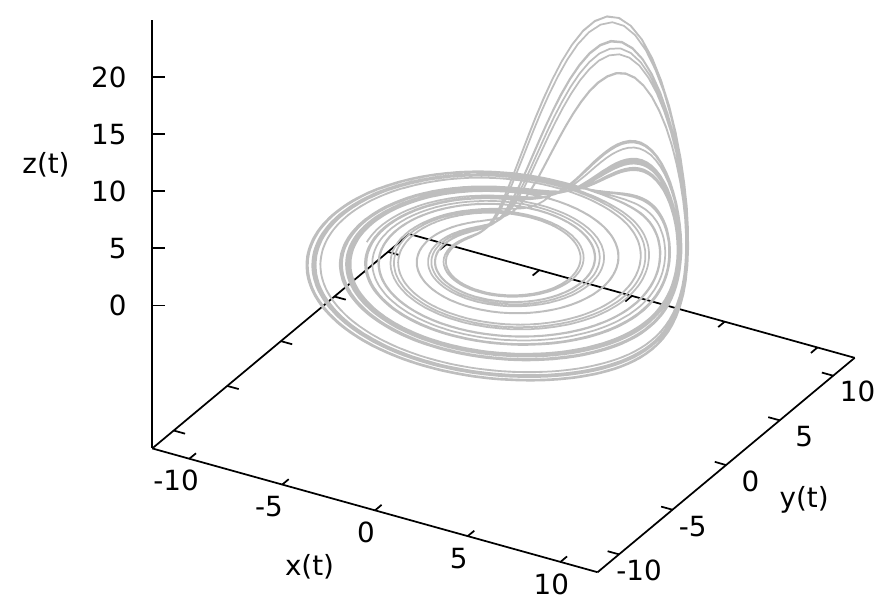}
}
\caption{\label{fig:attractor-rossler}Numerically observed
attractor in the R\"ossler ODE for classical values of parameters:
$a = b = 0.2$, $c = 5.7$. Picture generated by integrating forward
in time single trajectory for a long time.}
\end{figure}

We set the classical value of parameters $a = b = 0.2$,
$c = 5.7$ \cite{rossler}.  For those parameter values, 
an evidence of a strange attractor
was first observed numerically in \cite{rossler} ,
see Fig.~\ref{fig:attractor-rossler}. In \cite{PZ-rossler-henon-chaos}, 
it was proved by computer assisted argument that there is
a subset of the attractor which exhibit symbolic dynamics. 
A more recent results for R\"ossler system can also be found
in \cite{rossler-sharkovskii} (Sharkovskii's theorem) and 
the methodologies there should be easily adaptable in the 
context of delay perturbed systems presented in this paper. 

We are going to study a delayed perturbation of the R\"ossler system
\eqref{eq:rossler} of the following form:
\begin{equation}
\label{eq:rossler-perurbed}
v'(t) = f(v(t)) + \epsilon \cdot g(v(t - 1)),
\end{equation}
where parameter $\epsilon$ is small. We consider two toy examples:
first, where $g = f$ and the second one where $g$
is given explicitly as
\begin{equation}
\label{eq:perturbation-sin}
g(x, y, z) = \left(\sin(x \cdot y), \sin(y \cdot z), \sin(x \cdot z)\right).
\end{equation}
We expect that for any bounded $g$ there should be a sufficiently
small $\epsilon$ \cite{nasza-epsilon-perturb-dde} so that the dynamics 
of the perturbed system is preserved. However, in this work, we study 
explicitly given value for $\epsilon$. 

\begin{rem}
The source codes of the proof are generic. 
The interested reader can experiment with other forms 
of the perturbation by just changing the definition 
of the function $g$ in the source codes of the 
example. 
\end{rem}

\begin{figure}
\center{
\includegraphics[width=4.8cm]{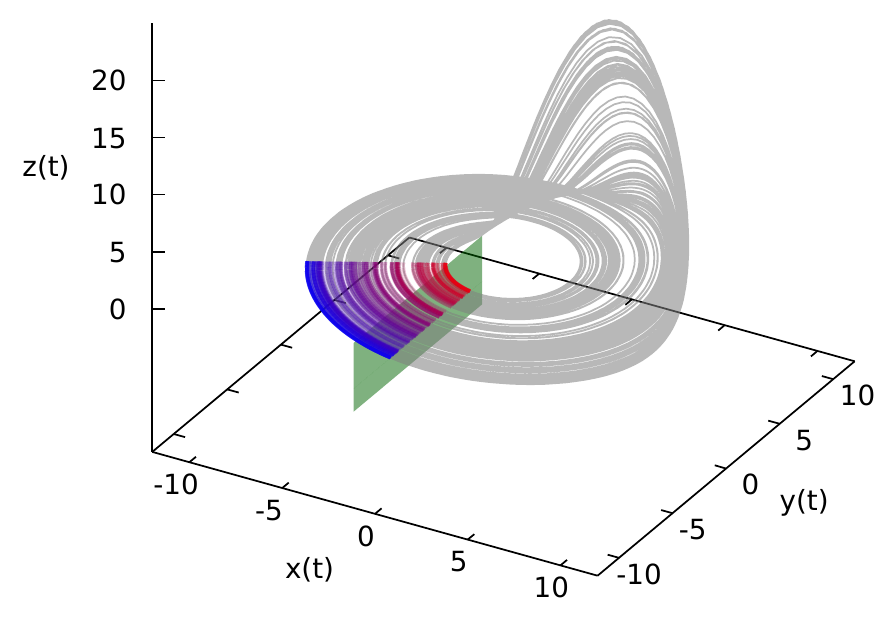}
\includegraphics[width=4.8cm]{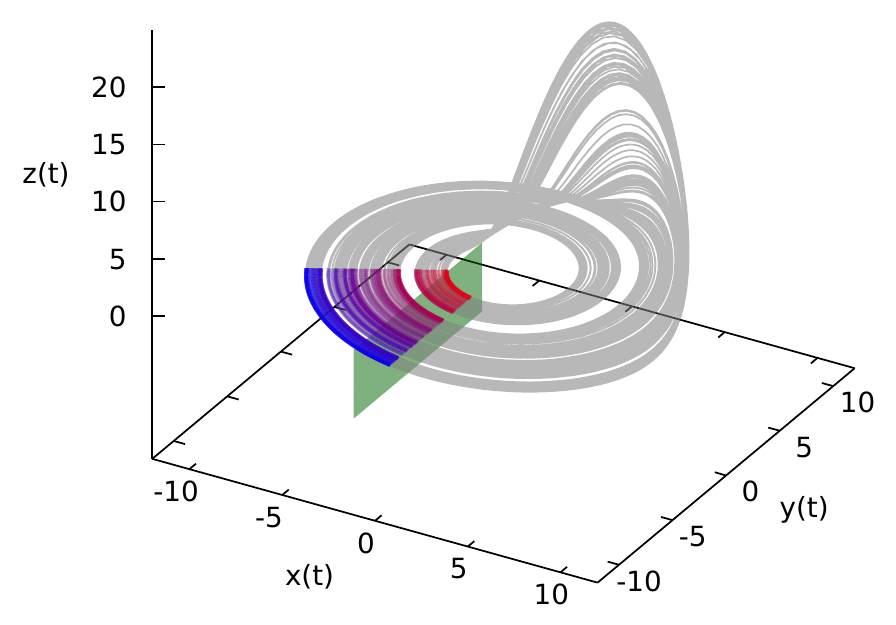}
\includegraphics[width=4.8cm]{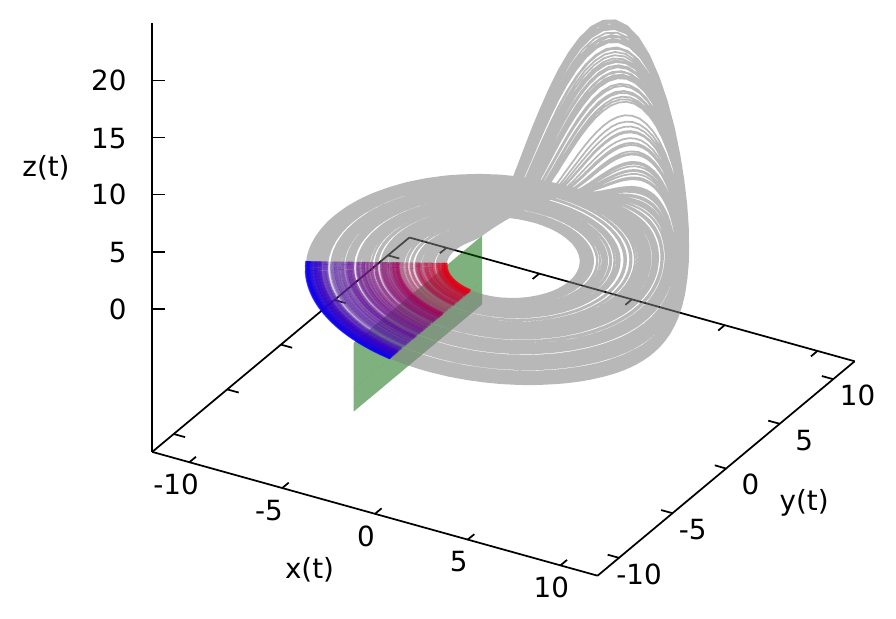}
}
\caption{\label{fig:attractor-rossler-delayed} The numerically observed
attractors for the system studied in Theorem~\ref{thm:rossler-delayed}.
The cases (a)-(c) are shown from left to right, respectively.
The grey attractor is the very long trajectory
$v(t)$ obtained for a single constant initial function. The section
$S_0$, represented as a green rectangle on the picture, spans in fact
across the space $\R^3$, as can bee seen
by the red to blue region that shows the segments of $v$
which lie on the section $S_0$,
i.e. the set $\{ v_t : \pi_x(v(t)) = 0 \}$.
The colours are assigned with ascending $\pi_y v(t)$ value.
Those segments are used to define the $W_u$ coordinate
in the set $X(A, \Xi)$.
}
\end{figure}

We will be studying the properties of a Poincar\'e map defined
on the section $S_0 \subset C^n_p$ given by:
\begin{equation*}
S_0 = \{ v \in C^n_p : \pi_x (v(0)) = 0 \}.
\end{equation*}
The section $S_0$ in an extension to $C^n_p$ of the section
$S = \{ v \in \R^3 : \pi_x v = 0 \} \subset \R^3$ used in the proofs
in \cite{PZ-rossler-henon-chaos}. The section $S$ is drawn in green 
in Fig~\ref{fig:attractor-rossler-delayed}, whereas the projection
of the attractor onto section $S_0$ is drawn as a blue-red gradient
(the solution segments $v$ with $\pi_x v(0) = 0$). 

In what follows, we set the parameters for the space $C^n_p$ to $p = 32$ and $n = 3$.
We prove, with the computer assistance, the following theorems:
\begin{theorem}
\label{thm:rossler-delayed} For parameter values $a = b = 0.2$, $c = 5.7$ 
in \eqref{eq:rossler} there exists sets
$X_A = X(A, \Xi), X_1 = X(N_1, \Xi), X_2 = X(N_2, \Xi) \subset S_0$
with explicitly given $A, N_1, N_2$ and $\Xi$, such that for
the system~\eqref{eq:rossler-perurbed} with $\epsilon = 10^{-3}$ and 
perturbations: (a) $g \equiv 0$ - original system treated as a DDE, 
(b) $g = f$ and (c) $g$ given as in in Eq~\eqref{eq:perturbation-sin} 
we have the following:

\begin{enumerate}
\item $P(X(A, \Xi)) \subset X(A, \Xi)$ and, in consequence, there
exists a non-empty invariant set in $X(A, \Xi)$ for the map $P : S_0 \to S_0$.

\item the invariant set $I = Inv(P^2, X_1 \cup X_2)$ of $ X_1 \cup X_2$
under the map $P^2$ on $I$ is non-empty and the dynamics of $P^2$ is
conjugated to the shift on two symbols ($\sigma : \Sigma_2 \to \Sigma_2$,
$\sigma(e_k) = e_{k+1}$), i.e. if we denote by $g : I \to \Sigma_2$
the function $g(x)_k = i \iff P^{2k}(x) \in X_i$, then we have
$g \circ P^2|_I = \sigma \circ g$.
\end{enumerate}
\end{theorem}

Before we present the proof(s), we would like to make
a remark on the presentation of the data from the
computer assisted part:
\begin{rem}[Convention used in the proofs]
\label{rem:description-sets}
The proofs of those theorems are computer assisted and the parameters
of the phase-space $C^n_p$ of representations are $d = 3$, $p = 32$, $n=3$,
giving in total the dimension of the finite dimensional part of
$M(d, p, n) = d \cdot (1 + p \cdot (n+1)) = 387$. Therefore it is not
convenient to present complete data of the proofs in the manuscript. Instead, we
assume the sets are explicitly given in the following forms (and the
interested reader is refereed to Appendix~\ref{app:data}
for the details on how they are constructed):
\begin{align*}
X_A = X(A, R): \quad & A = v_{ref} + C \cdot \{ 0\} \times W_u \times \Ball^{\|\cdot\|_\infty}_{M-2}(0, 1) \\
X_i = X(N_i, R): \quad & N_i = v_{ref} + C \cdot \{ 0\} \times W_i \times \Ball^{\|\cdot\|_\infty}_{M-2}(0, 1) \\
\Xi = \Ball^{\|\cdot\|_\infty}_{d \cdot p}(0, 1)  \quad & &
\end{align*}
with $v_{ref} \in S_0$, 
$W_u, W_1, W_2$ closed intervals such that $W_1 \cap W_2 = \emptyset$ and
$W_i \subset W_u \subset \R$, and we remind $\Ball^{\|\cdot\|_\infty}_{D}(0, 1)$
denotes the unit radius ball in the $\max$ norm in $\R^{D}$ centred at $0$.
Note, this description of sets makes it clear they are h-sets with tails
on $S_0$ (up to the scaling of nominally unstable direction $W$), where
$u = 1$ and $s_A = s_{N_i} = s = M(d, p,n)-2$, the support set
$|A| = \{ 0\} \times W_u \times \Ball^{\|\cdot\|_\infty}_{M-2}(0, 1)$ and
the affine coordinate change $c_A(\cdot) = v_{ref} + C(\cdot)$
with inverse change $c_A^{-1}(\cdot) = C^{-1}(\cdot - v_{ref})$.
Now, the computation of any Poincar\'e map $P : X_A \to S_0$ for
the initial data $X(A, \Xi)$ produces set $X(B, \Omega) = P(X(A, \Xi))$
and there exist sets
\begin{eqnarray*}
c_A^{-1}(B) &=\ B_c & \subset \ \{ 0 \} \times (B_c)_2 \times \Ball^{\|\cdot\|_\infty}_{M-2}(0, r_{B}) \\
&\ \Omega & \subset \ \Ball^{\|\cdot\|_\infty}_{d \cdot p}(0, r_{\Omega})
\end{eqnarray*}
for some $r_{B}, r_{\Omega} \in \R_{+}$. This allows to describe the geometry
of $X(A, \Xi)$ and (estimates on) $P(X(A, \Xi))$ by just a couple of numbers:
$W_u$, $\pi_2 B_c$ (the size of set $B$ in the nominally unstable direction), 
$r_{B}$ (upper bound on all coefficients in the finite nominally stable part) and 
$r_{\Omega}$ (upper bound on all $\xi$ in the tail part), which are suitable 
for a concise presentation in the manuscript. 
\end{rem}
The sets used in the computations are obtained 
by computing the appropriately enlarged enclosure on the set 
of segments of solutions to the unperturbed ODE~\eqref{eq:rossler}. 
We choose a set $\tilde{A} \subset \R^3$ such that 
$\tilde{A} \in \{ v \in \R^3 \pi_x v = 0\}$ is a trapping region for 
the Poincar\'e map of the unperturbed ODE: $P(\tilde{A}) \subset \tilde{A}$. 
Then we choose a set $X(A, \Xi)$ to contain the segments 
of $\tilde{A}$ propagated back in time for a full delay 
with the unperturbed ODE: 
\begin{equation*}
\left\{ v:[-1, 0] \to \R : v(0) \in A, v(s) = \varphi_0(s, v(0)) \right \} \subset A,
\end{equation*}
where $\varphi_0$ is the flow in $\R^3$ for \eqref{eq:rossler}.
Detailed procedure how the set $A$ was generated is described in 
the Appendix~\ref{app:data}.
The set $\tilde{A}$ 
was chosen to be $\{0 \} \times [-10.7, -2.2] \times [0.021, 0.041]$,
whereas the sets $\tilde{N}_1 = [-8.4, -7.6]$ and $\tilde{N}_2 = [-5.7, -4.6]$.
Finally, the orbit $v_0$ with $\pi_2 v_0(0) = -6.8$ is selected 
among the orbits in the attractor as the reference point of the sets $X_A, X_1, X_2$. 
The set $W_u$ is chosen as $W_u = \pi_2 A_c = \pi_y \tilde{A} - \pi_2 v_0(0) = [-3.9,4.6]$. 
The same is true for sets $N_1, N_2$, with $W_1 = [-1.6, -0.8]$, $W_2 = [1.1, 2.2]$. 

Now we can proceed to the proofs.

\vspace{1em}

\textbf{Proof o Theorem~\ref{thm:rossler-delayed}}
The proofs for parts (a), (b), and (c) follow the same methodology, 
therefore we present the details only for case (a) and then, only the estimates
from the other two cases. In principle, we will show that $P(X_A) \subset X_A$
and $X_i \cover{P^2} X_j$ for all $i, j \in \{ 1, 2 \}$ and then apply
Theorem~\ref{thm:symbolic-dynamics}.

The set $X(A,R)$ and two other sets are given as described
in Remark~\ref{rem:description-sets}. The computer
programs for the proof are stored in \verb#./examples/rossler_delay_zero#.
The data for which presented values were computed
is stored in \Verb#./data/rossler_chaos/epsi_0.001#.
See Appendix~\ref{app:data} for more information. Additionally
to the estimates presented below, the computer programs
verify that $t_P(x) > (n+1)$ (i.e. long enough for Poincar\'e maps
to be well defined) and that the function $t_P(\cdot)$ is well
defined. For details, see the previous work \cite{nasza-praca-focm}.

First, we prove that $P_c(X(A, \Xi)) \subset (A_c, \Xi)$. 
Let $(B_c, \Omega)$ will be output of the rigorous
program \Verb#rig_prove_trapping_region_exists# run for 
the system in case (a) such that 
$P_c(X(A, \Xi)) \subset (B_c, \Omega)$. It suffices to show 
the following:
\begin{itemize}
\item $\pi_2 P_c(X(A, \Xi)) = \pi_2 B_c \subset W_u = \pi_2 A_c$;
\item $\pi_i P_c(X(A, \Xi)) = \pi_i B_c < 1$ for all $i > 2$;
\item $\pi_{\Xi_i} P_c(X(A, \Xi)) = \pi_i \Omega < 1$ for all $i \in \{1, \ldots, p \cdot d \}$.
\end{itemize}
Indeed, we have:
\begin{itemize}
\item $\pi_2 P_c(X(A, \Xi)) = [-3.786230021035,3.92103823500285] \subset [-3.9,4.6] = W_u$; 
\item $\pi_i P_c(X(A, \Xi)) \le 0.910355124006778 < 1$, for $i > 2$; 
\item $\pi_{\Xi_i} P_c(X(A, \Xi)) \le 0.395102819146026 < 1$ for all $i$. 
\end{itemize}
Which finishes the proof of the first assertion.

For the second assertion we prove that we have a set of
full covering relations:
\begin{equation*}
X_i \cover{P^2} X_j, \quad i, j \in \{ 1, 2\}.
\end{equation*}
We remind that the sets $N_{i,c} = \{ 0 \} \times [W_{i}^{l}, W_{i}^r] \times \Ball^{\|\cdot\|_\infty}_{M-2}(0, 1)$
with $W_1 = [-1.6, -0.8]$, $W_2 = [1.1, 2.2]$.
The program \Verb#./rig_prove_covering_relations# produces the following inequalities: 
\begin{itemize}
\item (L1-L1) $\pi_2 P^2_c(X(N_{1}^{l}, \Xi)) = -1._{696238902429803}^{708946819732338} < -1.6 = \pi_2 N_{1,c}^{l} < \pi_2 N_{2,c}^{l}$
\item (R1-R2) $\pi_2 P^2_c(X(N_{1}^{r}, \Xi)) = 2.4_{09511664184434}^{17718805618395} > 2.2 = \pi_2 N_{2,c}^{r} > \pi_2 N_{1,c}^{r}$
\item (R2-L1) $\pi_2 P^2_c(X(N_{2}^{r}, \Xi)) = -1.83_{8887194518363}^{9215629292839} < -1.6 = \pi_2 N_{1,c}^{l} < \pi_2 N_{2,c}^{l}$
\item (L2-R2) $\pi_2 P^2_c(X(N_{2}^{l}, \Xi)) = 2.2_{69015891346912}^{70120359885664} > 2.2 = \pi_2 N_{2,c}^{r} > \pi_2 N_{1,c}^{l}$,
\end{itemize}
where sets $N^l$, $N^r$ etc. are defined as in Lemma~\ref{lem:one-unstable}.
It is ease to see that, those inequalities, together
with the existence of trapping region $X_A$, imply
that for each $i, j \in \{1, 2\}$ conditions (CC1)-(CC3)
in Lemma~\ref{lem:one-unstable} are satisfied, that is
$X_i \cover{P^2} X_j$, which finishes the proof
for the case (a) after applying Theorem~\ref{thm:symbolic-dynamics}.

For the cases (b) and (c) we only present estimates:

\begin{itemize}
\item
Case (b), $g = f$. Output from \Verb#rig_prove_trapping_region_exists# is:
\begin{itemize}
\item $\pi_2 P_c(X(A, \Xi)) = [-3.82791635121864,3.90123013871349] \subset [-3.9,4.6] = W_u$;
\item $\pi_i P_c(X(A, \Xi)) \le 0.960537051554584 < 1$, for $i > 2$; 
\item $\pi_{\Xi_i} P_c(X(A, \Xi))  \le 0.397264977921163 < 1 = r_\Xi$, for all $i$.
\end{itemize}
Output from program \Verb#./rig_prove_covering_relations# is:
\begin{itemize}
\item (L1-L1) $\pi_2 P^2_c(X(N_{1}^{l}, \Xi)) = -1.6_{68486957556001}^{84410417326001} < -1.6 = \pi_2 N_{1,c}^{l} < \pi_2 N_{2,c}^{l}$
\item (R1-R2) $\pi_2 P^2_c(X(N_{1}^{r}, \Xi)) = 2.4_{64065036803807}^{74268236696726} > 2.2 = \pi_2 N_{2,c}^{r} > \pi_2 N_{1,c}^{r}$
\item (R2-L1) $\pi_2 P^2_c(X(N_{2}^{r}, \Xi)) = -1.76_{7206286440370}^{9151140189891} < -1.6 = \pi_2 N_{1,c}^{l} < \pi_2 N_{2,c}^{l}$
\item (L2-R2) $\pi_2 P^2_c(X(N_{2}^{l}, \Xi)) = 2.36_{0282881761384}^{2685243092644} > 2.2 = \pi_2 N_{2,c}^{r} > \pi_2 N_{1,c}^{r}$
\end{itemize}

\item
Case (c), $g$ as in \eqref{eq:perturbation-sin}. Output from \Verb#rig_prove_trapping_region_exists# is:
\begin{itemize}
\item $\pi_2 P_c(X(A, \Xi)) = [-3.78710970137727,3.92188126709857] \subset [-3.9,4.6] = W_u$;
\item $\pi_i P_c(X(A, \Xi)) \le 0.951680057117636 < 1$, for $i > 2$; 
\item $\pi_{\Xi_i} P_c(X(A, \Xi)) \le 0.459753301095895 < 1$, for all $i$.
\end{itemize}
Output from program \Verb#./rig_prove_covering_relations# is:
\begin{itemize}
\item (L1-L1) $\pi_2 P^2_c(X(N_{1}^{l}, \Xi)) = -1._{695427259804897}^{714200213156898} < -1.6 = \pi_2 N_{1,c}^{l} < \pi_2 N_{2,c}^{l}$
\item (R1-R2) $\pi_2 P^2_c(X(N_{1}^{r}, \Xi)) = 2.4_{08774107762390}^{20396855111791} > 2.2 = \pi_2 N_{2,c}^{r} > \pi_2 N_{1,c}^{r}$
\item (R2-L1) $\pi_2 P^2_c(X(N_{2}^{r}, \Xi)) = -1.8_{38300180457653}^{41157932729915} < -1.6 = \pi_2 N_{1,c}^{l} < \pi_2 N_{2,c}^{l}$
\item (L2-R2) $\pi_2 P^2_c(X(N_{2}^{l}, \Xi)) = 2.2_{67377344403297}^{70144525622461} > 2.2 = \pi_2 N_{2,c}^{r} > \pi_2 N_{1,c}^{r}$
\end{itemize}
\end{itemize}
\qed

\begin{figure}
\center{
\includegraphics[width=6.5cm]{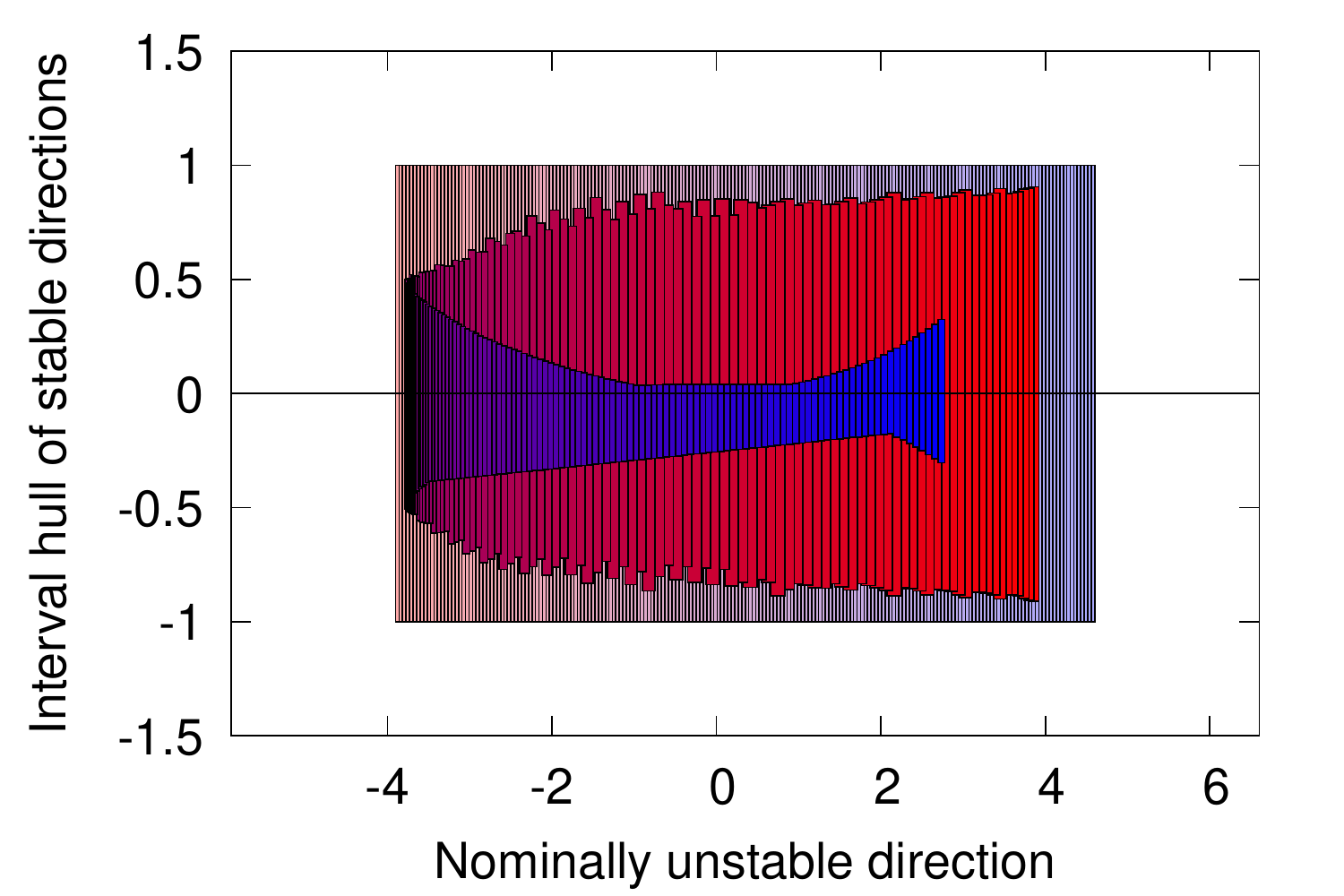}
\includegraphics[width=6.5cm]{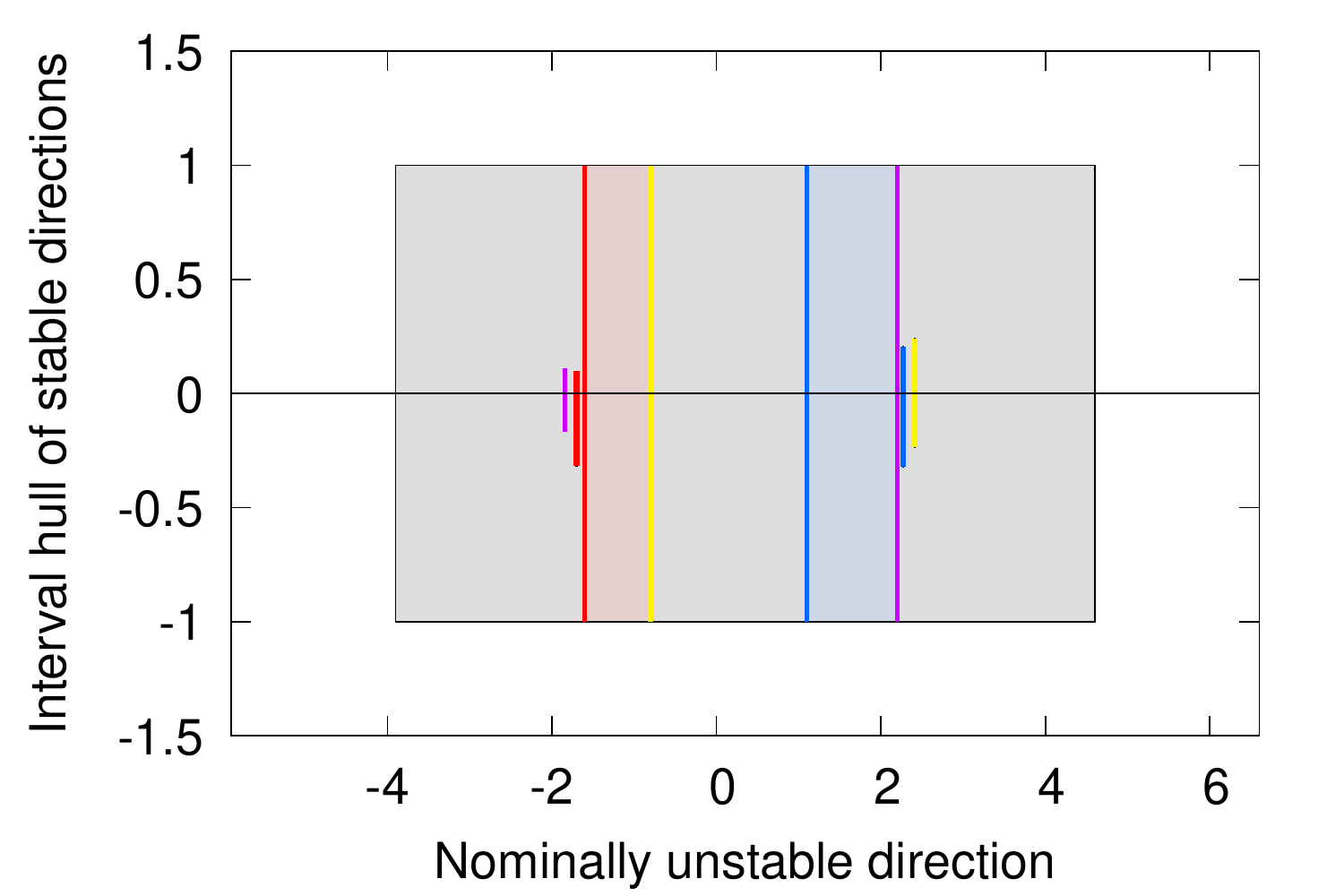}

\includegraphics[width=6.5cm]{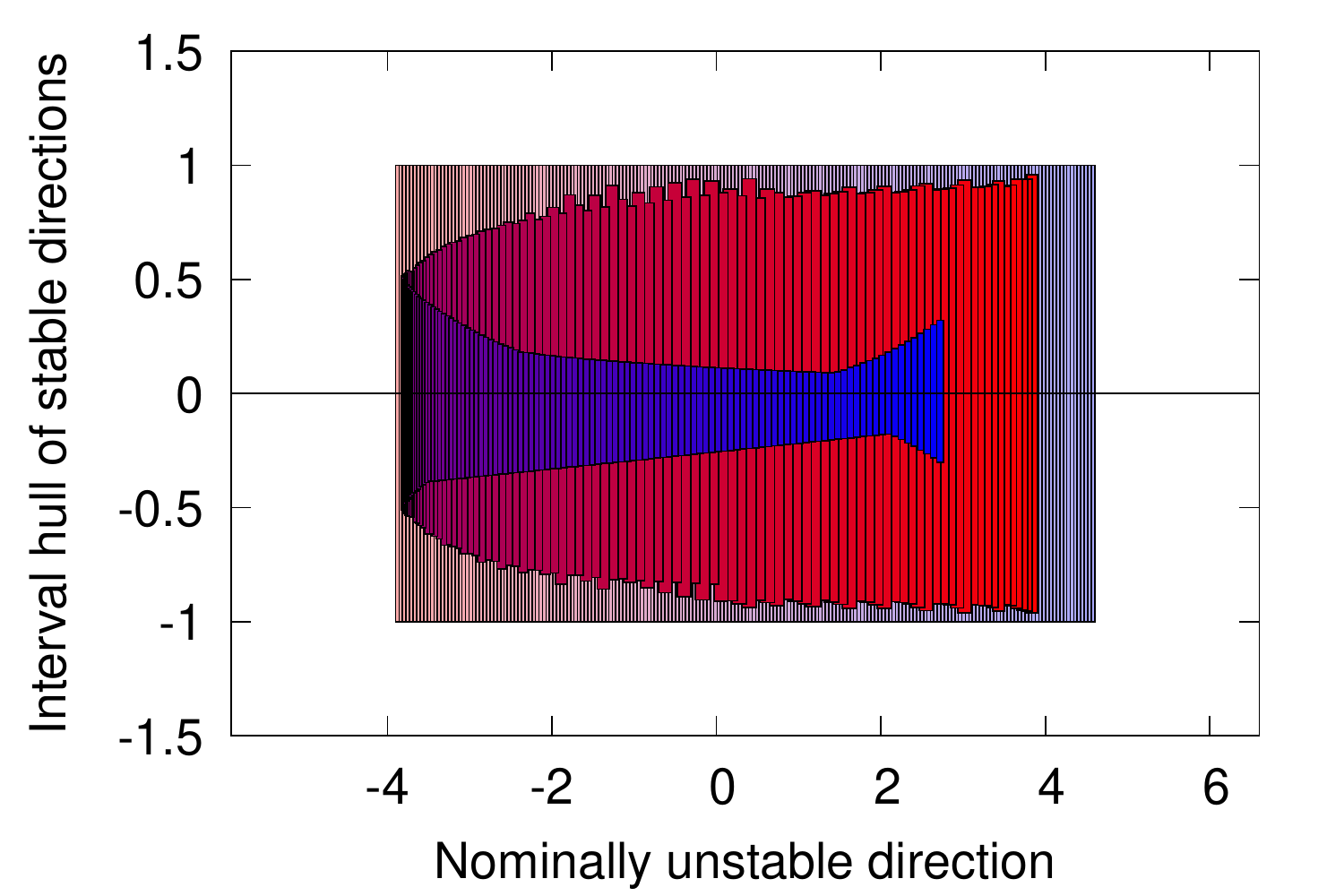}
\includegraphics[width=6.5cm]{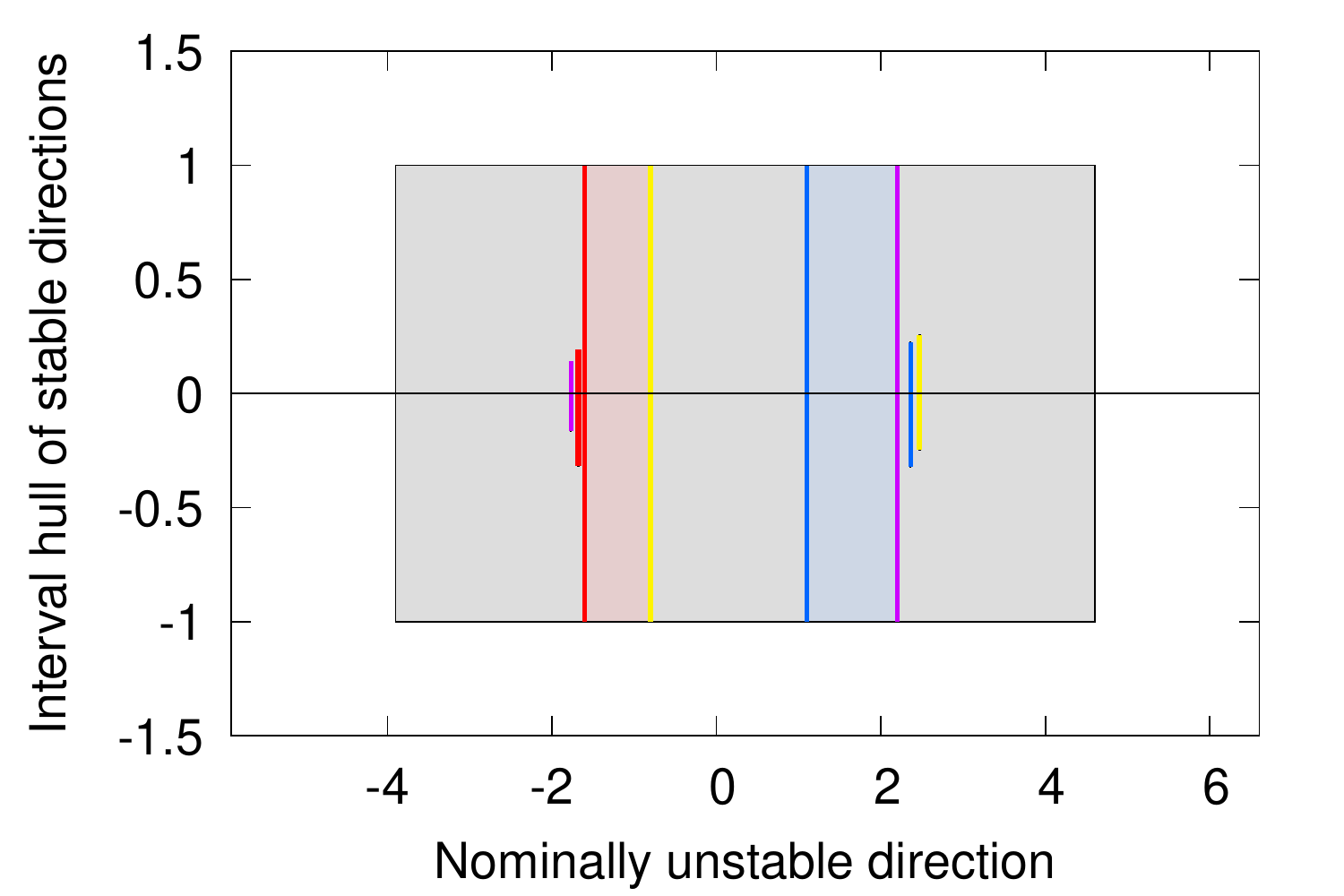}

\includegraphics[width=6.5cm]{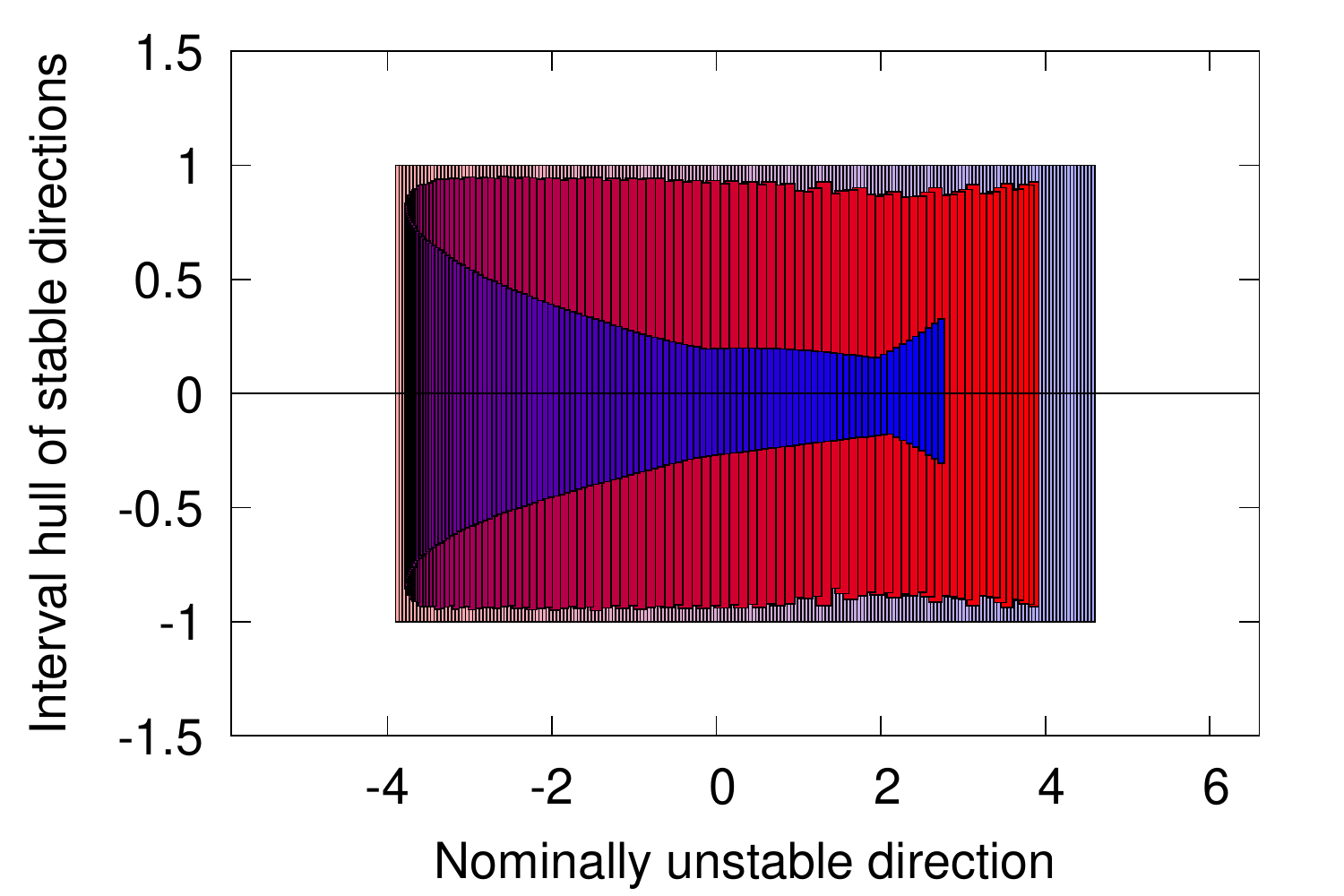}
\includegraphics[width=6.5cm]{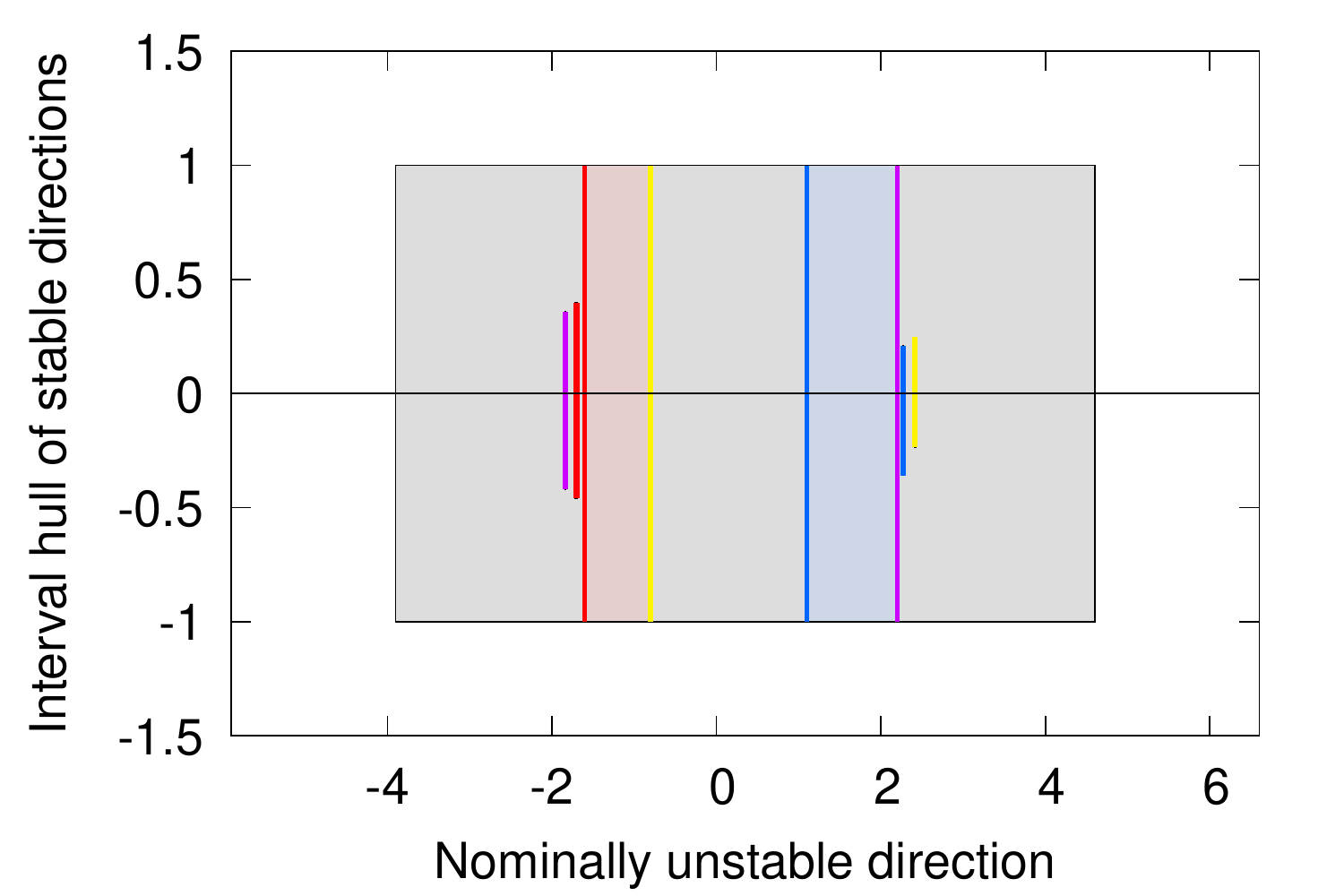}
}
\caption{\label{fig:proof-rossler-delayed}
The rigorous estimates obtained in the computer assisted
part of the proof of Theorem~\ref{thm:rossler-delayed}. The cases
(a)-(c) are presented top to bottom, respectively.
The left picture shows the representation of the computer assisted proof
of the trapping region $X_A$. The set is divided into 200 pieces $X_{A,i}$
along the $W_u$ direction, each piece is coloured according to
ascending number. Then for each piece $X_{A,i}$ the image $P(X_{A,i})$ is
computed and drawn in the same colour (but with increased intensity).
The dimension of the boxes in the $y$ coordinate represents the hull
of the nominally stable part of the set $P(X_{A,i})$, i.e.
the interval $I_i = [y_{lo}, y_{up}]$ such that all $\pi_{A_j} P_c(X) \subset I_i$,
for $j \in \{3, \ldots, M\}$ and $\pi_{\Xi_j} (P_c(X)) \subset I$ for
$j \in \{ 1, \ldots, p \cdot d \}$. Obviously, each
$I_i \subset B_1\left(0, \max (P(r_A), P(r_\Xi)) \right)$. A clear
evidence of the Smale horseshoe-like dynamics can be seen in
the picture, as the box is folding on itself under the map $P$.
On the right picture one there are represented the sets $X_1$
(light red, with red and yellow borders)
and $X_2$ (light blue, with blue and purple borders). The images
of the borders under the map $P^2$ are presented as
lines (in fact thin boxes) in the grey area outside
$X_1 \cup X_2$. It is evident that $P(W_{1,l})$ (red) and $P(W_{2,r})$ (purple)
are both mapped to the left of both sets and
$P(W_{1,r})$ (yellow) and $P(W_{2,l})$ (blue) are mapped
to the right. Therefore condition (CC2A) is satisfied between
the sets $X_1$ and any of $X_i$'s, and condition (CC2B) between $X_2$ and
any $X_i$, $i \in \{ 1, 2 \}$. Please consult online version of 
the plots for better quality.
}
\end{figure}

Fig.~\ref{fig:attractor-rossler-delayed} shows the numerical
representations of the apparent strange attractor in the respective systems,
while Fig.~\ref{fig:proof-rossler-delayed} depicts the
computed estimates of the proof in a human-friendly manner.
The total running time of the proof in (a) is around 16 minutes, and 
the cases (b) and (c) of around 23 minutes. Computations were done 
on a laptop with Intel\textregistered\ Core${}^{\mathtt{TM}}$ i7-10750H 2.60GHz CPU.
The majority of the computations is done in the proof of trapping region
$X_A$, which must be divided into 200 pieces along the vector $W_u$.
Those computations are easily parallelized (each piece
computed in a separate thread).
The data and programs used in the proofs are described in more details
in Appendix~\ref{app:data}, together with the links to source codes. 

\subsection{Unstable periodic orbits in Mackey-Glass equation}

\begin{figure}
\center{
\includegraphics[height=4.5cm]{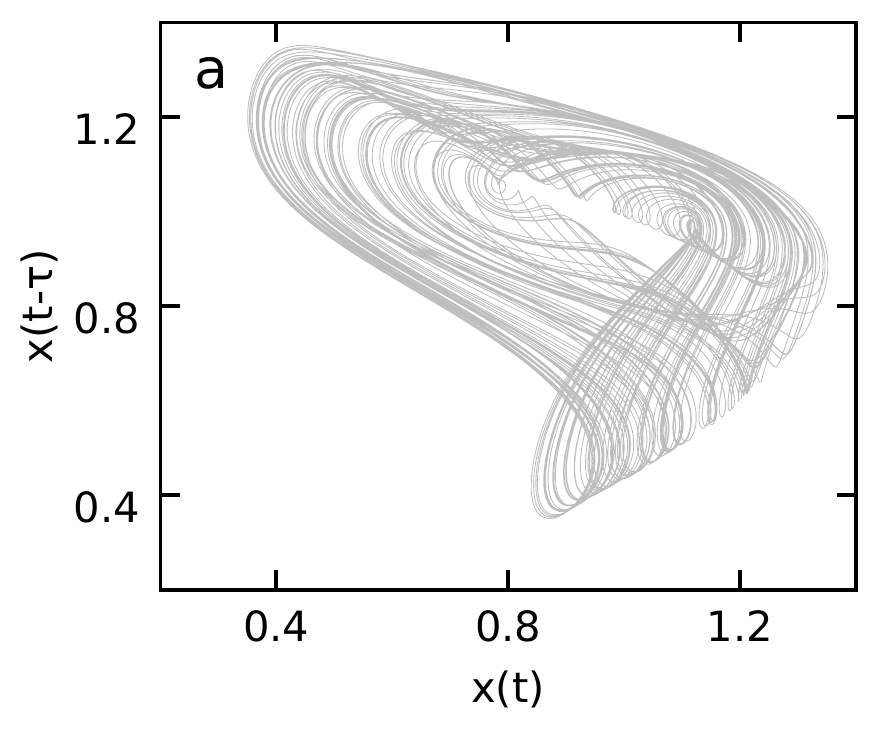}
\quad 
\includegraphics[height=4.5cm]{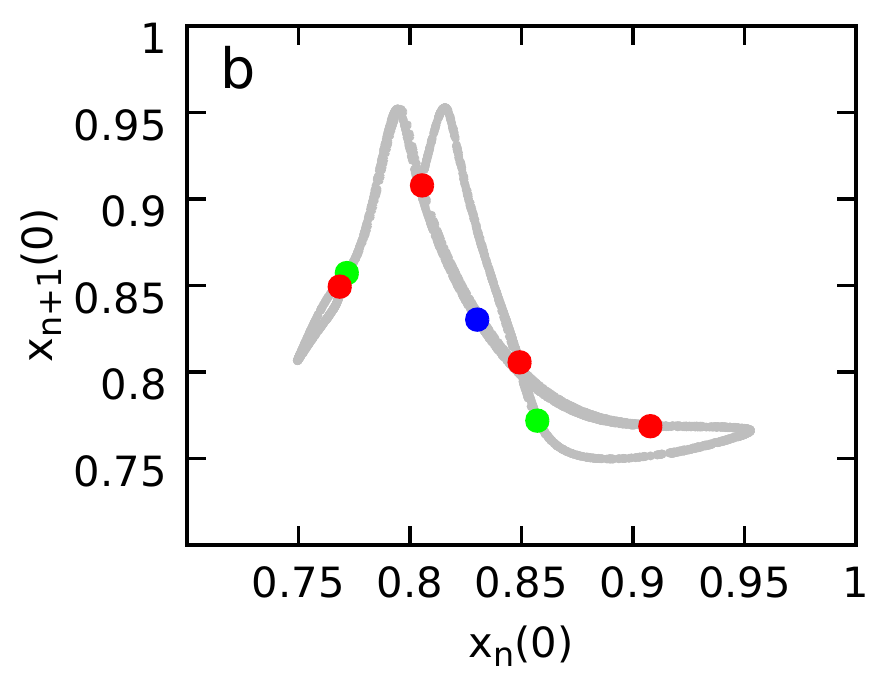} 

\includegraphics[height=4.5cm]{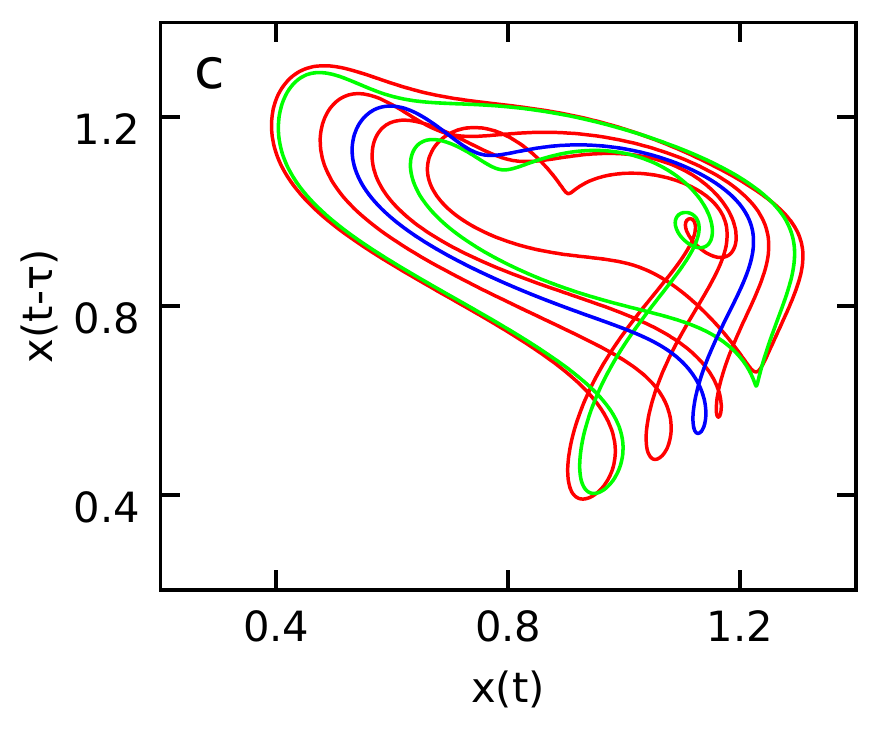}
\includegraphics[height=4.5cm]{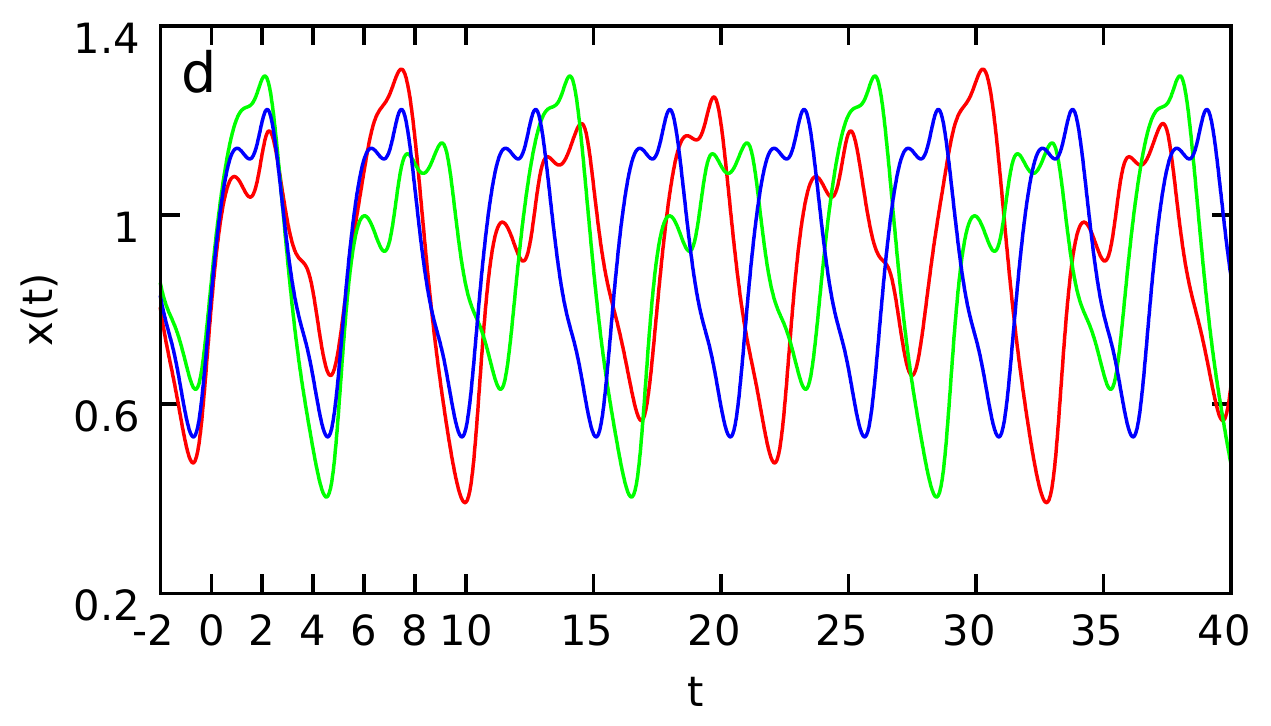}
}
\caption{\label{fig:mg-periodic} 
(a) The apparently chaotic attractor of the Mackey-Glass equation \eqref{eq:mg} 
for the classical parameter values $\tau = 2$, $n = 9.65$, $\beta = 2$, 
$\gamma = 1$ \cite{mackey-glass}. The attractor is drawn for a single very long 
solution, whose time-delay embedding coordinates $(x(t), x(t-\tau))$ are shown 
in the picture. (b) The representation of the attractor
drawn in the coordinates $\left(x_n(0), x_{n+1}(0)\right)$, where 
$x_{n+1} = P(x_n)$, $x_n, x_{n+1} \in C([-\tau, 0], \R)$. 
The map $P$ is constructed on the section $S = \{ x : x(t) = x(t-\tau), x(t) < 0.96 \}$,
see Figure.~13 in \cite{mackey-glass-scholarpedia}. The periodic points
$T^1, T^2, T^4$ of respective periods 1, 2 and 4 for map $P$ are drawn in colors
blue, green, red. (c) The same solutions are drawn in the time-delay embedding
of the attractor and (d) as the solutions over time long enough to contain
 basic periods of all presented solutions.
}
\end{figure}

In this application we study the following scalar equation:
\begin{equation}
\label{eq:mg}
x'(t) = -\gamma \cdot x(t) + \beta \cdot \frac{x(t-\tau)}{1 + (x(t-\tau))^n}.
\end{equation}
In the paper \cite{mackey-glass}, the authors shown
numerical evidence of chaotic attractor in that system, see
Fig.~\ref{fig:mg-periodic}a. In their work, Mackey and Glass
used the following values of parameters:
$\tau = 2$, $n = 9.65$, $\beta = 2$, $\gamma = 1$. In our
previous work \cite{nasza-praca-focm} we have shown existence of several (apparently)
stable periodic orbit for $n \le 8$. In this work we show that the new
algorithm, together with the fixed point index, can be used
to prove more diverse spectrum of results. We prove existence
of several (apparently) unstable periodic orbits for the classical
values of parameters, for which the chaotic attractor is observed, $\tau = 2$, 
$n = 9.65$, $\beta = 2$, $\gamma = 1$.

\begin{rem}
\label{rem:unit-delay}
In what follows we get rid of the variable delay $\tau$ and we 
rescale the system to have unit delay by the change of variables:
$y(t) = x(\tau \cdot t)$. It is easy to see, that the equation
\eqref{eq:mg} in the new variables becomes:
\begin{equation*}
\label{eq:mg-tau-1}
y'(t) = \tau \cdot f(y(t), y(t-1)),
\end{equation*} 
that is, we can remove parameter $\tau$ by rescaling 
$\beta$ and $\gamma$ to $\bar{\beta} = \tau \cdot \beta$
and $\bar{\gamma} = \tau \cdot \gamma$.
\end{rem}

We state the following:
\begin{theorem}
\label{thm:mg-periodic} Each of the three 
approximate solutions $\bar{T}^i$ shown in Figure~\ref{fig:mg-periodic}(c)-(d)
has a small, explicitly given vicinity $V_i \subset C^n_p$ with $n=4$ and $p=128$ 
of the initial segment $\bar{T}^i_0$ such that there exists a true periodic 
solution $T^i$ with the initial segment $T^i_0 \in V_i$ of the Mackey-Glass 
equation \eqref{eq:mg} for the classical parameter values $\tau = 2$, 
$n = 9.65$, $\beta = 2$, $\gamma = 1$ \cite{mackey-glass}.
\end{theorem}
\textbf{Proof of Theorem~\ref{thm:mg-periodic}}: we use the parameters $\beta = 4$ and $\gamma = 2$, $n = 9.65$ and $\tau = 1$ 
in \eqref{eq:mg} and we use Remark~\ref{rem:unit-delay}. The proof is similar 
to that of Theorem~\ref{thm:rossler-delayed} and boils down to checking 
appropriate covering relations. The initial segments $\bar{T}^i$ 
lie on the section $S = \left\{ x \in C([0,1],\R): x(0) = x(-1), x(t) < 0.96 \right\}$.
The index $i$ corresponds to the basic period of the solution $T^i$
as a periodic point of a map $P : S \to S$. In the coverings we
use map $P^2$ to guarantee that the return time $t_P$ to the section 
is \emph{long enough}. 

Each of the $V_i = X(N_i, \Xi_i)$ is given with $N_i = \bar{T}^i_0 + C_i \cdot r_i$
with $r_i = \{ 0 \} \times W_i^u \times [-1,1]^{M-2}$. Additionaly, in case of $T^4$ we 
have another set $V'_4 = X(N'_4, \Xi'_4)$ with $N'_4$ of the similar
form: $N'_4 = P^2(\bar{T}^4) + C'_4 \cdot r'_4$. In other words, the
origin point of the set $N'_4$ is the second iteration of the Poincar\'e map
$P^2$ of the initial segment of $T^4$.
The sets are obtained as described in Appendix~\ref{app:data}. Each of
these sets define a section $S_i = \{ x \in C^n_p : c_i \cdot (a(x) - \bar{T}^i_0) = 0 \}$ (different from $S$), where $c_i = (C_i)_{\cdot,1}$ - the first column of the matrix $C_i$.
The reason for that is described in the Appendix B, and boils down to assure that $\diam(t_P(X_i))$ is as small as possible.

We will show that:
\begin{equation}
\label{eq:mg-coverings}
V_1 \cover{P_{S_1 \to S_1}} V_1,
\quad\quad\quad V_2 \cover{P_{S_2 \to S_2}} V_2,
\quad\quad\quad V_4 \cover{P_{S_4 \to S'_4}} V'_4 \cover{P_{S'_4 \to S_4}} V_4,
\end{equation}
where the Poincar\'e maps $P_{S_i \to S_j}$ are derived from
the flow of Eq.~\eqref{eq:mg} and maps indicated
sections: $P_{S_i \to S_j} : S_i \to S_j$, with additional assumption that
the return time $t_P$ is long enough. We will drop
the subscripts if they are easily known from the context.

For $T^1$ we have:
\begin{itemize}
\item for all $i > 2$, $\left|\pi_i P_c(X(N_{1}, \Xi_{i}))\right| = 0.614451801967851 < 1$
\item for all $i$, $\left|\pi_{\Xi_i} P_c(X(N_{1}, \Xi_{i}))\right| = 0.999998174289212 < 1$
\item $\pi_2 P_c(X(N_{1}^{r}, \Xi)) = -{}_{3.845940820239275}^{4.514877050431105} < -1 = \pi_2 N_{c}^{l}$
\item $\pi_2 P_c(X(N_{1}^{r}, \Xi)) = {}_{3.827847664967472}^{4.496773405715568} > 1 = \pi_2 N_{c}^{r}$
\end{itemize}
For $T^2$ we have:
\begin{itemize}
\item for all $i > 2$, $\left|\pi_i P_c(X(N_{2}, \Xi_{2}))\right| \le 0.731193331043839 < 1$
\item for all $i$, $\left|\pi_{\Xi_i} P_c(X(N_{2}, \Xi_{2}))\right| \le 0.999996951451891 < 1$
\item $\pi_2 P_c(X(N_{2}^{r}, \Xi_2)) = {}_{3.995778903452447}^{5.033339010859840} > 1 = \pi_2 N_{2,c}^{r}$
\item $\pi_2 P_c(X(N_{2}^{r}, \Xi_2)) = -{}_{3.978765934264806}^{5.016322912452834} < -1 = \pi_2 N_{2,c}^{l}$
\end{itemize}
For $T_4$ we have:
\begin{itemize}
\item for all $i > 2$, $\left|\pi_i P_c(X(N_{4}, \Xi_{4}))\right| \le 0.999948121260377 < 1$
\item for all $i$, $\left|\pi_{\Xi_i} P_c(X(N_{4}, \Xi_{4}))\right| \le 0.956276660970399 < 1$
\item $\pi_2 P_c(X(N_{4}^{l}, \Xi_4)) = -{}_{1.1221122505976317}^{3.2368193002087367} < -1 = {N'_{4}}^{l}_{c}$
\item $\pi_2 P_c(X(N_{4}^{r}, \Xi_4)) = {}_{1.1239689716822263}^{3.2385726261859316} > 1 = {N'_{4}}^{r}_{c}$
\end{itemize}
and
\begin{itemize}
\item for all $i > 2$, $\left|\pi_i P_c(X(N'_{4}, \Xi'_{4}))\right| \le 0.898580326387734 < 1$
\item for all $i$, $\left|\pi_{\Xi_i} P_c(X(N'_{4}, \Xi'_{4}))\right| \le 0.952378028038733 < 1$
\item $\pi_2 P_c(X({N'_{4}}^{l}, \Xi'_4)) = {}_{1.6331410212899785}^{3.0485204456349866} > 1 = {N_{4}}^{r}_{c}$
\item $\pi_2 P_c(X({N'_{4}}^{r}, \Xi'_4)) = -{}_{1.6341550945779965}^{3.0495636957165507} < -1 = {N_{4}}^{l}_{c}$
\end{itemize}
All those inequalities satisfy appropriate assumptions of 
Lemma~\ref{lem:one-unstable}. Therefore all the coverings from \eqref{eq:mg-coverings}
exist and, from Theorem~\ref{thm:symbolic-dynamics}, we infer existence
of appropriate periodic points $T^i_0 \in V_i$. 
\qed

\input tables/table2-content.tex

The diameters of the sets expressed in commonly used functional norms 
are presented in Table~\ref{tab:diameters}.
The data and programs used in the proofs are described in more details
in Appendix~\ref{app:data}, together with the links to source codes. 

\subsection{A comment about the exemplary systems}

Both R\"ossler and Mackey-Glass systems studied as an exemplary application
in this work are chaotic for the parameters used. However, Mackey-Glass 
system is a scalar equation, so the chaos present in the system must be
a result of the infinite nature of the phase space and the delay
plays a crucial role here. It is not clear if the dynamics can be approximated
with a finite number of modes, and how to choose good coordinate
frame to embed the attractor. The 
R\"ssler system on the other hand is a 3D
chaotic ODE (for parameters specified), and the chaotic behaviour
is the result of the dynamic in this explicitly 
finite dimension space. The systems of the form 
\eqref{eq:rossler-perurbed} are small perturbations of the ODE
and thus one can expect the dynamics of the ODE persist in some sense,
at least for $\epsilon$ small enough \cite{nasza-epsilon-perturb-dde}. 
It is much easier to propose sets for the covering relations inherited directly
from the coverings in finite dimension for unperturbed system, 
see Appendix~\ref{app:data}, where we use the flow of unperturbed ODE 
to generate the apparently unstable direction for the trapping region
containing the attractor.

%% file: tables/table2-content.tex
\begin{table}
\begin{center}
\begin{tabular}{|c|c|c|c|c|}
\hline 
			& $L_\infty$ 	&  $L_2$ 	&  $H_4$ 	&  T (expressed in $\tau$)	 	\\
\hline 
$T^1_0$ & $2.74231097479455 \cdot 10^{-7}$ & $3.10483831050838 \cdot 10^{-6}$ & $21.9495663834241$ & $2.632897_{884924421}^{901501874}$ \\
\hline 
$T^2_0$ & $1.34240247683063 \cdot 10^{-7}$ & $1.52442781918968 \cdot 10^{-6}$ & $23.4410359636472$ & $5.982965_{269668710}^{324098800}$ \\
\hline 
$T^4_0$ & $2.09990758524436 \cdot 10^{-8}$ & $1.91536904191854 \cdot 10^{-6}$ & $26.9986770914825$ & $11.40640_{3860772954}^{4205303446}$  \\
\hline 
\end{tabular}
\end{center}
\caption{\label{tab:diameters}The basic period $T$ of each 
solution and the diameters of the sets $V_i$ estimated (upper bounds) 
in various functional norms: 
$\|x\|_{L_\infty} = \sup_{[-\tau, 0]}|x(t)|$, 
$\|x\|_{L_2} = \left(\int_{-\tau}^{0} (x(t))^2 dt\right)^{\frac{1}{2}}$, 
$\|x\|_{H_4} = \sum_{i=0}^{4} \|x^{(i)}\|_{L_2}$. 
Note that the period $T$ is expressed as the number of full
delays, and will be doubled for the original system
with $\tau = 2$, $\beta = 2$, $\gamma = 1$ and $n = 9.65$.
}
\end{table}

%% file: sections/6_appendix_A_lohner.tex

\section{\label{app:lohner}Lohner-type algorithm for control of the wrapping effect}

In this Appendix we present technical details of the implementation of an
efficient Lohner-type control of the wrapping effect.

\subsection{Lohner's algorithms and Lohner's sets - basic idea}

\begin{figure}
    \centering{
        \includegraphics[width=60mm]{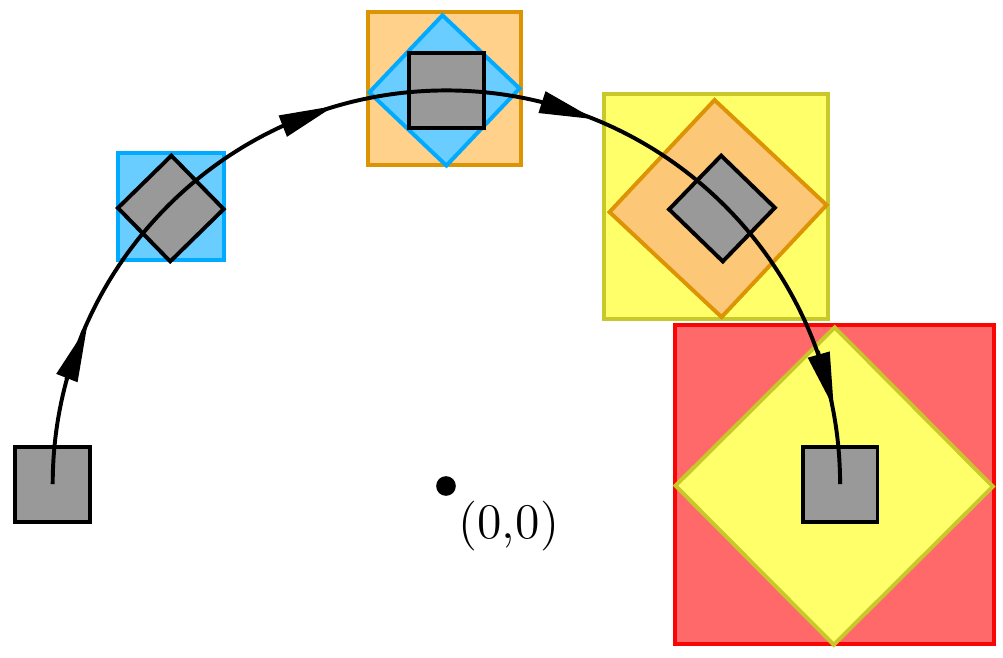}
    }

    \caption{\label{fig:wrapping-effect}
	An illustration taken from \cite{nasza-praca-focm}
	of the wrapping effect problem for a classical,
	idealized mathematical pendulum ODE $\ddot{x} = -x$. The picture
	shows a set of solutions in the phase space $(x, x')$. The grey
     boxes shows the set of initial conditions (a box)
	moved by the flow. The coloured boxes present the wrapping effect
	occurring at each step when we want to enclose the moving points in a
	product of intervals in the basic coordinate system. For example, the
     blue square on the left encloses the image of the first iteration. Its image is
     then presented with blue rhombus which is enclosed again by an orange
     square. Then the process goes on. We see that the impact of the wrapping
     effect rapidly becomes overwhelming. 
    }
\end{figure}

Lohner \cite{capd-lohner} proposed, in the case of finite dimensional
maps $G : \R^M \to \R^M$, to use a decomposition of the rigorous
method for $G$ into the numerical (approximate) part $\Phi : \R^M \to \R^M$, 
that can be explicitly differentiated w.r.t. initial value $x$, and 
the remainder part of all the errors $\Rem$, such that
$G(x) \in \Phi(x) + \Rem(X)$ for all $x \in X$. The Lohner's original
idea was to use Mean Value Form of the $\Phi$ part to ,,rotate'' the
coordinate frame to reduce the impact of the so called the \emph{wrapping
effect} encountered when using interval arithmetic. Without the change 
of local coordinate frame for the set, the result of computations 
would be represented as an Interval box in $\R^M$ and big over-estimates 
would occur, see for example Fig.~\ref{fig:wrapping-effect}. The Lohner's idea
 allows to reduce this problem significantly.

In a case of a general map $G$ one can 
use the mean value form  for $\Phi$ to get the following:
\begin{align}
\label{eq:mean-value-form}
\Phi(z) \in \Phi(x) + \left[D \Phi([X])\right] \cdot (X - x)
\end{align}
for all $z \in X \subset \R^M$, and the point $x \in X$ is just any point, 
but usually chosen to be the centre of the set $X$. 
Here $[X] \in \I^M$ is an interval hull of $X$ and $\left[D\Phi([X])\right]$ is
an interval matrix that contains the true Jacobians $D\Phi(z)$
at all $z \in [X]$. Then, the strategy to reorganize operations
depends on the shape of the set. In the simplest case let assume
\begin{equation}
\label{eq:basic-lohner-set-structure}
X = x + C \cdot r_0 + r 
\end{equation}
where $C$ is a linear transform $\R^M \to \R^M$, 
$x \in \R^M$, and with interval vectors $r_0, r \in \I^M$ centred at $0$. 
Using \eqref{eq:mean-value-form} we have:
\begin{align}
\notag
\Phi(z) & \in \Phi(x) + \left[D \Phi([X])\right] \cdot \left(C \cdot r_0 + r\right), \\
\label{eq:lohner-before-reorganisation}
& = \Phi(x) + \left(\left[D \Phi([X])\right] \cdot C\right) \cdot r_0 + \left[D \Phi([X])\right] \cdot r
\end{align}
It is now evident, that the result set has structure similar to 
\eqref{eq:basic-lohner-set-structure}:
\begin{equation}
\label{eq:lohner-one-step}
G(z) \in Y := \bar{x} + \bar{C} \cdot r_0 + \bar{r}.
\end{equation}
With some additional reorganizations to keep $x$ and $C$ 
as thin as possible (e.g. point vector and matrix) and 
vectors $r$ and $r_0$ centred at $0$, we arrive at
the following Lohner-type algorithm:
\begin{align}
\label{eq:Phi-x}			\bar{x}	& := \midpoint(\Phi(x) + \Rem(X)) \\ 
\label{eq:dominant-mul} 	S 		& := \left[D \Phi([X])\right] \cdot C \\
\label{eq:Phi-C}			\bar{C} 	& := \midpoint(S) \\
\label{eq:Phi-r}			\bar{r} 	& := (S - \midpoint(S)) \cdot r_0 + \left[D \Phi([X])\right] \cdot r + \left(\Phi(x) + \Rem(X) - \midpoint\left(\Phi(x) + \Rem(X)\right)\right),
\end{align}
where $\midpoint(\cdot)$ is an operation that returns the middle 
point of the interval vector or matrix. The terms 
in \eqref{eq:Phi-r} might require some comments:
the first term is the error left from the part $S \cdot r_0$ 
introduced by taking midpoint of matrix $S$ as $\bar{C}$
in \eqref{eq:Phi-C}. Second term is just applying mean value
form on the $r$ part. Third term is the error introduced
after taking midpoint of the sum in \eqref{eq:Phi-x} as the
new reference $\bar{x}$. If the matrix $\left[D \Phi([X])\right]$
and the term $\Rem(X)$ are ,,thin'' (i.e. their entries as 
intervals have small diameter) then we hope the newly
introduced errors should be small in comparison
to the term $\bar{C} \cdot r_0$. 

This is just one of the proposed shapes of the
set in Lohner's algorithm, the so called ,,parallelepiped 
($C \cdot r_0$ part) with interval form of the remainder (the $r$ part is an interval box in $\I^M$)''. 
A more general approach is the ,,doubleton set'':
\begin{equation}
\label{eq:rect2-lohner-set-structure}
X = x + C \cdot r_0 + B \cdot r 
\end{equation}
where matrix $B$ is chosen in some way (to be described later).
The Lohner algorithm is more involved in this case: 
\begin{eqnarray}
\notag						\bar{x}	& := & \midpoint(\Phi(x) + \Rem(X)) \\ 
\notag				 		S 		& := & \left[D \Phi([X])\right] \cdot C \\
\notag						\bar{C}	& := & \midpoint(S) \\
\label{eq:rect2-QR}		Q \cdot R & := & \midpoint\left(\left[D \Phi([X])\right] \cdot B\right) \\
\label{eq:rect2-B}		\bar{B} & := & Q \\
\label{eq:rect2-Phi-r}	\bar{r}	& := & \left( Q^{-1} \cdot \left[D \Phi([X])\right] \cdot B\right) \cdot r \ + \\
\notag & + & \left(Q^{-1} \cdot (S - \midpoint(S))\right) \cdot r_0 \ + \\
\notag & + & \left(Q^{-1} \cdot \left(\Phi(x) + \Rem(X) - \midpoint\left(\Phi(x) + \Rem(X)\right)\right)\right).
\end{eqnarray}
The difference from the previous algorithm in \eqref{eq:Phi-r} is 
in Eqs.~\eqref{eq:rect2-QR}-\eqref{eq:rect2-Phi-r}. The idea of
the improvement over the previous version is that one hope
the first term in \eqref{eq:rect2-Phi-r} to have some wrapping 
effect controlled by the matrix $Q^{-1}$, when doing interval enclosure. 
The choice about $Q$ and $Q^{-1}$ is done in Eq.~\eqref{eq:rect2-QR} and 
depends on the algorithm implementation.
Ideally, we should set $R = Id$, so that
\begin{equation}
Q = \midpoint\left(\left[D \Phi([X])\right] \cdot B\right),
\end{equation}
just as in case of \eqref{eq:Phi-C}. However, we need to compute 
rigorous inverse of this matrix, which might be
either computationally expensive, very difficult or even impossible.
On the other hand, we can choose $Q = Id$, which transforms
the algorithm into the previous one (for sets with the interval 
form of the remainder, i.e. defined as \eqref{eq:basic-lohner-set-structure}). 
Finally, the most commonly used method
is to apply (rigorously) any QR decomposition in \eqref{eq:rect2-QR}
so that the matrix $Q^{-1} = Q^{T}$ is easily obtainable. 
This strategy will be crucial later to get better results for
DDEs in the case of $d > 1$ (systems of equations).  

One last remark, before we move on to the application
of the Lohner's idea in the context of DDEs, is that the method
can be applied also to functions $G : \R^{M_1} \to \R^{M_2}$ where
the dimensions of the domain and the image is different: 
$M_1 \ne M_2$. Formulas~\eqref{eq:Phi-x}-\eqref{eq:Phi-r} are all valid, 
but one must be very careful about dimensions of all vectors and matrices 
involved in the computations. 

\subsection{\label{sec:optimization-idea}Lohner's algorithm - complexity and optimization idea}

Lohner's algorithm complexity is dominated by the
two main factors: computation of $\left[D \Phi([X])\right]$
used in \eqref{eq:dominant-mul} and multiplication of matrices. 
Additionally, there might be some set-structure dependent complexity,
such as the need to compute the QR decomposition and the inverse 
of the matrix $Q$ in \eqref{eq:rect2-QR}.
All other operations such as matrix-vector multiplication and 
matrix-matrix and vector-vector additions have lower computational 
complexity. Computation of $\left[D \Phi([X])\right]$ cannot be
avoided and has complexity depending on the complexity of
the formula for $\Phi$. The complexity of matrix-matrix 
multiplication is $O(M^3)$, not taking into account 
the possible faster (and more complicated) multiplication 
algorithms (e.g. Strassen's algorithm and similar). In the rest 
of the appendix we will discuss the possible simple and 
effective optimization of those dominant operations based on 
the sparse structure of $\left[D \Phi([X])\right]$.
We will recall from \cite{nasza-praca-focm}
that the $\left[D \Phi([X])\right]$ is very sparse in the case of the 
integration algorithm $\mathcal{I}$ for DDEs. We will extend and 
provide nicer description for the ,,fast matrix multiplication'' method 
presented in \cite{nasza-praca-focm} that is easily generalized
for any used variables $u$ in the case of multiple delays. Moreover, 
later on, we will discuss possible shape of the matrix $B$ in 
\eqref{eq:rect2-lohner-set-structure} which will provide better results 
but without significant cost in the computational complexity.

The matrix multiplication optimization idea was first proposed 
in \cite{nasza-praca-focm} for a specific case of DDEs with one delay
i.e. of the form \eqref{eq:dde}. Now we propose a more elegant and more
general implementation, that will be suitable for implementation
of the problem \eqref{eq:dde-many-delays} and with $d > 1$ (systems
of equations).  
The idea is based on the decomposition of the computation of 
$A \cdot B$, $A \in \matrices(M, M)$, 
$B \in \matrices(M, N)$ into consecutive computation of 
products of $A_{i,\cdot}$ - $i$-th row of $A$ and $B$. 
In our case we think of $A$ as $A = \left[D \Phi([X])\right]$. 
Let assume that $A_{i,\cdot}$
has a lot of zeros (it is sparse). Let denote by $u(\cdot)$ 
(name conflict intentional) the following projection 
\begin{equation*}
u = (\pi_{l_1}, \pi_{l_2}, \ldots, \pi_{l_k}), \quad \forall l_k : A_{i,l_k} \ne 0
\end{equation*} 
The function $u : \R^M \to \R^k$ and reduces the dimension of the
vectors from $M$ to $k$, so we will call it a \emph{reduction}.
For matrix $B$ we define:
\begin{equation}
\label{eq:matrix-projection}
u(B) = \left(\begin{array}{c}
B_{l_1,\cdot} \\
B_{l_2,\cdot} \\
\vdots \\
B_{l_k,\cdot}
\end{array}\right),
\end{equation}
that is, $u(B) \in \matrices(k, N)$ contains 
all rows corresponding to the variables used in the reduction $u$.
It is now easy to see that
\begin{equation}
\label{eq:reduced-multiplication}
A_{i, \cdot} \cdot B = u(A_{i, \cdot}) \cdot u(B),
\end{equation}
and the complexity of the operation is reduced from $O(M \cdot N)$ to 
$O(k \cdot M)$. We can now apply the multiplication in a loop
for all $i$ separately, changing the $u$ accordingly (or
using the same $u$ for some coordinates and do
multiple rows of $A$ at the same time). 

We note that, in the simplest case, when $u = (\pi_l)$ 
(only one non-zero element in the $i$-th row of $A$),
and $A_{i,l} = 1$ we have:
\begin{equation}
\label{eq:fast-mul-Id}
A_{i,\cdot} \cdot B = 1 \cdot B_{l, \cdot}
\end{equation}
and we can skip multiplication completely, changing
it to a shift (selection of a given row). This will  
be used when $A$ has a large $Id$ block in its structure. 

\subsection{Lohner-type algorithm for DDEs integrator - preparation}

Now we apply the Lohner strategy to our rigorous DDE integrator $\mathcal{I}$. 
We decompose the general method for many delays from Section~\ref{sec:integrator-many-delays}  into 
the numerical procedure $\Phi : \R^n_{p,q} \to \R^n_{p,q+1}$ 
and the remainder $\Rem : \R^n_{p,q} \times \I^{d \cdot p} \to \I^n_{p,q+1} \times \I^{d \cdot p}$
in the following way:
\begin{align}
\label{eq:Phi-neta}    n     		& := n(\eta, f) \\
\label{eq:Phi-y} y(u(x)) 			& := \left(z(x), w_{n+1} * F^{[n]}\left(u(x)\right) \right) \\
\label{eq:Phi-all} \Phi(a(x)) 		& := \left(\taylor(y(u(x)); h), y(u(x)), j_2(x), \dots, j_{p-1}(x)\right) \\
\label{eq:Rem-A} 	\Rem_A(x, \ixi) 	& := \left([F]_{[n+2]} \cdot [0,h] \cdot h^{n+1}, 0, \dots, 0\right) \in \I^n_{p,q+1} \\
\label{eq:Rem-R} 	\Rem_R(x, \ixi) 	& := \left(\frac{[0,h]}{n+2} \cdot [F]_{[n+2]}, \ixi_2(x), \dots, \ixi_{p-1}(x)\right) \in \I^{d \cdot p}
\end{align}
where $F^{[n]}$ as in \eqref{eq:jet-reccurence-general}, 
$[F]$ as in \eqref{eq:FU}, and $a(x) = (z(x), j(x))$ is the
finite dimensional part of the description $(z, j, \xi)$ of $x$. 
The order $n$ of the new jet \eqref{eq:Phi-neta} comes 
from \eqref{eq:neta} in the algorithm, see details there.
The intermediate variable $y$
is defined in \eqref{eq:Phi-y} to shorten \eqref{eq:Phi-all} 
and underline the dependence on the ,,used variables'' $u(x)$. 
We remind that the ,,used variables'' vector $u(x)$ is defined
for DDE~\eqref{eq:dde-many-delays} with $m$ delays 
$\tau_1 = p_1 \cdot h = \tau$ (i.e. $p_1 = p$), 
$\tau > \tau_i = p_i \cdot h > \tau_j = p_j \cdot h$
for $i, j \in \{2,.., m\}$, $i < j$, $p_i, p_j \in \{1,..,p-1\}$ as:
\begin{eqnarray*}
u(x) = (z(x), j_{p_1}(x), j_{p_2}(x), \ldots, j_{p_m}(x)).
\end{eqnarray*}
Please note that, with some abuse of notation, we can think
of $u$ as a vector in $\R^{dim(u)}$. If $x \in C^n_p$
then $dim(u) = d(1 + m \cdot (n+1))$.

First we observe that the map $\Phi$ is well defined
map from $\R^M \to \R^{M + d}$ with $M =M(d, p, \eta)$ and 
it can be differentiated w.r.t. $a$ if $f$ is smooth enough, 
for example as in our simplifying assumption $f \in C^\infty$. 
Therefore, the Lohner algorithm might be applied ,,as it is'' to
the algorithm in the pair of Eqs~\eqref{eq:Phi-all}-\eqref{eq:Rem-A} 
(the $A$-part of the set). However, this approach would be highly
ineffective in applications, we will demonstrate now why.

\subsection{Naive, straightforward implementation and the structure of $D\Phi$}
For simplicity, let assume we deal with the interval representation
of the error term $B \cdot r = Id \cdot r$ in the Lohner set 
\eqref{eq:rect2-lohner-set-structure} for $X = X(A,R) \subset C^n_p$.
In that case, it is easy to observe that the dominant operation in the Lohner's 
algorithm (in terms of computational complexity) is the matrix-matrix multiplication 
in Eq.~\eqref{eq:dominant-mul}.
The application of the standard naive matrix multiplication leads to the 
computational complexity of
$O(M^3) = O( (d \cdot n \cdot p)^3) $ since the matrix dimensions of 
both $D \Phi(x)$ and $C$ dimensions are of the order of $O(d \cdot n \cdot p)$. 
This is also true (under some assumption) if the size of the representation 
$M$ grows as the algorithm is iterated.  
Indeed, let consider $X_0 = X(x + C \cdot r_0 + r, R)$ with $C \in \matrices(M, N)$
$r_0 \in \I^N$, $r \in \I^M = \I^n_{p,0}$. Usually $N = M$, but set-up with 
$N \le M$ might be beneficial in some applications. 
Let now consider the chain of sets $X_i = \Phi(X_{i-1})$ represented as
Lohner's sets \eqref{eq:rect2-lohner-set-structure}. We have, that in the $i$-th step 
($i \ge 1$) the sizes of the matrices involved in Eq.~\eqref{eq:dominant-mul} are 
$D\Phi(X_i) \in \matrices\left(M + d \cdot i, M+ d \cdot (i-1)\right)$,
 $C \in \matrices\left(M+d \cdot (i-1), N\right)$ and the result 
 matrix $S \in \matrices(M+d \cdot i, N)$. 
So the naive multiplication complexity is proportional to
\begin{equation*}
(M+d \cdot i) \cdot (M+d\cdot(i-1)) \cdot N \in O(M^3),
\end{equation*} 
provided that both $N, i \in O(M)$ - this is usually the case, as
$N > M$ does not make sense and $i >> M$ is not feasible computationally.

Please note that, for $M$ used in applications, we usually have $M \approx 1000$.
Therefore the matrix-matrix multiplication in the naive implementation 
of Lohner's algorithm does enormous $O(10^9)$ operations per integration
step. On the other hand, investigating Eqs.~\eqref{eq:Phi-all}-\eqref{eq:Rem-R} 
reveals that the dynamics on a lot of coefficients is simply 
a \emph{shift to the past}. Therefore, $\left[D \Phi([X])\right]$
has a following nice block structure:
\begin{equation}
\label{eq:DPhi-structure}
D\Phi(v) = \left(\begin{array}{cccc}
J_{11}(v)	& J_{12}(v)	& J_{13}(v) \\
J_{21}(v)	& J_{22}(v)	& J_{23}(v) \\
0 			& Id			& 0
\end{array}\right).
\end{equation}
The matrix $J_{11}(v) \in \matrices(d, d)$ corresponds
to the derivative $D_z \Phi_z(v)$, i.e. the 
derivative of the $z$-th component (value of the solution $x$ at 
current time $t = h$) w.r.t. to the change in $z(x)$ - the value of $x$ 
in the previous step (at $t = 0$). Likewise, 
$J_{13}(v) \in \matrices(d, d \cdot (n+1))$ corresponds
to the $D_{j_p} \Phi_z(v)$, $J_{21}(v) = D_{z} \Phi_{j_1}(v) \in \matrices((n+2) \cdot p, d)$,
and so on. We will denote the matrix $(J_{11}, J_{12}, J_{13})$
as $D_u \Phi_z (v)$ and $(J_{21}, J_{22}, J_{23})$ as 
$D_u \Phi_{j_1} (v)$, respectively. Here, we use the convention
that subindex such as $j_i$, $z$, etc. denotes the corresponding 
set of variables from the description of the function $x = (z, j, \xi)$.

Investigating the matrices $J_{12}(v)$ and $J_{22}(v)$ we see they 
correspond to the derivatives of $\Phi$ w.r.t. values at all 
intermediate delays $\tau_{p_i}$, $i > 1$, so they 
might also contain a large number of zeros (if the equation
does not depend on a particular $\tau_i$). 
When we are dealing with only one delay ($m = 1$), 
then $J_{12}(v) = 0$ and $J_{22}(v) = 0$. 
In that case, we can apply idea proposed in the 
previous section~\ref{sec:optimization-idea} to
get enormous reduction in the computational 
complexity. We will additionally introduce
the structure to the matrix $B$ defined
in \eqref{eq:rect2-lohner-set-structure} to
help with wrapping effect in the error part $B \cdot r$.

\begin{rem}
All matrices $J_{ij}$ in the actual implementation
of the method are computed using Automatic 
Differentiation techniques. Those techniques can be readily 
applied to any equation of the form \eqref{eq:dde-many-delays}
as long as $f$ is a composition of simple (well known) functions 
like $\sin$, $\exp$, etc. and standard algebraic operations 
$\times$, $\div$, $+$, $-$. We do not discus details of this 
matter in the article. 
\end{rem} 

\subsection{Lohner algorithm using $D_u\Phi_{z}$ and $D_u\Phi_{j_1}$ directly}

Let $X(A,R) \subset C^n_{p,q}$ be an fset such that 
\begin{equation}
\label{eq:A-structure}
A = x + C \cdot r_0 + B \cdot r
\end{equation}
as  in the Lohner structure \eqref{eq:rect2-lohner-set-structure} where 
$C \in \matrices(M, N)$, $M = \dimension(\R^n_{p,q})$. The matrix $B$
will have a special block-diagonal:
\begin{equation}
\label{eq:B-structure}
B = \left( \begin{array}{ccccc}
B_{z}		& 0 				& \cdots & 0						\\
0 			& B_{j_{1,[0]}}	& 0		  & \ddots  				\\
\vdots 		& 0					& \ddots & \ddots					\\
0	 		& \vdots		   	& \ddots & B_{j_{p,[\eta_p]}}
\end{array} \right),
\end{equation} 
where each $B_{b,b} \in \matrices(d, d)$. 

Now, we can apply \eqref{eq:lohner-one-step} to the pair of methods $(\Phi, \Rem_A)$
in Eqs.~\eqref{eq:Phi-all}-\eqref{eq:Rem-A} to get a new fset of
the same structure $Y = X(\bar{x} + \bar{C} \cdot r_0 + \bar{r},\Rem_R(X)) \subset C^n_{p,q+1}$ 
so that for all $z \in X(A, R)$ we have $\varphi(h, z) \in Y$.
Please note, that $\bar{C} \in \matrices(M + d, N)$, $\bar{r} \in \R^n_{p,q + 1} = \R^{M + d}$
and $r_0 \in \R^N$ stays the same as in the original Lohner's algorithm (this is important).
The extra $d$ rows in matrix $\bar{C}$ are due to the extra Taylor coefficient
computed at $t=0$. In general, in $i$-th iteration of the algorithm
the matrix $C_i$ will be of the dimension $\matrices(M_0+d \cdot i, N)$ and
the error term $r$ will be of dimension $\R^{M_0+d \cdot i}$, $B \in \matrices(M_0+d \cdot i, M_0+d \cdot i)$, 
where $M_0 = M(d, p_0, \eta_0)$ is the dimensional of the initial set $X_0 \in C^{\eta_0}_{p_0}$ at the 
beginning of the integration process.
In applications, we usually set $N = M_0$. 

\begin{rem}
There is a slight abuse of notation here, as we are using $x$ to
denote the base point of the set $A$ and, at the same time, usually it denotes
the segment of the solution $x \in X$. 
However, the two are used in a different context, so it should not create confusion 
(one is the Lohner's set of the $A$ part in $X(A, R)$, second is as an element of $X(A, R)$). 
We will state explicitly if $x \in X$ otherwise $x$ always denotes the mid point 
of $A$. Please also note, that, by definition, if $x \in X$, then naturally 
$a(x) \in A = x + C \cdot r_0 + B \cdot r$.
\end{rem}

Now, the crucial part is to look at each $d$-dimensional variable
$z(X)$ and $j_{i,[k]}(X)$ as a \emph{separate Lohner's set with
its own structure} inherited from the full set $X = X(A, R)$ and
apply the Lohner's algorithm separately on each part, together with  
the optimization idea from \ref{sec:optimization-idea}.

\subsubsection{The convention}
As with the $u$ in \eqref{eq:matrix-projection},
for a matrix $C \in \matrices(M, N)$ we define 
$z(C)$ and $j_{i,[k]}(C)$
as the matrix containing all the appropriate 
rows from $C$. Each $z(C)$ and $j_{i,[k]}(C)$
is therefore a matrix in $\matrices(d, N)$.

It is easy to see that if the set $X = X(A, R)$, with
$A$ as in \eqref{eq:A-structure}, then 
\begin{equation*}
z(X) = z(x) + z(C) \cdot r_0 + B_z \cdot z(r), 
\end{equation*}
where $B_z \in \matrices(d, d)$ given as in \eqref{eq:B-structure} and 
$z(C) \in \matrices(d, M)$. Similarly $j_{i,[k]}(X) = j_{i,[k]}(x) + j_{i,[k]}(C) \cdot r_0 + B_{j_{i,[k]}} \cdot j_{i,[k]}(r)$.

\begin{rem}
The use of the abstract operations $z(\cdot)$, $j_{i,[k]}(\cdot)$, and $u$
allows for a more general implementation of the methods,
independent of the actual storage organization of the data
in computer programs. 
\end{rem}

\subsubsection{The shift part}
First consider easy case of computing $j_{i}(\bar{X})$
in $\bar{A} = \Phi(a(X))$ for $i > 1$. We observe
that 
\begin{equation*}
D_{j_l} \Phi_{j_i}(a(X)) = \begin{cases}
Id_{d \times d} & l = i - 1 \\
0_{d \times d} & otherwise
\end{cases}.
\end{equation*}
as this is the case of 
the shift to the past in Eq.~\eqref{eq:shift-taylor}.
The procedure is exact (i.e. $\Rem_A(X)_{j_i} = 0$, see Eq.~\eqref{eq:Rem-A}) 
and no extra errors are introduced. Therefore:
\begin{equation*}
j_i(\bar{X}) = j_{i-1}(X),
\end{equation*}
and using observation \eqref{eq:fast-mul-Id} we have
for all appropriate $k$:
\begin{eqnarray}
\notag j_{i,[k]}(\bar{C}) &=& j_{i-1,[k]}(C) \\
\notag j_{i,[k]}(\bar{x}) &=& j_{i-1,[k]}(x) \\
\label{eq:B-shift} \bar{B}_{j_{i,[k]}} &=& B_{j_{i-1,[k]}} \\
\label{eq:r-shift} j_{i,[k]}(\bar{r}) &=& j_{i-1,[k]}(r).
\end{eqnarray} 
With a proper computer 
implementation those assignment
operations could be avoided completely, 
for example by implementing some form of 
pointers swap or just by designing
the data structures to be easily extended
to accommodate new data. This last approach is 
implemented in our current source 
code so that the computational complexity is 
negligible. 

What is left to be computed are two parts: $j_{1}(\bar{X})$
and $z(\bar{X})$.

\subsubsection{The $\Phi_{j_1}$ part} 

From \eqref{eq:Phi-all} we have 
\begin{equation*}
\Phi_{j_{1,[k]}}(a(x)) = \left(y\left(u(x)\right)\right)_{[k]} = \left(z(x), w_{n+1} * F^{[n]}\left(u(x)\right) \right)_{[k]}.
\end{equation*}
It is obvious that $\Phi_{j_{1,[k]}}$ 
as a function of the variables $a$ is in fact a function only of 
the subset of variables $u$, so is the function
\begin{equation*}
\Phi_{j_1} = \left(\Phi_{j_{1,[0]}}, \Phi_{j_{1,[1]}}, \ldots, \Phi_{j_{1,[n(f, \eta)]}}\right)
\end{equation*} 
Therefore, with some abuse of notation, we can define $D_u \Phi_{j_1}(u)$
for all $u \in u(X)$. This is a matrix $\matrices(K, dim(u))$ with 
$K = (1+n(f, \eta)) \cdot d$ and is given by:
\begin{equation*}
D_u\Phi_{j_1} = \left(
\begin{array}{c}
D_u \Phi_{j_{1,[0]}} \\
 D_u\Phi_{j_{1,[1]}} \\
  \vdots \\ 
  D_u\Phi_{j_{1,[n(f, \eta)]}}
\end{array}
\right).
\end{equation*} 
This way we can skip the computation of 
many entries in $D\Phi_{j_{1,[k]}}$.

Applying the trick from Eq.~\eqref{eq:reduced-multiplication}
on each row of $D \Phi_{j_1} \cdot C$ we get:
\begin{equation}
\label{eq:Du-apply-trick}
D \Phi_{j_1}(u) \cdot C = D_u \Phi_{j_1}(u) \cdot u(C).
\end{equation}
The dimension of matrices taking part in the multiplication
on the right side are $\matrices(K, dim(u))$ and 
$\matrices(dim(u), N)$. Therefore, the cost of the computation
is $O(K \cdot dim(u) \cdot N)$ instead of 
$O(K \cdot M \cdot N)$ when performing the 
multiplication on the left.
 
\subsubsection{The $\Phi_{z}$ part} 

Similarly as before, we treat $\Phi_z(a(x))$
as a function of only used variables $\Phi_z(u(x))$.
It has an explicit formula:
\begin{equation*}
\Phi_z(u(x)) = \taylor(y(u(x)); h) = \sum_{k=0}^{n(f, \eta)} \Phi_{j_{1,[k]}}(u(x)) \cdot h^k,
\end{equation*}
so that the Jacobian $D_u \Phi_z(u)$ has a known form
expressed in terms that are already computed:
\begin{equation}
\label{eq:jacobian-sum}
D_u \Phi_z(u) = \sum_{k=0}^{n(f, \eta)} D_u\Phi_{j_{1,[k]}}(u(x)) \cdot h^k.
\end{equation}
Therefore, the computation of $D_u \Phi_z$ is computationally
inexpensive in comparison to the matrix-matrix multiplication.

Applying again the trick from Eq.~\eqref{eq:reduced-multiplication}
on each row of $D \Phi_{z} \cdot C$ and the Eq.~\eqref{eq:jacobian-sum}
we get:
\begin{eqnarray*}
D \Phi_{z}(u) \cdot C &=& \sum_{k=0}^{n(f, \eta)} \left(D_u \Phi_{j_{1,[k]}}(u) \cdot u(C)\right) \cdot h^k \ =  \\
& = & \sum_{k=0}^{n(f, \eta)} \left(j_{1,[k]}\left(D_u \Phi_{j_1}(u) \cdot u(C)\right)\right) \cdot h^k,
\end{eqnarray*}
where the matrix term $D_u \Phi_{j_1}(u) \cdot u(C)$ is already computed in 
Eq.~\eqref{eq:Du-apply-trick}. Therefore the cost of this operation consists 
only of additions and scalar-matrix multiplication, thus may be neglected
in comparison to other operations in the Lohner-type algorithm. 

\subsubsection{Summary computational cost of $D\Phi \cdot C$ multiplication}

Taking all into account, we get that computing the matrix-matrix multiplication
of the Lohner-type algorithm applied to the integration scheme for DDEs is
dominated only by the multiplication $D_u \Phi_{j_1}(u) \cdot u(C)$
in \eqref{eq:Du-apply-trick}. Its cost is $O(K \cdot dim(u) \cdot N)$
with $K = (1+n(f, \eta)) \cdot d$ which is a big reduction
from $O(M^2 \cdot N)$. To better see this, note that $dim(u) \le M$ and assuming $\eta_i = n$ 
for all $i$ (for simplicity) we have $M = d \cdot (1 + (n+1) \cdot p) = O(d\cdot n\cdot p)$. 
In that case $K = d \cdot (n+2)$ and we get the upper estimate on the complexity 
$O(d \cdot n \cdot M \cdot N) = O(d^2 \cdot n^2 \cdot p \cdot N)$.
The naive implementation has the complexity $O(d^2 \cdot n^2 \cdot p^2 \cdot N)$.
Moreover, if we assume a constant and small number of delays 
used in the definition of r.h.s. $f$ of equation \eqref{eq:dde-many-delays},
$m = const << p$ then we get complexity of order $O(d^2 \cdot n^2 \cdot N)$
- a reduction of the factor $p^2$. Noting that $p$ is usually the biggest
of the parameters $d, p, n$ (see applications), we get an enormous
reduction in the computation times, making the algorithm feasible to be applied
for a variety of problems. To see how big is the reduction let assume $p = 128$,
$n=4$, $d = 1$ and $N = M$ as in the Mackey-Glass examples. We get 
$(d p n)^3 = 134217728$ of order $10^8$, 
whereas $d^2 \cdot n^2 \cdot M = 10256$ of order $10^4$.

\subsubsection{QR decomposition on $\Phi_z$, $\Phi_{j_1}$ parts in the case $d > 1$}

Up to now we only focused \eqref{eq:dominant-mul} in the Lohner-type
algorithm for DDEs. Now, we need to return to the problem of 
managing local errors - the part $B \cdot r$ in \eqref{eq:rect2-lohner-set-structure},
and the formulas in Eqs.~\eqref{eq:rect2-QR}-\eqref{eq:rect2-Phi-r}.

First we note that the structure of matrix $B$ in \eqref{eq:A-structure}
is block diagonal \eqref{eq:B-structure}. We want to create a matrix 
$\bar{B}$ in the representation of $\Phi(X) + \Rem(X)$ of the same structure. 
The choice of the block-diagonal structure of $B$ is dictated by
the need to compute $Q^{-1} = \bar{B}^{-1}$ in \eqref{eq:rect2-Phi-r}.
We cannot hope to be able to rigorously compute decomposition $Q \cdot R$ 
or rigorously invert a big and full matrix $\bar{B}$, as those operations 
are ill-conditioned and very costly ($O(M^3)$). The sparse diagonal matrix 
$B$ removes both those problems with the trade-off in a form of
more complicated algorithm and some extra error terms. 

We have already used the structure of $B$ in Eqs.~\eqref{eq:B-shift}-\eqref{eq:r-shift}
to reduce problem  complexity significantly for the shift part, i.e. 
computing $\bar{B}_{j_i,[k]}$ for all $i > 1$. What is left
to compute is $\bar{B}_{j_1,[k]} \in \matrices(d, d)$ and $\bar{B}_{z} \in \matrices(d, d)$. 

To define appropriate $Q_{j_{1,[k]}} = \bar{B}_{j_{1,[k]}}$ and 
a new $j_{i,[k]}\left(\bar{r}\right)$ we investigate 
the term $\left[D \Phi([X])\right] \cdot B \cdot r$ from Eq.~\eqref{eq:rect2-Phi-r}.
Taking projection onto the $j_{1,[k]}$-th coordinate 
and using the $u$-variable trick we get:
\begin{eqnarray*}
j_{1,[k]}\left( \left[D \Phi([X])\right] \cdot B \cdot r \right) &=& \left[D_u \Phi_{j_{1,[k]}}([X])\right] \cdot u(B) \cdot u(r) \ = \\
   & =: & D \cdot u(B) \cdot u(r),
\end{eqnarray*}
with $D = \left[D_u \Phi_{j_{1,[k]}}([X])\right]$ for a shorter notation. 
A close inspection reveals that $D \cdot u(B) \in \matrices(d, dim(u))$.
Such a matrix is not suitable to apply the QR strategy of the Lohner's set.
We expand further:
\begin{eqnarray*}
D \cdot u(B) = \left(\begin{array}{ccccc}
	D_z \cdot B_z &  
	D_{j_{p_1,[0]}} \cdot B_{j_{p_1,[0]}} & 
	D_{j_{p_1,[1]}} \cdot B_{j_{p_1,[1]}} &
	\cdots &
	D_{j_{p_m,[\eta_{p_m}]}} \cdot B_{j_{p_m,[\eta_{p_m}]}}
\end{array}\right),
\end{eqnarray*}
where $D_{j_{q,[s]}} =  \left[ D_{j_{q,[s]}} \Phi_{j_1,[k]}([X])\right] \in \matrices(d, d)$.
Now, the term $D \cdot u(B) \cdot u(r)$ can be computed as follows:
\begin{equation}
D \cdot u(B) \cdot u(r) = \left(D_z \cdot B_z\right) \cdot r_z + \sum_{j_{i,[k]} \in u} \left(D_{j_{i,[k]}} \cdot B_{j_{i,[k]}}\right) \cdot r_{j_{i,[k]}} 
\end{equation}
Now, a decision has to be made, as to which $I \in u = (z, j_{p_1,[0]}, j_{p_1,[1]}, \ldots)$
to choose for the QR decomposition:
\begin{equation*}
Q_{j_{1,[k]}} \cdot R_{j_{1,[k]}} = \midpoint(D_I \cdot B_I).
\end{equation*}
and to compute $\bar{r}_{j_{1,[k]}}$ according to \eqref{eq:rect2-Phi-r}:
\begin{eqnarray}
\bar{r}_{j_{1,[k]}} 
\label{eq:fast-phi-r-QR}		& = & \left(Q_{j_{1,[k]}}^{-1} \cdot D_z \cdot B_z\right) \cdot r_z + \sum_{j_{i,[k]} \in u} \left(Q_{j_{i,[k]}}^{-1} \cdot D_{j_{1,[k]}} \cdot B_{j_{i,[k]}}\right) \cdot r_{j_{i,[k]}} \ + \\
\label{eq:fast-phi-r-S}		& + & \left(Q_{j_{1,[k]}}^{-1} \cdot j_{1,[k]}\left(S - \midpoint(S)\right)\right) \cdot r_0 \ + \\
\label{eq:fast-phi-r-PhiRem}	& + & Q_{j_{1,[k]}}^{-1} \cdot \left( \Phi_{j_{1,[k]}}(x) + \Rem_{j_{1,[k]}}(x) -\midpoint\left(\Phi_{j_{1,[k]}}(x) + \Rem_{j_{1,[k]}}(x)  \right) \right)
\end{eqnarray}
The matrix-matrix and matrix-vector operations are 
done in the order defined by parentheses.
Note, that all operations are well defined. The dimensions
of the matrices are as follows: $Q_{j_{1,[k]}}^{-1}, D_J, B_J \in \matrices(d,d)$, and 
$j_{i,[k]}(S) \in \matrices(d, N)$ for all variables $J \in u$. The vectors are:  
$r_0 \in \I^N$ while $r_J, \Phi_{j_{1,[k]}}(x), \Rem_{j_{1,[k]}}(x) \in \I^d$.

\begin{rem}
The same algorithm might be used to compute $\bar{r}_{z}$. 
We only change the projection $j_{1,[k]}$ to $z$ in the presented
formulas. 
\end{rem}

\subsubsection{Complexity of handling doubleton set structure}

The computational cost of using QR strategy with the doubleton 
set structure \eqref{eq:rect2-lohner-set-structure} in comparison
to the interval form of the error terms in the basic structure
\eqref{eq:basic-lohner-set-structure} is as follows. In the basic 
set structure the operation \eqref{eq:Phi-r} is exactly realized
by the presented algorithm when we take $Q = Id_{d \times d}$
and we use fact that each $B_J = Id$.
We can of course skip multiplication by $Id$. Therefore,
the cost of operations is (we count scalar multiplications):
\begin{itemize}
\item for \eqref{eq:fast-phi-r-QR}: $\frac{dim(u)}{d} \cdot d^2 = d \cdot \dim(u) $
\item for \eqref{eq:fast-phi-r-S}: $d \cdot N$,
\item for \eqref{eq:fast-phi-r-PhiRem}: $0$ (no matrix-matrix and matrix-vector multiplications).
\end{itemize}
In total, we get that computing $\bar{r}_J$ for each 
$J \in \left\{z, j_{1,[0]}, \ldots, j_{1,[n]}\right\}$ is 
$O\left(d \cdot (dim(u) + N)\right)$.
Taking into account $dim(u) \le M$ and $d << M$, together with the assumption 
$N = M$ (in applications) we get the complexity $O(d \cdot M) = O(M)$
under assumption that $d = const$, small. 

The algorithm for the doubleton set \eqref{eq:rect2-lohner-set-structure}
with non-trivial QR-decomposition has the following 
complexity for each $\bar{r}_J$:
\begin{itemize}
\item cost of computing QR decomposition for a matrix $D_I \cdot B_I \in \matrices(d,d)$, usually it is $O(d^3)$ multiplications,
\item cost of computing $Q^{-1}$, should not exceed $O(d^3)$, but it is usually $O(1)$ - if $Q$ is chosen to be orthogonal,
\item for \eqref{eq:fast-phi-r-QR}: $\frac{dim(u)}{d} \cdot \left(d^3 + d^3 + d^2\right) = O(d^2 \cdot \dim(u))$
\item for \eqref{eq:fast-phi-r-S}: $d^2 \cdot N + N \cdot d = O(d^2 \cdot N)$,
\item for \eqref{eq:fast-phi-r-PhiRem}: $d^2$.
\end{itemize}
In total, the complexity is $O\left(d^2 \cdot (d + dim(u) + N)\right)$.
Under the same assumptions as before, we estimate that in applications
the complexity is $O(d^2 \cdot M)$. Therefore, handling
the proposed doubleton structure is not much more costly
than using the interval form of the remainder (at least for small $d$). 

\begin{rem}
Please note, that the current strategy does not help
in the scalar case $d = 1$. The two sets are in this 
case equivalent and the computational cost is basically
the same. 
\end{rem}

\subsubsection{Choice of the matrix $D_I \cdot B_I$ for the QR procedure}

As to the selection of the matrix $D_I \cdot B_I$ used
in the QR decomposition procedure, 
in our current implementation we always use $I = z$. 
The motivation is as follows: in our applications
to the R\"ossler system we apply the method
to a perturbed system $x'(t) = f(x(t)) + \epsi \cdot g(x(t-\tau))$
with $\epsi$ small. Therefore, we expect that the influence
$D_J \cdot B_J$ for $J \ne z$ will be small. Taking $I = z$
allows to compare the method to the ODE version of 
the proofs. Indeed, if we set $\epsi = 0$ and integrate
the problem with our code, all $D_J = 0_{d \times d}$
and the method \eqref{eq:fast-phi-r-QR}-\eqref{eq:fast-phi-r-PhiRem}
reduces to that of the ODE (we only do operations on $z$ coordinate). 

Other choices of $I$ are easily implementable and one might
want to pursue other forms of the matrix $B$, for example
using various size blocks $B_{j_i}$ that does QR decomposition
on more than $d$ dimensions.

%% file: sections/7_appendix_B_data.tex

\section{\label{app:data}Description of the data and computer programs}

In this appendix we present details of the methodology to generate
initial sets for computer assisted proofs. As the data sets are large, 
one cannot hope to select good initial set candidates ,,by hand'', 
as can be done sometimes in the context of low-dimensional ODEs or maps. 
Instead, some kind of automatic or semi-automatic procedure must be used. 

\subsection{Source codes and a virtual machine}
The compressed archive of the source codes can be downloaded
from the web page \cite{www-praca-DDE-2-sources}. A file \Verb#README.txt#
from the main directory contains information on dependencies, compiling process 
and running the programs on the user's own computer. For a users
not wanting to compile files by themselves, we made an image of a 
virtual machine (VM) with Linux system, all compiler tools, and the source
codes compiled to executables that were used to produce data for this 
paper. It can be downloaded from \cite{www-praca-DDE-2-vm}, where
one also find the instructions for running the virtual machine. 
A computer running \Verb#Docker# VM on a Linux system is needed
to use the virtual machine image. We recommend following
instructions on the official web page: 
\Verb#https://docs.docker.com/engine/install/# and
selecting the user's system. 

\subsection{List of programs used in computer-assisted proofs}

All the programs used in proofs reside in the subdirectories 
placed under the root directory of the compressed archive
or in the directory \Verb#/home/user/DDEs# in the VM image. 
This root directory is common for all programs, and in what follows
we give paths relative to this root directory. The programs
can be found in the \Verb#./examples# subdirectory. 
The data for the proofs used in this paper can be found
in \Verb#./results# subdirectory. 

\subsubsection{The program used to produce data in Table~\ref{tab:benchmark}}

The programs used for benchmark in  Table~\ref{tab:benchmark}
can be found in \Verb#./examples/benchmark#. The program can be compiled
by issuing the following command in the main directory of the
source codes:

\vspace{0.5em}
\Verb#make benchmark#
\vspace{0.5em}

To obtain data from table~\ref{tab:benchmark} one 
needs to invoke the following command in the \Verb#./bin/examples/benchmark#
directory:

\vspace{0.5em}
\Verb#./benchmark diam='[0,0]' xi='[0,0]' dirpath='table-1'#
\vspace{0.5em}

Results of the computations will be stored in \Verb#./bin/examples/benchmark/table-1# and
will consist of several files (\Verb#.tex#, \Verb#.pdf#, \Verb#.png#, \Verb#.dat#, etc.).
One need \Verb#latex# and \Verb#gnuplot# packages installed 
in the host system for the program to work correctly. 
One can make various tests by changing the parameters.
The full list of parameters is presented below:

\vspace{0.5em}
\Verb#./benchmark \ # \newline
\Verb#	initial='{[1.1,1.1]}' \ # \newline
\Verb#	dirpath='.' \ # \newline
\Verb#	prefix='benchmark' \ # \newline
\Verb#	N='[8,8]' \ #\newline
\Verb#	n='4' \ # \newline
\Verb#	p='128' \ #  \newline
\Verb#	epsi='50\%' \ # \newline
\Verb#	diam='[-1e-06,1e-06]' \ # \newline
\Verb#	xi='[-0.1,0.1]' # \newline
\vspace{0.5em}

The values on the right of \Verb#=# are the default values. 
It is important to put the parameters into single quotes \Verb#'...'#.
The system modelled by the program is Mackey-Glass equation with $\gamma = 1$,
$\beta = 2$, $\tau = 2$ and $n = $ \Verb#N# (do not confuse with $n$ in the definition
of $C^n_p$. Parameter \Verb#initial# represents $x_0$ in the definition of the initial 
Lohner set $X(A_0, R_0)$, $A_0 = x_0 + Id \cdot r_0$, $r_0$ given by \Verb#diam# (see later). 
Parameter \Verb#initial# can either be an interval, in that case the program 
will treat $x_0$ as a representation of a constant initial function $x_0 \equiv $ \Verb#initial#; 
or it can be a path to a file containing a vector describing the $(z, j)$ part of the initial segment. For examples
of such files, please refer to initial data in the computer assisted proofs. 
Parameter \Verb#dirpath# describes the directory of the output files from the program.
It is advised to select non-existing folder, as the program overwrites existing 
files without asking. Parameter \Verb#prefix# will be appended in front of all 
filenames. 
Parameters \Verb#n#,  \Verb#p# corresponds to $n$, $p$ in $C^n_p$, and are the order
of the representation and the number of grid points on the base interval $[-2, 0]$, respectively.
The full step $h = \frac{2}{p}$.
Parameter \Verb#epsi# corresponds to $\epsi$ step done to simulate computation
of the Poincar\'e map, and can be given as a percentage of the full step $h$
or as an explicit interval. Parameter \Verb#diam# is the diameter $r_0$ of the $A_0 = x_0 + Id \cdot r_0$ 
part in $X(A_0, R_0)$, while \Verb#xi# is the diameter of $R_0$.

\subsubsection{The programs used in the proof of Theorem~\ref{thm:rossler-delayed}}
The programs used in the proof of Theorem~\ref{thm:rossler-delayed}
reside in the following subdirectories:
\begin{itemize}
\item 
\Verb#./examples/rossler_delay_zero# for the R\"ossler original 
ODE~\eqref{eq:rossler}, but studied in the extended space $C^n_p$ over 
the base delay interval $[-1,0]$. Programs in this directory are used 
to generate common set for all the proofs.  

\item 
\Verb#./examples/rossler_delay_rossler# for the delayed perturbation 
of the form $g = f$.

\item 
\Verb#./examples/rossler_delay_other# for the delayed perturbation 
$g$ given in \eqref{eq:perturbation-sin}.
\end{itemize}

Each of the directories contain the following programs
\begin{itemize}
\item \Verb#nonrig_attractor# (non-rigorous, approximate), it is used 
to generate initial set of functions to be used in the program 
\Verb#nonrig_coords#. It generates also plot to view the structure 
of the apparent attractor.  

\item \Verb#nonrig_coords# (non-rigorous, approximate) it computes first
approximation of the coordinate frame for the set $A$ in the definition of the
set $X(A, R)$. The program uses the output of the program
\Verb#nonrig_attractor# to generate a set $S$ of several hundred solution
segments lying on the section $S_0$. Due to the nature of R\"ossler attractor,
those segments are contained in a thin strip over $(x, y)$-plane (set is thin in
$z$ direction), see Fig.~\ref{fig:attractor-rossler-delayed}. Then, a reference
solution $v_{ref}$ segment is selected as the one closest to the centre of this
collection. Let denote by $|S|$ the number of solution segments in $S$. We
define:
\begin{eqnarray*}
w^i & = & v^i - v_{ref}, \quad v^i \in S \setminus \{ v_{ref} \}, \\
u^i & = & \frac{u^i}{u^i_2}, \\
u & = & \frac{\sum u^i}{|S| - 1}     
\end{eqnarray*}
or, in other words, $u$ is the mean vector that spans
the intersection of the apparent R\"ossler attractor
with the section $S_0$. The vector is normed in such a
way that $u_2 = 1$. This corresponds to $\pi_y u = 1$.
The initial coordinate frame is chosen to be:
\begin{equation}
\tilde{C} = 
\left(
\begin{array}{ccccc}
1			&		& 0 & \cdots 		& 0  	\\
0			& 		& 0 & \cdots 		& 0  	\\
\vdots 	& u^T	&   &     	  	& 		\\
\vdots 	&		&   &	Id_{(M-2) \times (M - 2)}	& 		\\
0 			&		&   &        		& 
\end{array}
\right),
\end{equation}
that is, first column corresponds to the normal vector to the section 
hyperplane $S_0 = \{ v : \pi_x v(0) = 0\}$, the second column corresponds 
to the nominally unstable direction $u$, and the rest of coordinates are 
just the canonical basis in $\R^{M-2}$. The set $X(A, \Xi)$ is then defined 
with: $A = v_{ref} + \tilde{C} \cdot r_0$, with $r_0$ to be defined 
by the next program \Verb#rig_find_trapping_region#.

\item \Verb#rig_find_trapping_region# (rigorous) as input, it takes the width 
$W = [W_l, W_r]$ of the set in the nominally unstable direction $u$ and the
coordinate frame $v_{ref}$ and $\tilde{C}$. Then it starts with 
$A^0 = v_{ref} + \tilde{C} \cdot r$ with $r_1 = 0$, $r_2 = W$, 
$r_i = (-\epsilon, \epsilon)$ for some small $\epsilon$ for $i > 2$. 
Then it tries to obtain the set $X(A, \Xi)$ iteratively ,,from below'', i.e. at
each step $k = 0,1,...$ it computes image $P(X(A^k, \Xi))$ and checks if it is
subset of $X(A^k, \Xi)$. If the test is passed, the program stops, otherwise it
takes $A^{k+1} = \hull(A^{k}, \pi_{A} P(X(A^k, \Xi)))$ and continues to the next
step. The computation of $P(X(A^k, \Xi))$ is done by dividing the input set into
$N$ pieces along the nominally unstable direction $r_2$. For proofs used in this
work $N = 200$.

Finally, when the set $A$ is found, the coordinate frame $\tilde{C}$
is changed (by rescaling) to $C = \tilde{C} \cdot \mathrm{Diag}(r)$ so that 
$A = v_{ref} + C \cdot \left(\{ 0 \} \times W \times \Ball^{\|\cdot\|_\infty}_{M-2}(0,1)\right)$.
The matrix $\mathrm{Diag}(r)$ denotes the diagonal matrix with $r_i$'s 
on the diagonal.

\item \Verb#rig_prove_trapping_region_exists# (rigorous) the program computes
the image of the set $X(A, \Xi)$ under the Poincar\'e map $P : S_0 \to S_0$.
If the \Verb#rig_find_trapping_region# is successful in finding 
the rigorous candidate $A$, then this program must succeed, as it
computes the image $P(X(A, \Xi))$ in the same way, dividing
the set into the same $N$ pieces. 

\item \Verb#rig_prove_covering_relations# (rigorous) the program checks
the conditions (CC2A) and (CC2B) of Lemma~\ref{lem:one-unstable} 
on sets $X(N_1,\Xi)$ and $X(N_2, \Xi)$.
The sets are defined as the restrictions of the set $X(A, \Xi)$ 
on the nominally unstable direction $r_2$. The user
can manipulate the definitions of sets changing the values
in the configuration file \Verb#rig_common.h#. 
\end{itemize}

\begin{rem}
To prepare the set $X(A, \Xi)$ for the Theorem~\ref{thm:rossler-delayed},
we run \Verb#rig_find_trapping_region# for the system
without the delay first. Then we use it again on the resulting set
for the delayed systems. In this way we obtain three sets
$A_0$ for the unperturbed system, $A_f$ for $g=f$ and $A_g$ 
for the system with $g$ given as in \eqref{eq:perturbation-sin}.
As the final set $A$ we take $A = \hull(A_0, A_f, A_g)$. The computer
assisted proofs show that this set is a trapping region for
all systems.  
\end{rem}

The data used in the proofs can be found for each of the systems
in the respective directories:
\begin{itemize}
\item 
\Verb#./results/work3_rossler_delay_zero#,

\item 
\Verb#./results/work3_rossler_delay_rossler#,

\item 
\Verb#./results/work3_rossler_delay_other#.
\end{itemize}

\subsubsection{The programs used in the proof periodic orbits in the Mackey-Glass equation}

The programs used in the proof of Theorem~\ref{thm:mg-periodic}
can be found in \Verb#./examples/mackey_glass_unstable#. 
There is a single program for each of the orbits: 
\Verb#prove_T1#, \Verb#prove_T2#, and \Verb#prove_T4#, respectively.
The data for the proofs is stored in \Verb#./results/work3_mackey_glass_proofs#.
Additionally, a set of generic programs to generate the flow, coordinates
and sets for proofs are available in \Verb#./examples/mackey_glass_finder#.
The programs comprise a full set of tools that once compiled
can be used to find candidates for periodic solutions to Mackey-Glass
equation for any set of parameters and later to prove their existence. 
With a little effort the programs can be adjusted to work
with any scalar DDE - small changes in the file
\Verb#constans.h# should suffice. 
The collection of programs is as follows:
\begin{itemize}
\item 
	\Verb#attractor-coords-nonrig# non-rigorous program to generate good 
	coordinates for presentation of the Mackey-Gass attractor (used
	for parameters in chaotic regime).
\item 
	\Verb#compare-rig# a program to compare two interval sets and
	print the comparison in a human-friendly manner. 
\item 
	\Verb#draw-nonrig# non-rigorous program to draw solutions
	in various ways. 
\item 
	\Verb#find-nonrig# a non-rigorous Newton-like method
	to refine the non-rigorous candidate for a periodic solution
	with high accuracy.
\item 
	\Verb#jac-poincare-nonrig# a program to compute non-rigorous
	image of the Poincar\'e map $P$ for a single initial condition $x$, 
	together with the (approximate) Jacobian $DP(x)$.
\item 
	\Verb#periodic-coords-nonrig# a program to generate the ,,good''
	basis for the set $X(A,R)$. The procedure is described
	in more details later in this Appendix.
\item 
	\Verb#poincare-nonrig# a program to compute non-rigorous
	Poincar\'e map for a collection of initial conditions. 
\item 
	\Verb#poincare-rig# a program to compute rigorously
	the image of the Poincar\'e map $P$ on a (p,n)-fset 
	$X(A, R) \subset C^n_p$.
\item 
	\Verb#simple-coords-nonrig# alternative program to generate
	very simple coordinates, where majority of base vectors
	comes from the cannonical basis. 
\end{itemize}

An exemplary process of finding good candidates and initial sets for 
orbits $T1$, $T2$ and $T4$  
can be found in the subdirectory \Verb#./results/work3_mackey_glass_finder#. 
The scripts contained there was used to generate data for this paper
 and they are as follows (\Verb#$i# $\in \{1, 2, 4\}$):
\begin{itemize}
\item 
	\Verb#setup.sh# a common setup script for all other scripts.
	In principle, user only edit this file. 
\item 
	\Verb#T$i-1-make-coords.sh# uses some aforementioned programs to
	refine the candidate solution $x_0$ (mid point of the set) 
	and generate good coordinates for the fset $X(A,\Xi)$,
	with $A = x_0 + C \cdot r_0$.
	
	In case of the solution $T4$, the script is more complicated
	as it finds a pair of sets $X_1 = X(A_1,\Xi_1)$, $X_2 = X(A_2,\Xi_2)$
	for the covering $X_1 \cover{P_1} X_2 \cover{P_2} X_1$, where $P_1, P_2$
	are two Poncar\'e maps, defined between two different sections.
\item 
	\Verb#T$i-2-find-start.sh# try to do the first step of an algorithm to find
	appropriate $r_0$	and $\Xi$ such that $X(A, \Xi) \cover{P} X(A, \Xi)$. 
	It starts with a thin set $x_0 + C \cdot \{0\} \times W \times [0,0]^M-2$ 
	and build the set $X_1 = P(X)$. The sets are set to have the first coordinate
	of width $W$ (guessed by the user, similarly to the building sets in 
	the R\"ossler case).
	
	The procedure is more involved in the case of $T4$, as there are two
	sets $X_1$ and $X_2$ on different sections such that we hope 
	$X_1 \cover{P_1} X_2 \cover{P_2} X_1$. The program
	tries to find both sets at the same time. 
	
\item 
	\Verb#T$i-3-find-once.sh# subsequent iteration of the 
	previous algorithm. The user needs to run this script until
	satisfactory $r_0$ and $\Xi$ is found. This part of the finding
	procedure is semi-automatic, as the decision when to stop the procedure
	is left to the user. 
	
	The data generated by the authors of this manuscript
	can be found in the following directory: \Verb#./data/examples/mg#.
\item 
	\Verb#T$i-4-proof.sh# the script runs the final check
	of the covering relations in the proof of Theorem~\ref{thm:mg-periodic}.
	If the $W$ width in the unstable direction was set improperly
	in the previous steps (too narrow), then the script might fail
	to prove the covering relation. 
	
	In case of $T1$ and $T2$ the program checks the simplest
	covering relation $X \cover{P} X$.
	In case of $T4$ the covering relation to check is more involved: 
	$X_1 \cover{P_1} X_2 \cover{P_2} X_1$, where $P_1, P_2$.	
\end{itemize}

\subsubsection{Idea of the coordinate selection}
The following procedure is adopted in the program \Verb#periodic-coords-nonrig#
for choosing right coordinates for periodic orbits proofs.
Each set $V_i = X(A_i, \Xi_i)$ with $A_i = \bar{T}^i_0 + C_i \cdot r_i$.
The matrix $C_i$ must be computed carefully, as it was shown in \cite{nasza-praca-focm}
and in more details in \cite{recent-advences-poincare}. In short, one expects
that the linearized dynamics near the stationary point $\bar{T}^i_0$ of the 
Poincar\'e map $P^{i}$, $i \in 1, 2, 4$ should decompose into invariant 
subspaces $E_c \oplus E_u \oplus E_s$ with $dim(E_c) = dim(E_u) = 1$. The 
subspace $E_c$ corresponds to the direction along the flow, where $E_u$ 
is the unstable space of solutions that have a backward in-time limit 
at the fixed point, where $E_s$ is the space of those solutions 
that approach the fixed point as $t \to +\infty$. The matrix $C_i$ 
is chosen as a composition of the bases of appropriate subspaces:
\begin{equation}
C_i = 
\left(
\begin{array}{ccccc}
c^T & u^T & s_3^T & \ldots & s_M^T
\end{array}
\right),
\end{equation}
where we have:
\begin{itemize}
\item $c$ is the vector defining a section and is chosen in the program
as a left eigenvector (i.e. eigenvector of the transposed matrix) 
of the approximate matrix $\overline{DP}(\bar{T}^i_0)$
corresponding to the eigenvalue $\lambda_2 = 1$. 
\item $u$ is the eigenvector of $\overline{DP}(\bar{T}^i_0)$
corresponding to the largest and unstable eigenvalue $\lambda_1$ with $|\lambda_1| > 1$.
\item the set of vectors $s_j$ is the basis of the (finite projection) of the stable
subspace $E_s$. It is obtained as an orthonormal basis orthogonal to the vector
$\tilde{u}$ - the left eigenvector corresponding to the unstable eigenvalue 
$\lambda_1$. 
\end{itemize}
The reason why those are chosen as described is explained
in details in \cite{nasza-praca-focm,recent-advences-poincare}
and we skip the details here. We only hint that the selection
of $c$ guarantees to have a very thin interval $[t_p(X)]$
in rigorous computations of Poincar\'e maps, where
selection of $s_3, \ldots s_M$ gives a hope that 
the finite dimensional projection of the $P(X)$ onto stable
subspace $E_s$ could be mapped inside the stable part
of the initial set $X$ ($\pi_{E_s} P(X) \subset \pi_{E_s}X$)
without the need to resort to a set subdivision in computations.

\subsubsection{Utility programs}

The programs in directory \Verb#./examples/converter# are
used in some other scripts to do conversion between 
different kinds of data, for example making interval
versions of vector/matrices from their \Verb#double# counterparts.
One important program is the matrix converter \Verb#convmatrix# that can
compute rigorously in high precision the rigorous inverse
of interval matrix. For more information how to use those
programs, see their source code documentation. 
The list of utility programs is as follows:
\begin{itemize}
\item 
	\Verb#convmatrix# a conversion between various formats of matrices. 
	It can compute rigorous inverse of a matrix in high precision. 
\item 
	\Verb#convvector# a conversion between various formats of vectors. 
\item 
	\Verb#growvector# converts an interval vector $[x]$ into 
	$[w] = \midpoint([x]) + r \cdot \left([x] - \midpoint([x])\right)$, 
	for some $r \in \R$.
\item 
	\Verb#splitvector# converts an interval vector $[x]$ into 
	$[w] = \midpoint([x])$ and $[v] = [x] - \midpoint([x])$.
\item 
	\Verb#crmatrix# takes a matrix $C$ and a vector $r$ and makes a rigorous 
	matrix $C_r$ such that $C \cdot r \subset C_r \cdot [-1,1]^M$.
\item 
	\Verb#rmatrix# takes a vector $r$ and makes a rigorous 
	matrices $R$ and $R^{-1}$ such that $R_{i,i} = r_i$, $R_{i,j}=0$ otherwise.	
\item 
	\Verb#invmatrixtest# performs a test if the pair of matrices $A$
	and $B$ given on input has the property that $Id \subset A \cdot B$ 
	and $Id \subset B \cdot A$. Computes also the width (maximal diameter
	of all entries) of $A \cdot B$ and $B \cdot A$, allowing to assess
	the quality of computed inverses.
\item 
	\Verb#matrixcmp# compares two (large) matrices in a human-friendly manner.
\item 
	\Verb#midmatrix# compute $\midpoint(A)$ for a matrix $A$.
\item 
	\Verb#vectorcmp# compares two (large) vectors in  a human-friendly manner.
\item 
	\Verb#vectorhull# compute an interval hull of all vectors given on input.
\end{itemize}